\documentclass[11pt,amsfonts]{article}
\usepackage{amssymb,epsfig}
\def\cz{{\Bbb C}}
\def\rz{{\Bbb R}}

\def\nz{{\Bbb N}}
\def\gz{{\Bbb Z}}
\def\pr{{\rm I \! P}}
\def\qz{{\Bbb Q}}
\newcommand{\skp}{\hspace{1pt}}
\newcommand{\eps}{\varepsilon}
\newcommand{\ph}{\varphi}

\newcommand{\id}{{\mbox{\rm Id}\skp}}
\newcommand{\st}{{\mbox{\rm\scriptsize st}}}
\newcommand{\sing}{{\mbox{\rm\scriptsize sing}}}
\newcommand{\reg}{{\mbox{\rm\scriptsize reg}}}

\newcommand{\can}{{\mbox{\rm\scriptsize can}}}
\newcommand{\hol}{{\mbox{\rm\scriptsize hol}}}

\newcommand{\rk}{{\mathop{\rm rk}\skp}}
\newcommand{\codim}{{\mbox{\rm codim}\skp\skp}}
\newcommand{\coker}{{\mbox{\rm coker}\skp}}
\newcommand{\prj}{{\mbox{\rm pr}\skp}}
\newcommand{\img}{{\mbox{\rm im}\skp}}

\newcommand{\ind}{{\mbox{\rm ind}\skp}}

\newcommand{\GL}{{\mbox{\rm GL}\skp}}
\newcommand{\Aut}{{\mathop{\rm Aut}\skp}}
\newcommand{\Hom}{{\mbox{\rm Hom}\skp}}
\newcommand{\Ext}{{\mbox{\rm Ext}\skp}}

\newcommand{\Spec}{{\mbox{\rm Spec}\skp}}

\newcommand{\End}{{\mbox{\rm End}\skp\skp}}

\newcommand{\Or}{{\mbox{\rm Or}\skp}}
\newcommand{\supp}{{\mbox{\rm supp}\skp\skp}}
\newcommand{\Vol}{{\mbox{\rm vol}\skp\skp}}
\newcommand{\LUS}{{\mbox{$\call$\rm US}\skp\skp}}

\def\GW{{\calg\calw}}

\newcommand{\diam}{{\mbox{\rm diam}\skp}}

\newcommand{\llongrightarrow}{{\begin{picture}(25,6)(2.5,-4)
  \unitlength 1pt\put(0,0){\vector(1,0){30}}\end{picture}}}
\newcommand{\llongleftarrow}{{\begin{picture}(25,6)(2.5,-4)
  \unitlength 1pt\put(29,0){\vector(-1,0){30}}\end{picture}}}
\newcommand{\longdownarrow}{{\begin{picture}(0,0)
  \unitlength 1pt\put(0,16){\vector(0,-1){30}}\end{picture}}}
\newcommand{\longuparrow}{{\begin{picture}(0,0)
  \unitlength 1pt\put(0,-14){\vector(0,1){30}}\end{picture}}}

\newcommand{\diagl}[1]%
  {\makebox[0cm]{${\scriptstyle#1\ }\longdownarrow
  \phantom{\scriptstyle#1\ }$}}
\newcommand{\diagr}[1]%
  {\makebox[0cm]{$\phantom{\ \scriptstyle#1}
  \longdownarrow{\ \scriptstyle#1}$}}
\newcommand{\updiagl}[1]%
  {\makebox[0cm]{${\scriptstyle#1\ }\longuparrow
  \phantom{\scriptstyle#1\ }$}}
\newcommand{\updiagr}[1]%
  {\makebox[0cm]{$\phantom{\ \scriptstyle#1}\longuparrow
  {\ \scriptstyle#1}$}}
\newcommand{\di}{\partial}
\newcommand{\dbar}{{\bar\di}}

\newcommand{\calc}{{\cal C}}

\newcommand{\calf}{{\cal F}}
\newcommand{\calg}{{\cal G}}
\newcommand{\calh}{{\cal H}}

\newcommand{\call}{{\cal L}}
\newcommand{\calm}{{\cal M}}

\newcommand{\calo}{{\cal O}}

\newcommand{\calq}{{\cal Q}}
\newcommand{\cals}{{\cal S}}
\newcommand{\calt}{{\cal T}}
\newcommand{\calu}{{\cal U}}
\newcommand{\calv}{{\cal V}}
\newcommand{\calw}{{\cal W}}

\newtheorem{prop}{Proposition}[section]
\newtheorem{theorem}[prop]{Theorem}
\newtheorem{lemma}[prop]{Lemma}
\newtheorem{corollary}[prop]{Corollary}
\newtheorem{rem}[prop]{Remark}

\newtheorem{defi}[prop]{Definition}
\newcommand{\pf}{{\em Proof. }}
\newcommand{\qed}{{{\hfill$\diamond$}\vspace{1.5ex}}}

\title{Gromov-Witten invariants for general symplectic manifolds}
\author{Bernd Siebert}
\date{December 7, $1997^2$}

\begin{document}

\maketitle

{\small \tableofcontents}
\vspace{1cm}


\addcontentsline{toc}{section}{Introduction}
\noindent
{\Large\bf Introduction}
\vspace{1.5ex}

\noindent
Invariants related to (pseudo) holomorphic curves on symplectic
manifolds with a compatible almost complex structure arose for the
first time in the celebrated work of Gromov \cite{gromov}. More refined
invariants, now known under the name of Gromov-Witten (GW) invariants,
turned up somewhat later in physics as correlation functions in
Witten's topological sigma models (\cite{witten0}, \cite{witten},
\cite{vafa}), and were put on a solid mathematical basis by Ruan and
Tian using a construction reminiscent of Donaldson theory \cite{ruan},
\cite{ruantian1}, \cite{ruantian2}, cf.\ also \cite{duffsal}. The mathematical
theory is based on transversality and genericity arguments with the
purpose of identifying certain moduli spaces of pseudo holomorphic
curves as manifolds, supplemented by the Gromov compactness
theorem and a dimension estimate on the ``bad part'' at infinity.  This
leads to quite satisfactory results on the most important classes of manifolds
in this context, notably of Fano and Calabi-Yau type. The restriction to
manifolds with certain positivity properties of the first Chern class is due
to the dimension estimate.

On the other hand, from the expected formal properties of the invariants
(axiomatized in \cite{kontmanin}) nothing hints on such limitations of
GW-theory. Deeper support in favour of general validity of GW-theory come
from the development of the analogous invariants in algebraic geometry, that
has been achieved recently by Li and Tian \cite{litian}, and (loosely based
on their ideas) by Behrend and Fantechi \cite{behrend}, \cite{behrendfantechi}.
Ultimately it would be desirable to have a similarly general theory in
the symplectic category, together with a comparison theorem linking to
the algebraic theory, in order to combine the computational power
and the strength of handling singular situations from algebraic
geometry with the full symplectic topological significance and range of
applications of the invariants.

The purpose of this paper is to solve the problem on the symplectic
side with a fully general theory of symplectic GW-invariants, without
appealing to genericity and  properly taking the ``bad locus'' at
infinity into account from the very formulation. The comparison with the
algebraic theory is deferred to a separate paper \cite{si}, cf.\ \cite{si2} for a
sketch.
\vspace{1ex}

\noindent
{\em Outline of paper.}\\
Our approach is based on the following simple topological fact: In
finite dimensions, given a section $s$ of an oriented topological
vector bundle $E$ over an oriented manifold, a Poincar\' e dual to the
Euler class of $E$ can be represented by a class with support in the
zero locus of $s$. The same works for differentiable Fredholm sections
of Banach bundles over Banach manifolds, as has been demonstrated
recently by Brussee \cite{brussee}. Now moduli spaces of pseudo
holomorphic curves on an almost complex manifold $(M,J)$ actually do
arise as zero loci of differentiable sections of Banach bundles over
mapping spaces. The basic idea is to extend this description over
the compactification, which inevitably means admitting changes of the
domain through the occurrence of bubbles.

A major difficulty is the possible change of topology of the domain. For
example, a sequence of holomorphic maps from the Riemann sphere to
a K\"ahler manifold $M$ can converge to a map from two copies of
the sphere to $M$. Explicitely let $M=\pr^2$ and $\ph_t: \pr^1
\rightarrow M$ map $[u,v]$ to $[u^2,uv,t\cdot v^2]$.
The image is the sequence of smooth conics $XZ-tY^2=0$
converging to the two lines $X=0$, $Z=0$ as $t$ tends to $0$.
To include such phenomena we look at spaces of maps
from possibly nodal Riemann surfaces ("prestable
algebraic curves") to $M$. As observed by Kontsevich, to make such
spaces Hausdorff one has to require the domain as small as possible.
That is, the map shall not factor through a map from a prestable curve
with fewer irreducible components (stability condition). We will in fact endow
spaces of stable continuous maps to any Hausdorff topological space
with a Hausdorff topology (Section~\ref{C(X)_topology}).

The definition of charts for Banach manifolds of stable maps
from prestable algebraic curves to $M$ has two ingredients.
First we study the deformation theory of prestable curves per se (Section~\ref
{section_deformation}). This is a generalization of the classical deformation
theory of Riemann surfaces and a special case of the well-established but
intricate deformation theory of compact complex spaces. The main point of
our discussion in this section is that this deformation can be effectively
decomposed into deformations of the irreducible components (classical
deformation theory of Riemann surfaces) and the obvious deformations
$zw=t$ of the nodes $zw=0$. The deformations fit together into a proper
holomorphic map of manifolds $q: \calc\rightarrow S$ with $q^{-1}(0) =C$.
So to any $s$ in the base $S$ of the deformation is associated a prestable curve
$C_s=q^{-1}(s)$. There is a Lipschitz retraction $\kappa:\calc\rightarrow C$
that restricted to fibers $C_s$ is a diffeomorphism away from the singular
points of $C$ and whose fibers over the nodes of $C$ are either points or
circles. Then given a map $\ph: C\rightarrow M$ and $s\in S$ we obtain a map
$\ph\circ\kappa_s: C_s\rightarrow M$ with the same image as $\ph$.

On the other hand, for fixed domains coordinate charts can be
modelled on appropriate spaces of sections of $\ph^*T_M$ via the
exponential map (having chosen a Riemannian metric on $M$): To
$v\in \Gamma(\ph^*T_M)$ one associates the map $z\mapsto
\exp_{\ph(z)} v(z)$. To combine this with the previous reasoning
one needs to identify spaces of sections for varying domains. The solution
to this problem is at the heart of our approach. Roughly, the identification
is simply by pull-back via $\kappa$ for $L^p$ spaces (where we introduce
certain weights at the nodes). Next, if $M$ is almost complex, $\ph^*T_M$
can be made into a holomorphic vector bundle by \cite{holizan}, \cite{ivshev}
(Section~\ref {hol_structure}). We show that the graph of the
corresponding $\dbar$-operator is in fact a Banach bundle
over the parameter space $S$ (Chapter~\ref{chapter_dbar}).
The technical tools are the inhomogeneous Cauchy integral formula,
a Calderon-Zygmund inequality in weighted $L^p$-spaces and
Banach bundles of bounded holomorphic \v Cech cochains.
Since we can only deal with one derivative by this process our
method is essentially restricted to maps of Sobolev class $L_1^p$.

A slight complication is due to the fact that if the domain $C$ has
non-trivial infinitesimal automorphisms the base $S$ does not effectively
parametrize deformations of $S$. In fact there is a local action of the
automorphism group of $C$ on $\calc$ respecting the fiber structure.
This shows that our candidates for charts have higher dimensional
subsets parametrizing the same maps. The induced action on $S$
typically has factors isomorphic to the $\cz^*$-action
on the unit disk, so a quotient would not be Hausdorff. Nevertheless, it
turns out that by stability the induced actions on the charts are free. In the algebraic
geometric situation a slice to this action is defined by intersection with
transversal hyperplanes in $M$, and the same works in spaces of continuously
differentiable maps. Unfortunately this is not true for $L_1^p$-spaces.
Instead we use integrals that compute centers of gravity of the pull-back
of bump functions on $M$ to achieve the same (Section~\ref{sect_rigidification}).

The parametrization of mapping spaces might still not be effective though,
because a possible finite group of biholomorphism of $C$ intertwining the map
acts on the charts. As we know from unobstructed situations in algebraic
geometry this can not be avoided and leads to the fact that our space of
maps will only have the structure of a Banach orbifold. This is usually merely
a matter of discomfort. However in our case the construction of the localized
Euler class requires the existence of a finite rank orbibundle with effective
actions of the local groups on its fibers. A general existence proof that I
gave in a previous version of this paper turned out to be essentially flawed
in discussions with Ruan. Instead we now mimick the construction of such
bundles in the algebraic geometric context (Section~\ref{section_F}). I
would like to point out that the existence of this bundle establishes the existence
of a ``virtual tangent bundle" for the moduli space of pseudo-holomorphic
curves, as an element of the $K$-group of orbibundles. Characteristic classes
of this bundle lead to new symplectic invariants.

The most subtle issue in the whole story, that caused a lot of delay
in finishing this paper, concerns regularity questions. The Banach orbifold
of maps is inherently {\em not} differentiable! This is already true for
spaces of maps with any kind of finite differentiability from a non-trivial
family of Riemann surfaces. The reason is that the action of the diffeomorphism
group on $L^p_k$- or $C^{k,\alpha}$-spaces is only continuous, not
differentiable. So different choices of local trivializations of the family
of Riemann surfaces lead to different differentiable structures. Basically
the same kind of regularity is true for our mapping spaces as well as for the
relevant objects for construction of the localized Euler class, the Banach
bundle and its section. Differentiability relative a finite dimensional base
space would still be enough to construct localized Euler classes. By the
appearance of slices in our charts the base space does not however have
a global meaning. Careful choices of charts and abstract perturbations
are necessary. This is formalized in the notion of ``Kuranishi section''
(Section~\ref{Kuranishi_structure}).
\vspace{1ex}

Our main results can be summarized as follows. For details we refer to
Theorems~\ref{virfc}, \ref{orbifold} and \ref{gw_are_invts}.
\begin{theorem}
  Let $(M,\omega)$ be a closed symplectic manifold with a tame almost
  complex structure $J$. Then the space $\calc(M;p)$ of stable
  parametrized marked complex curves in $M$ of Sobolev class $L_1^p$
  (Definition~\ref{cplx_curve}) is a Banach orbifold.
  Moreover, there is a Banach orbibundle $E$ over
  $\calc(M;p)$ with fiber $\check L^p(C;\ph^*T_M
  \otimes{\bar\Omega}_C)$ at $(C,{\bf x}= (x_1,\ldots
  ,x_k),\ph:C\rightarrow M)$ with an oriented Kuranishi section
  $s$ (Definition~\ref{Kuranishi_structure}) with $\hat
  s(C,{\bf x},\ph) =\dbar_J\ph$. The zero locus of $s$
  is the set $\calc^\hol(M,J)$ of stable pseudo holomorphic curves in $(M,J)$
  (Definition~\ref{def.ps.hol}), which is a locally finite dimensional Hausdorff space
  with compact components.
  
  Let $\calm_{g,k}$ be the moduli space of
  Deligne-Mumford stable $k$-marked algebraic curves of genus $g$, with the
  convention $\calm_{g,k}=\{{\rm pt}\}$ whenever $2g+k <3$.
  The localized Euler class $\GW_{g,k}^{M,J}\in H_*(\calc^\hol(M,J))$ associated to
  $(E,s)$ (Theorem~\ref{virfc}) gives rise to
  GW-correspondences (Definition~\ref{GW_correspondence})
  \[
    GW_{g,k}^{M,J}:H^*(M)^{\otimes k}\longrightarrow H_*(\calm_{g,k})\,,
  \]
  that are invariants of the symplectic deformation type of $(M,\omega)$.
  They coincide with the ones defined in \cite{ruantian2} in case $(M,\omega)$ is
  semi-positive.
\end{theorem}

The author believes that a definition of Floer homology for general symplectic
manifolds (needed to prove the Arnold conjecture on the numbers of fixed
points of neondegenerate Hamiltonian symplectomorphisms) should be
possible by similar methods. The ambient space in this case,
a space of stable connecting orbits, would be a Banach orbifold with corners.
\vspace{1ex}

While I was struggling with the mentioned regularity problems
there were announcements by Fukaya/Ono, Li/Tian, Hofer/Salamon
and Y.\ Ruan of similar results. The first to finish their paper
were Fukaya and Ono \cite{fukaya}, followed by Li and Tian \cite{litian2},
myself (dg-ga/9608005) and Ruan \cite{ruan2}). The paper of
Hofer and Salamon was not yet available at the time of this writing.
The closely related problem of generalizing Floer homology to general symplectic
manifolds is treated in \cite{fukaya}, \cite{liu} \cite{ruan2}. The methods in the
mentioned papers are essentially different from the one given here.


\section{Localized Euler classes on topological orbifolds}
For differentiable Fredholm sections of differentiable Banach bundles
over differentiable Banach manifolds a theory of localized Euler
classes has been developed recently by Brussee in \cite{brussee}
(motivated by previous work in this direction of Pidstrigatch and
Tyurin). The object of the present chapter is a refinement thereof
with regard to two aspects. We are going to\ \ 1)\ weaken
differentiability to local differentiability relative a finite
dimensional parameter space, and simultaneously implement a device to
globalize local transversality constructions\ \ 2)\ include the
orbifold case.


\subsection{Banach orbifolds}
The notion of orbifold is due to Satake \cite{satake}. Orbifolds
generalize manifolds in that as local models open subsets of vector
spaces are replaced by quotients of such by linear actions of finite
groups. The transfer to infinite dimensions and a topological instead of a
differentiable setting is straightforward, but for the sake of completeness
and to fix notations we give full details here. Partly we adopt the more
modern point of view of
\cite{kawasaki}.
\begin{defi}\rm
Let $\LUS$ be the category whose objects consist of tuples
$(q:\hat U\rightarrow U,G,\alpha)$ ({\em local uniformizing systems})
with
\begin{itemize}
\item
  $\hat U$ is an open set in some Banach space $T$
\item
  $\alpha:G \rightarrow\GL(T)$ is  a (not necessarily faithful) continuous
  linear action of the finite group $G$ on $\hat U$
\item
 $q$ is a topological quotient for $\hat U$ by the action of $G$.
\end{itemize}
For shortness we often write $U=\hat U/G$ for such local uniformizing
system. $U$ is called {\em support} of $U=\hat U/G$ and $q$
{\em structure map}.

A {\em morphism $f:(V={\hat V}/G') \rightarrow (U=\hat U/G)$ of local
uniformizing systems} is a homomorphism $\gamma: G'\rightarrow G$
and a $\gamma$-equivariant continuous map $\hat f:{\hat V}
\rightarrow \hat U$. The underlying continuous map $V\rightarrow U$ is
denoted $\bar f$.

$f$ is called {\em open embedding} (of $V={\hat V}/G'$ into $U=\hat U/G$)
if $\hat f$ is an open embedding and $\gamma$ is a monomorphism
fulfilling the maximality condition
\[
  \img(\gamma)=\{g\in G\mid\hat f(\hat V)\cap g\cdot\hat f(\hat   
  V)\neq\emptyset\}\, .
\]
\vspace{-8ex}

\qed
\end{defi}
There is a forgetful functor $\calq:\LUS\rightarrow\cals\mbox{ets}$
to the category of sets, sending a local uniformizing system $U=\hat U/G$
to its support $U$.

Now let $X$ be a Hausdorff space. If $\calu= \{U_i\}_{i\in I}$ is a
covering of $X$ by open sets we can form a category $\calt
(\calu)\subset \cals{\rm ets}$ with objects $U_i$ and a morphism $U_i
\rightarrow U_j$ for any inclusion $U_i \subset U_j$. If for any $i,j$
there exists $k$ with $U_k \subset U_i\cap U_j$ we call $\calu$ {\em
fine}.
\begin{defi}\rm\label{def_orbifold}
A {\em (Banach) orbifold structure} on a Hausdorff space $X$ is a fine
covering $\calu$ of $X$ and a functor $\calo:\calt(\calu)
\rightarrow\LUS$ with $\calq\circ\calo =\id_{\calt(\calu)}$ and with
$\calo(\iota)$ an open embedding for any $\iota\in\Hom(\calt(\calu))$.

A {\em morphism between (Banach) orbifold structures}
$F:(Y,\calv=\{V_i\}_{i\in I},\calo')\rightarrow(X,\calu=\{U_j\}_{j\in
J},\calo)$ is a map $\kappa:I\rightarrow J$ such that $V_i\subset V_j$ implies
$U_{\kappa(i)}\subset U_{\kappa(j)}$, together with morphisms of local
uniformizing systems $f_i:V_i ={\hat V}_i/G'_i\rightarrow
U_{\kappa(i)}={\hat U}_{\kappa(i)}/G_{\kappa(i)}$ such that
for any $i,j$ with $V_i\subset V_j$ there exists $g_{ij}\in
G_{\kappa(j)}$ with
\[
  \rho_{ji}\circ{\hat f}_i\ =\ g_{ij}\cdot{\hat f}_j\circ\sigma_{ji}\,,
\]
where $\rho_{ji}=\calo(U_{\kappa(i)}\hookrightarrow U_{\kappa(j)})$,
$\sigma_{ji}=\calo' (V_i\hookrightarrow V_j)$. $F$ induces a well-defined
continuous map from $Y$ to $X$, that we denote by $\bar F$.  If all the
$f_i$ are open embeddings, $F$ is also called {\em open embedding} (of
(Banach) orbifold structures).

Two orbifold structures $(\calu',\calo')$ and $(\calu,\calo)$ are {\em
pre-equivalent} if there is an open embedding $F:(X,\calu',\calo')
\rightarrow (X,\calu,\calo)$ with $\bar F=\id_X$. A {\em (Banach)
orbifold} $(X,[\calu,\calo])$ is a Hausdorff space $X$ and an
equivalence class $[\calu,\calo]$ of Banach orbifold structures
with respect to the equivalence relation generated by this notion
of pre-equivalence.

Likewise, a {\em morphism of Banach orbifolds} $F:(Y,[\calv,\calo'])
\rightarrow (X,[\calu,\calo])$ is an equivalence class of morphisms of
orbifold structures in the classes $[\calv,\calo']$ on $Y$ and
$[\calu,\calo]$ on $X$, the equivalence generated by compatible
open embeddings of orbifold structures in domain and image.

A {\em submersion} $P:X\rightarrow S$ of Banach orbifolds is a morphism
that is represented by splittable projections, i.e.\ for any locally
representing $p:U=\hat U/G\rightarrow U_S={\hat U}_S/G_S$, there exists
a decomposition $\hat U=U_S\times U_P$ such that $p$ is projection to
the first factor.  \qed \end{defi}
\begin{rem}\rm
1)\ \ The definition of orbifold is slightly complicated by the fact that two
local models may only be comparable on smaller sets. As an example
consider $S^2$ with cyclic quotient singularities of coprime orders
$m,n>1$ at the poles $0, \infty \in S^2$. This means that $S^2$ is covered
via stereographic projection by two local uniformizers
\[
  \hat V=\cz\longrightarrow V=\cz/ \gz_m \simeq \cz\simeq S^2
  \setminus\{0\},\quad
  \hat U=\cz\longrightarrow V=\cz/ \gz_n \simeq
  \cz\simeq S^2\setminus\{\infty\}
\]
with $\gz^k$ acting on $\cz$ by $k$-th roots of unity. But since
$(m,n)=1$ there is no open embedding of the restrictions to $U$ of
$V= \hat V/ \gz^m$ to $U=\hat U/ \gz^n$. So unlike for manifolds,
comparison of two local uniformizing systems might be possible only
locally. The introduction of categories of open sets is a convenient way to
handle this problem.\\[1ex]
2)\ \ Other definitions of orbifolds usually require
the action of the finite group to be faithful. This is inappropriate for cases
as $\calm_{g,k}$ of moduli spaces of $k$-marked curves of genus $g$ with
$(g,k)\in\{(1,1),(2,0)\}$ where the generic curve has automorphisms.

The category of orbifolds with effective actions also does not lead to
satisfactory treatment of certain natural constructions, e.g.\ zero
sections of orbibundles, cf.\ Definition~\ref{def_orbibundle} below.
\qed
\end{rem}
As a matter of notation we usually refer to representatives (local
uniformizers) of any objects on local uniformizing systems by adding a
hat. And if we deal with a fixed orbifold $(X,[\calu,\calo])$ and talk
about local uniformizing systems, we always mean local uniformizing
systems belonging to an orbifold structure in $[\calu,\calo]$.  We will
often suppress representatives of the orbifold structure and just talk
of $X$ being an orbifold. When the structure map and action of local
groups are understood we occasionally call $\hat U$ a {\em chart} of $X$.

To any $x\in X$ belongs a group $G_x$, unique up to isomorphism, the
stabilizer of any lifting $\hat x\in \hat U$ of $x$. In fact, by the
maximality condition, any open embedding $\sigma:U'={\hat
U}'/G'\rightarrow U=\hat U/G$ of local uniformizing systems establishes
an isomorphism of stabilizers $G'_{\hat x} \simeq G_{\hat\sigma(\hat
x)}$. A local uniformizing system $U=\hat U/G$ is said to be {\em
centered at $x$} if $G=G_x$ (i.e.\ $\sharp q^{-1}(x)=1$, $q:\hat U
\rightarrow U$) and $\hat x=0\in T$.

If $(q:\hat U \rightarrow U,G,\alpha)$ is a local uniformizing system
centered at $x$ and $V\subset U$ is an open neighbourhood of some point
$y\in U$, the restriction $(q:q^{-1}(V) \rightarrow V,G,\alpha)$ of $U=
\hat U/G$ to $V$ is not in general centered in $y$. Let $\hat y \in\hat
U$ be a lift of $y$. We claim that by taking the connected component
$\hat V$ of $q^{-1}(V)$ containing $\hat y$ and $V$ sufficiently small,
$(q:\hat V \rightarrow V, G_{\hat y}, \alpha|_{G_{\hat y}})$ is a local
uniformizing system centered at $y$ and compatible with $U=\hat U/G$
(here the translation of the modelling Banach space $T$ about $-\hat y$
that is necessary to achieve $\hat y=0$ is understood). In fact,
choose representatives $g_1=e,g_2,\ldots,g_r$ of $G/G_{\hat y}$. Then
there is an open neighbourhood $W$ of $\hat y$ with $g_i W\cap W=
\emptyset$ for $i>1$. Set $\hat V:= \bigcap_{g\in G_{\hat y}} g\cdot W$.
Letting $V=\hat V/G_{\hat y} \hookrightarrow U=\hat U/G$ be defined by
the inclusions $\hat V \subset \hat U$, $G_{\hat y}\subset G$, the only
thing to be checked is that for any $\hat z,\hat w\in \hat V$ with $g\cdot
\hat z=\hat w$ it follows $g\in G_{\hat y}$. Writing $g=g_i\cdot h$, $h\in
G_{\hat y}$, we obtain $\hat w\in g_i\hat V\cap\hat V \subset g_i W\cap
W$, and hence $g_i=e$. We call $V=\hat V/G_y$ a {\em proper restriction}
of $U=\hat U/G$, written $(U=\hat U/G)|_{y,V}$, where it is understood that
$V$ is chosen sufficiently small.

We may thus restrict to such orbifold structures in the sequel having for
any $x\in X$ local uniformizing systems centered in $x$. $x\in
X$ is called {\em regular} if there exists a local uniformizing system
$U=\hat U/G$ centered in $x$ with $G_x=\ker(G\rightarrow\GL(T))$,
else {\em singular}.  The corresponding subsets of $X$ will be denoted
$X_\reg$ and $X_\sing$ respectively.
\begin{rem}\rm\label{construct_orbifold}
To construct an orbifold structure on a set $X$ the following data
suffices:
\begin{itemize}
\item
  a covering of $X$ by supports of local uniformizing systems $U_i
  ={\hat U}_i /G_i$, $i\in I$, with $\codim U_{i,\sing}\ge 2$ (!)
\item
  for any $i,j\in I$, $x\in U_i\cap U_j$, a sufficiently small open
  neighbourhood $V\subset U_i\cap U_j$ of $x$ (such that a proper
  restriction $V=\hat V/G_x$ of $U_i={\hat U}_i/G_i$ to $V$ exists) and
  an open embedding
  \[
    \sigma:(V=\hat V/G_x) \longrightarrow (U_j={\hat U}_j/G_j)\ ,
  \]
  with $\bar\sigma: V\hookrightarrow U_j$ the inclusion.
\end{itemize}
Take as elements of the covering $\calu$ all open sets $U\subset V$ for
some $V$ as before, with $x\in U$ and with $\hat U \subset\hat V$ {\em
starlike} with respect to $\hat x\in\hat V$. Define $\calo:
\calt(\calu)\rightarrow \LUS$ on objects by sending $U\in \calu$ to
$(V=\hat V/G_x)|_{x,U}$. And if $U',U\in\calu$ with $U'\subset U$, an
open embedding of $\calo(U')$ into $\calo(U)$ can be obtained by gluing
the existing open embeddings of restrictions of $\calo(U')$ into
$\calo(U)$. This is possible by the following Lemma~\ref{glue_lus}.

Conversely, any orbifold $X$ can be defined in a neighbourhood of any
compact set $K\subset X$ in this way, with a finite index set.
\qed
\end{rem}
\begin{lemma}\label{glue_lus}
  Let $U'={\hat U}'/G'$, $U=\hat U/G\in\LUS$, with $U'\subset U$ open,
  $\codim{\hat U}'_\sing\ge2$ and with ${\hat U}'\subset T'$ starlike
  with respect to $0\in T'$. Let also be given a covering of $U'$ by connected
  open sets $W_j\subset U$, $j\in J$, and for each $j$ a point $z_j\in
  W_j$, and an open embedding $\sigma_j$ of the proper restriction
  $(W_j={\hat W}_j/G'_{z_j}) = (U'={\hat U}'/G')|_{z_j,W_j}$ into
  $U=\hat U/G$.

  Then there are $g_j\in G$ such that $g_j\cdot{\hat\sigma}_j$ fit
  together to an open embedding $\sigma:(U'={\hat U}'/G') \hookrightarrow
  (U=\hat U/G)$. $\sigma$ is unique up to left-multiplication by $g\in G$.
\end{lemma}
\pf
Choose $j_0\in J$ with $0\in{\hat W}_{j_0}$. Then $G'_{{\hat
z}_{j_0}}=G'$, so $\sigma_{j_0}$ defines a monomorphism
$\gamma:G'\hookrightarrow G$. The other assertions concern only the
action of $G'$ and $G$ on ${\hat U}'$ and $\hat U$ respectively, so we
may henceforth assume the actions to be faithful. To construct $g_j$,
cover the path $q(\{t\cdot{\hat z}_j\mid t\in[0,1]\})$ by finitely many
of the $W_j$, say $W_{j_0},W_{j_1},\ldots,W_{j_k}$ with $W_{j_\nu}\cap
W_{j_{\nu+1}}\neq\emptyset$. Possibly by replacing ${\hat W}_{j_\nu}$
by $g'_\nu\cdot {\hat W}_{j_\nu}$ we may assume the ${\hat W}_{j_\nu}$
connnect ${\hat z}_j$ with $0$, i.e.\ ${\hat W}_{j_\nu}\cap {\hat
W}_{j_{\nu+1}} \neq\emptyset$ and ${\hat z}_j\in{\hat W}_{j_k}$.  Since
for any $j$, $(W_j)_\sing$ does not disconnnect $W_j$ (this is clear in
finite dimensions and here by restriction to sufficiently big,
finite-dimensional, $G$-invariant linear subspaces of $T'$), it holds
${\hat\sigma}_{j_\nu} ={\bar g}_{\nu+1}\cdot {\hat\sigma}_{j_{\nu+1}}$ on
${\hat W}_{j_\nu}\cap{\hat W}_{j_{\nu+1}}$ for uniquely determined ${\bar
g}_\nu\in G$, $\nu=1,\ldots,k$ (derive uniqueness first over $U_\reg$ as in
\cite[Lemma~1]{satake} and extend by continuity). We set $g_j:={\bar
g}_1\cdot\ldots\cdot{\bar g}_j$. In view of the effectiveness of the
action of $G$, $g_j$ does not depend on choices other than $j_0$, and
the $g_j\cdot{\hat\sigma}_j$ coincide on common domains of definition.
It remains to be shown that the map $\hat \rho:\hat U' \rightarrow \hat U$
thus obtained is $G'$-equivariant and injective. Away from $U_\sing$,
$\hat \rho$ is a map of unbranched coverings. The covering groups are $G'$
and $G$ respectively. Since the inclusion of covering groups is nothing but
$\gamma$, we conclude both $G'$-equivariance and injectivity away from
$U_\sing$. Equivariance extends to the whole domain by continuity. We infer
injectivity by observing that since $\hat \rho$ is an open map, the set of injective
points is in fact closed.
\qed
\begin{rem}\rm\label{U_reg}
The codimension of the singular locus of $U=\hat U/G$ is at least two
iff $G$ does not act by any reflections. In fact, the action being
linear, ${\hat U}_{\mbox{\scriptsize sing}}$ can be stratified by
locally closed submanifolds, which are one-codimensional iff there
exists $g\in G$ with $\codim\mbox{Fix}(g)=1$.  This happens iff $g$
comes from an action on $\rz$ by an epimorphism $T\rightarrow\rz$,
hence acts by a reflection (choose a $G$-invariant complement to
$\ker(T\rightarrow T,v\mapsto g\cdot v)$).

As a consequence, the hypothesis of the lemma is fulfilled for complex
linear actions.
\qed
\end{rem}
\vspace{1ex}

We shall also use the notion of orbi (vector) bundles.
\begin{defi}\rm\label{def_orbibundle}
Let $p:E\rightarrow X$ be a continuous surjection of topological
spaces. A {\em (Banach) orbibundle structure} on $p$ is a morphism of
Banach orbifold structures
\[
  P:(E,\{p^{-1}(U)\}_{U\in\calu},\calo^E) \longrightarrow (X,\calu,\calo)
\]
on $E$ and $X$ with
\begin{itemize}
\item
  $\kappa=\id_I$
\item
  if $U\in\calu$ and $\calo(U)= (U=\hat U/G)$ then
  $\calo^E(p^{-1}(U)) =(\hat U \times E_0)/G^E$ with $E_0$ a Banach
  space and $G^E$ acting diagonally on $\hat U\times E_0$ via an
  epimorphism $G^E \rightarrow G$ and a continuous representation
  $G^E \rightarrow \GL(E_0)$.
\end{itemize}
We call $\calo^E(U)$ a {\em local (uniformizing) trivialization} of the
orbibundle structure over $U$, written $\hat U\times E_0/G^E
\downarrow U=\hat U/G$.

A {\em morphism of Banach orbibundle structures} $P'$ on $p':
F\rightarrow Y$ and $P$ on $p:E\rightarrow X$ is a morphism of Banach
orbifold structures
\[
  H:(F,\{{p'}^{-1}(V)\}_{V\in\calv},\calo^F) \longrightarrow   
  (E,\{p^{-1}(U)\}_{U\in\calu},\calo^E)
\]
with local uniformizers of the form
\[
  h:({\hat F}_V:=\hat V\times F_0)\longrightarrow ({\hat E}_U:= \hat
  U\times E_0), \quad (x,w) \longmapsto (\hat f(x),A_x(w))
\]
for some $\hat f:\hat V\rightarrow\hat U$ and $A: \hat V\rightarrow
L(F_0,E_0)$, i.e.\ $A_x$ linear for any $x\in \hat V$.

A {\em (Banach) orbibundle} (by abuse of notation denoted
$p:E\rightarrow X$, or just $E$ if the base $X$ is understood) is an
equivalence class of Banach orbibundle structures with respect to the
equivalence generated by open embeddings of $E$ into itself inducing the
identity on the underlying sets. Similarly, {\em morphisms of Banach
orbibundles} are defined as equivalence class of morphisms of Banach
orbibundle structures respecting the linear structure.

The Banach orbifold $E$ naturally associated to a Banach orbibundle $p:
E\rightarrow X$ (i.e.\ forgetting the structure relative $X$) is called
{\em total space}. The projection $p:E\rightarrow X$ of sets underlies
a submersion of Banach orbifolds (also denoted $p$), locally
contracting ${\hat E}_0$ to a point.

If $X$ is already a Banach orbifold, $p:E \rightarrow X$ is said to be
a Banach orbibundle {\em over} $X$ if the orbifold structure on $X$
induced by the maps $\hat f: {\hat U}'\rightarrow \hat U$ is in the class of
the given orbifold structure on $X$.

A {\em section} of $E$ is a morphism $s:X\rightarrow E$ of orbifolds
with $p\circ s= \id_X$. If all epimorphisms $G^E\rightarrow G$ are
split there is a canonical morphism $s_0:X\rightarrow E$, which is a
section, the {\em zero section}. Otherwise there is only a morphism
$s_0:\tilde X\rightarrow E$ with $p:\tilde X\rightarrow X$ a bijective
submersion of Banach orbifolds (replace $G$ by $G^E$ in the local
uniformizing systems defining $X$). The {\em zero locus} $Z(s)$ of a
section $s$ is the set $Z(s):=s^{-1}(\img s_0)$.
\qed
\end{defi}
\begin{rem}\rm\label{cp-open-top}
We require here only continuity of the map of total spaces $h:{\hat
F}_V \rightarrow {\hat E}_U$. Since our spaces are all metric, hence
compactly generated, this is equivalent to continuity of the associated map
\[
  \Phi_h:\hat U\rightarrow L({\hat F}_0,{\hat E}_0)
\]
with $L({\hat F}_0,{\hat E_0})$ endowed with the compact-open topology
(for general spaces one needs to go over to the compactly generated
topology on the mapping space). We refer to \cite[\S1.4]{whitehead} for an
appropriate discussion of these matters.

The tempting condition of continuity of $\Phi_h$ with respect to the strong
topology (the operator norm topology) is considerably more restrictive.
In fact, an instance where we can not work with norm topologies will be the
bundles $q^p_*E$ and ${q_1^p}_*E$ in Chapter~\ref{chapter_dbar}.
\qed
\end{rem}
Ordinary Banach bundles over the topological space underlying an
orbifold can of course be viewed as orbibundles (with trivial
representations $G^E \rightarrow \GL(E_0)$). More interesting are
orbibundles that do not arise in this way. Then there are fibers
$E_x=p^{-1}(x)$ (not to be confused with their uniformizers ${\hat
E}_{\hat x} \simeq E_0$) that arise as quotients of $E_0$ by
non-trivial group actions. Usually this will destroy the additive
structure on $E_x$. Typical examples are tangent bundles of
{\em differentiable} orbifolds (i.e.\ where the open embeddings are
required to be differentiable).
\begin{rem}\rm
To construct an orbibundle structure on a set of maps $p:E\rightarrow
X$ we may go along the lines of Remark~\ref{construct_orbifold}. The
data needed are now coverings $\{U_i={\hat U}_i/G_i\}$ of $X$ and
$\{E_{U_i}= {\hat U}_i\times E_0/G_i^E\}$ of $E$, $G_i$ epimorphic
image of $G_i^E$, such that $p$ is locally uniformized by the
projection ${\hat U}_i\times E_0 \rightarrow {\hat U}_i$. Comparison
maps (the second point in Remark~\ref{construct_orbifold}) are required
to respect the linear structure on $E_0$.

In fact, by Remark~\ref{construct_orbifold} we get orbifold structures
on $E$ and $X$. By construction these come together with a morphism of
orbifolds with underlying map $p$, which is of the form specified in
Definition~\ref{def_orbibundle}.
\qed
\end{rem}

Banach orbibundles give rise to associated orbibundles in the familiar
way; e.g.\ if $T$, $E$ are Banach orbibundles over $X$, there is the
bundle of continuous linear maps $L(T,E)$ with ${(L(T,E))_x^\wedge}
=\Hom({\hat T}_x,{\hat E}_x)$. And if $f:Y\rightarrow X$ is a morphism of
Banach orbifolds, an orbibundle $E$ on $X$ pulls back to an orbibundle
$f^*E$ over $Y$ such that there is a fiberwise isomorphic morphism
$f^*E\rightarrow E$ over $f$. But note that if $H\in L(T,E)$, a dual
$H^\vee: E^\vee=L(E,\underline{\rz}) \rightarrow T^\vee= L(T,
\underline{\rz})$ does not in general exist unless $H$ specifies
isomorphisms of the local groups of $T$ and $E$. Here $\underline{\rz}$
stands for the trivial bundle $X\times \rz \rightarrow X$.


\subsection{Localized Euler classes in finite dimensions}
For a locally compact space $X$ let $H_*(X)$ denote Borel-Moore
homology with coefficients in the ring $\qz$ \cite{bredon},
\cite{iversen} (or, equivalently in our case, singular homology of the
second kind, i.e.\ with locally finite singular chains,
cf.\ \cite{sklya}), and $H_\ph^*(X)$ (sheaf) cohomology with supports
$\ph$, where we always stick to rational coefficients. The choice of
homology theory is suggested by the need of fundamental classes for
non-compact manifolds (or rather orbifolds). $H_*(X)$ has functorial
properties that are closer to cohomology than homology: For inclusions
$i:U\hookrightarrow X$ of an open subset there is a restriction map
$i^!:H_*(X) \rightarrow H_*(U)$; push-forwards $f_*$ exists only for
{\em proper} maps $f:X \rightarrow Y$. If $Z$ is a closed subset of $X$
there is a remarkable long exact sequence \[
  \ldots\longrightarrow H_{k+1}(U)\longrightarrow H_k(Z) 
  \stackrel{j_*}{\longrightarrow} H_k(X)
  \stackrel{i^!}{\longrightarrow} H_k(U) \ldots\ ,
\]
where $j:Z\hookrightarrow X$ and $i:U:=X\setminus Y \hookrightarrow X$
are the inclusions. And for a closed subset $Z\subset X$ there are
(localized) cap products $H_{d+k}(X)\otimes H_Z^k(X)\rightarrow
H_d(Z)$, making $H_*(X)$ into a {\em right} module over $H^*(X)$
($Z=X$). The following vanishing theorem will also be used: If the
(covering) dimension of $X$ does not exceed $d$ then $H_k(X)=0$ for any
$k>d$. We refer to \cite{bredon}, \cite{iversen} and \cite{sklya}
for a presentation of these and related facts to be used in the sequel
without further notice.

Some remarks on fundamental and Thom classes of orbifolds and
orbibundles in finite dimensions are now in order:  Let $X$ be a
(connected for simplicity) topological orbifold of dimension $d$. We
call $X$ {\em locally orientable} (as orbifold) if for any $x\in X$,
$G_x$ acts orientation preservingly on ${\hat U}_x$ (this does not
depend on choices). Since this property implies $\dim X_\sing\le d-2$,
$X$ is also locally orientable as a topological space, i.e.\ the sheaf
$\calh_d(X)$ associated to the presheaf $U\mapsto H_d(U)$ is locally
free of rank 1. In fact, from
\[
  H_d(U_\sing)\longrightarrow H_d(U)\longrightarrow
  H_d(U_\reg) \longrightarrow H_{d-1}(U_\sing)
\]
and the dimension assumption we get $H_d(U)\simeq H_d(U_\reg)$, while
$U_\reg$ is a connected manifold iff $U$ was connected, and orientable
for $U$ sufficiently small.  An {\em orientation} for $X$ is the choice
of a trivialization of $\calh_d(X)$. Given an orientation, one may
combine the above isomorphism applied to $U=X$ and Poincar\'e duality
on $X_\reg$ to conclude
\[
  H_d(X)\simeq H_d(X_\reg,\calh_d(X)^\vee)\simeq H^0(X_\reg)
  \simeq\qz\, .
\]
In this case the {\em fundamental class} $[X]\in H_d(X)$, defined as
corresponding to $1\in H^0(X_\reg)$, restricts to a generator of
$H_d(U)$ for any connected $U\subset X$ (this all works with integral
coefficients). Notice that an orientation on $X$ is equivalent to the
giving of orientations on $\hat U$ for a defining set of local
uniformizing systems $U=\hat U/G$ (in the sense of
Remark~\ref{construct_orbifold}) with $G$ acting orientation
preservingly, compatible with respect to open embeddings.

Likewise, if $({\hat F}_U=\hat U\times{\hat F}_0)/G\downarrow U=\hat
U/\bar G$ is a local trivialization with $\rk{\hat F}_U=r$, $\Pi:{\hat
F}_U\rightarrow{\hat F}_0$ the projection, choosing an orientation for
${\hat F}_0$ selects an integral generator $\Pi^*\delta$ of $H^r_{\hat
U}({\hat F}_U)\simeq\qz$, the {\em (local) Thom class} $\Theta_{{\hat
F}_U}$, where $\delta\in H^r_{\{0\}}({\hat F}_0)$ is Poincar\'e dual to
$[0]\in H_0(\{0\})$. It is characterized by $[{\hat F}_U]\cap
\Theta_{{\hat F}_U} =[\hat U]$ for any choice of orientation on $\hat
U$ (taking the product orientation on ${\hat F}_U$). Now the quotient
$q_0:{\hat F}_0\rightarrow F_0={\hat F}_0/G$ induces a transfer
isomorphism $q_0^*:H^r_{\{0\}}(F_0)\simeq H^r_{\{0\}}({\hat F}_0)^G$
(cf.\ \cite[Thm.19.1]{bredon}; here we need rational coefficients). The
superscript ``G'' indicates the submodule of $G$-invariant classes.
If $G$ acts on ${\hat F}_0$ in an orientation preserving way,
$\Theta_{{\hat F}_U}$ is $G$-invariant, and so there exists
$\bar\delta\in H^r_{\{0\}}(F_0)$ with $\delta=q_0^*\bar\delta$.  We put
$\Theta_{F_U}:=(b/\bar b)\cdot{\bar\Pi}^*\bar\delta$, where $\bar\Pi:
F_U={\hat F}_U/G\rightarrow F_0$ is the map induced by $\Pi$, and
$b=\sharp G/\ker (G\rightarrow\GL({\hat F}_U))$, $\bar b=\sharp \bar G/
\ker (\bar G\rightarrow\GL(\hat U))$ are the covering degrees of ${\hat
F}_U\rightarrow F_U$ and $\hat U\rightarrow U$ respectively. Note that
$b/\bar b$ is nothing but the covering degree of ${\hat F}_{\hat
x}\rightarrow F_x$ for $x$ generic. Hence $\Theta_{F_U}$ is natural
with respect to open embeddings of local uniformizing systems. 
The constant has been chosen in such a way that if $\bar G$ acts
orientation preservingly on $\hat U$ then $[F_U]\cap\Theta_{F_U}=[U]$ for
any choice of orientation on $\hat U$, as one easily checks using the
projection formula.

Let us call an orbibundle $F$ over $X$ {\em locally orientable} if
$G_x$ acts orientation preservingly on ${\hat F}_x$ for any $x\in X$.
The transfer isomorphism shows that in this case the sheaf
$\calh_X^r(F)$ associated to the presheaf $U\mapsto H^r_U(F_U)$ is
locally free of rank 1. An {\em orientation} for $F$ is a
trivialization of $\calh_X^r(F)$. Obviously, $F$ admits an orientation
iff the orbibundle $\Lambda^r F$ is a trivial real line bundle. If
$U,V\subset X$ are open sets with one of $U,V\in\calu$, the
Mayer-Vietoris sequence
\[
  H^{r-1}_{U\cap V}(F_{U\cap V})\longrightarrow H^{r}_{U\cup V}
  (F_{U\cup V})\longrightarrow H^r_U(F_U)\oplus
  H^r_V(F_V)\longrightarrow H^r_{U\cap V}(F_{U\cap V})
\]
has vanishing term on the left (by transfer), and thus classes
$\alpha\in H_U^r(F_U)$, $\beta\in H^r_V(F_V)$ glue (uniquely) iff
$\alpha|_{F_{U\cap V}}=\beta|_{F_{U\cap V}}$.  So provided $F$ is
oriented we may glue the local Thom classes (with appropriate signs
matching the orientation) to a unique {\em (global) Thom class}
$\Theta_F\in H^r_X(F)$, a global generator of $\calh_X^r(F)$. Note that
if $X$ is oriented too, then so is the total space $F$ and it holds
\[
  [F]\cap\Theta_F=[X]\, .
\]
By compatibility of cap products with open embeddings it is enough to
check this locally, where it is clear by the foregoing. Another way to
see this is to observe that any locally oriented orbifold obeys
rational Poincar\'e duality. This follows from the computation of local
cohomology groups using transfer. $\Theta_F$ is then the class dual on
$F$ to $[X]\in H_d(X)$, $X$ embedded into $F$ as zero section.

Now if $s$ is a section of an oriented orbibundle $F$ of rank $r$ over
an oriented $n$-dimensional topological orbifold $X$, the {\em
localized Euler class of $s$} is the class
\[
  [F,s]:=[X]\cap s^*\Theta_F\in H_{n-r}(Z(s))\, .
\]
The point is that $[F,s]$ lives on the zero set of $s$. If $s$ is
locally given by a map $\kappa:\hat U\rightarrow{\hat F}_x$, $\hat
U\subset\rz^n$ uniformizing an open neighbourhood $U$ of $x\in X$, then
by compatibility of cap products with open embeddings
\[
  [F,s]|_U=q_*([\hat U]\cap\kappa^*\delta)\,,
\]
where as above $\delta\in H^r_{\{0\}} ({\hat F}_x)$ is Poincar\'e dual to
$[0]\in H_0(\{0\})$ and $q:\hat U\rightarrow U$ is the structure map. But
note that the local classes do not determine $[F,s]$ uniquely unless $\dim
Z(s)=n-r$.


\subsection{Abstract transversality}
In infinite dimensions extra conditions are needed to reduce to a
finite dimensional situation. The basic tool is the implicit function
theorem, that we state here in the spirit of \cite{kuranishi}, extended
to the case of equivariant families of differentiable maps with
splittable differentials.
\begin{prop}\label{diff_kuranishi}
Let $T$ and $E$ be Banach spaces, $S$ a Hausdorff space, $V\subset T$ 
an open neighbourhood of $0$ and $H:S\times V\rightarrow E$
continuous and differentiable in the second variable with
differential $D_V H$ uniformly continuous at $(s_0,0)$. Let also
be  given decompositions $T\simeq K \times P$, $E\simeq
C\times Q$ with $\sigma:P\rightarrow Q$ an isomorphism,
$\sigma=D_V H(s_0,0)|_{\{0\}\times P}$.

Then there are open sets $V'\subset T$, $S'\subset S$, an open
embedding $\Phi:S'\times V' \rightarrow S'\times V$ of the form
$\Phi(s,k,p) =(s,k,\ph(s,k,p))$ with $\Phi(s_0,0) = (s_0,0)$ and a
continuous map $\kappa:S'\times V'\rightarrow C$ with 
\[
  H\circ\Phi=(\kappa,\sigma\circ\prj_2):S'\times V'\longrightarrow
  C\times Q\,,
\]
$\prj_2: T\simeq K\times P \rightarrow P$ the projection.

Moreover, if $H$ is $G$-equivariant for a finite group $G$ and $K$ is
$G$-invariant, then $P$ and $C$ may be chosen $G$-invariant, in the
case of which $\Phi$ and $\kappa$ become $G$-equivariant.
\end{prop}
\pf
We may assume $V=V_K\times V_P$, $V_K\subset K$, $V_P\subset P$.
Let $q: E\rightarrow Q$ be the projection along $C$. Then the derivative
in the $P$-direction of
\[
  q\circ H: S\times V_K\times V_P\longrightarrow Q
\]
is norm-continuous at the point $(s_0,0)$ with bijective differential.
Applying the ordinary implicit function theorem (cf.\ e.g.\
\cite[Thm.4B]{zeidler}) we obtain, possibly after shrinking $S$, the
injection $\Phi: S\times V' \rightarrow S\times V$ with $q\circ H\circ
\Phi =q\circ\sigma$. Note that the last equation determines the germ
of $\Phi$ at $(s_0,0)$ uniquely. One defines $\kappa:= H\circ\Phi -\sigma$,
which by the former equation indeed takes values in $C$.

In the equivariant case we may assume $(s_0,0)$ to be a fixed point of
the action. Let $\pi:T\rightarrow K$ be the projection along $P$. Then
\[
  \bar\pi:T\longrightarrow K\,,\quad x\longmapsto
  \frac{1}{|G|}\sum_{g\in G}g^{-1}\pi(gx)
\]
is $G$-equivariant with $\bar\pi|_K=\id$, and so by replacing $P$ by
$\ker\bar\pi$ we achieve $P$ to be $G$-invariant. In turn $Q=\sigma(P)$
is $G$-invariant too, and thus, by an analogous argument, has a
$G$-invariant complement $C$. Equivariance of $\Phi$ and $\kappa$ are
then immediate from the construction.
\qed
\vspace{1ex}

We will not be able to use the implicit function theorem directly in our
problem, since differentiability relative the natural base space
$\calm_{g,k}$ fails when bubbles (cf.\ Chapter~3) occur. However,
differentiability holds locally up to some finite-dimensional space. This
leads to the following definition that for purposes not pursued in this
paper is formulated relative to some possibly infinite dimensional base
space.
\begin{defi}\label{Fredholm_structure}\rm
Let $p:X\rightarrow S$ be a submersion of Banach orbifolds, $s$ a
section of a Banach orbibundle $E$ over $X$. A {\em Fredholm structure}
for $s$ {\em relative $S$} is a choice of orbifold structures for $X$,
$S$ and $E$ such that any local trivialization $\hat E_U =\hat U\times
E_0/ G^E\downarrow U=\hat U/G$ centered in some $z\in Z(s)$ has the
form
\[
  \hat U\ =\ W\times L\times V
\]
with
\begin{itemize}
\item
  the projection $\hat U\rightarrow W$ locally uniformizes $p$
\item
  $L$ is an open subset of a finite dimensional vector space
\item
  $\Pi\circ\hat s: W\times L\times V \rightarrow E_0$ is differentiable
  relative $W\times L$ with relative differential
  \[
    D_V(\Pi\circ\hat s):\ W\times L\times V\ \longrightarrow\ L(T,E_0),
    \quad T=T_0 V
  \]
  continuous at 0 and with $\sigma= D_V(\Pi\circ \hat s)(0)$ Fredholm.
\end{itemize}
We will refer to a pair $(s, \mbox{Fredholm structure for $s$ relative
$S$})$ by saying that $s$ is an {\em $S$-Fredholm section}. Local
trivializations and local uniformizing systems belonging to the
Fredholm structure will be called {\em distinguished}.
\qed
\end{defi}
\begin{rem}\rm
The concept of Fredholm structure comprises two things: First,
differentiability of the section relative $S$ up to some finite
dimensional factor. And second, a Fredholm condition for the relative
differential. Since $L$ is finite dimensional this condition is not
affected by making $L$ bigger.

Two extreme cases should be kept in mind here. One is the
finite dimensional situation, where both conditions are empty. The
other is that of differentiable sections over differentiable
Banach orbifolds with Fredholm linearizations. The latter case is familiar
from global analysis, e.g.\ in gauge theory. For the case of stable curves
in manifolds we have well-defined differentiable structures only after
fixing the complex structure of the domain, which is why we need the
hybrid concept.
\qed
\end{rem}
Given the notion of Fredholm structures we can say when a section
should be considered transverse to the zero section.
\begin{defi}\label{transverse}\rm
An $S$-Fredholm section $s$ of $E$ over $p:X \rightarrow S$
is {\em transverse} along a closed subset $A\subset Z(s)$ iff for any
$z\in A$ there exists a distinguished local trivialization centered in $z$
with $\sigma$ surjective.
\qed
\end{defi}
By applying the implicit function theorem locally relative $W\times L$,
we see that the zero locus of transverse sections is a finite dimensional
topological orbifold, locally uniformized by $W\times L \times K$,
$K=\ker \sigma$ (hence submerges onto $S$):
\begin{prop}\label{transverse_section}
  Let $s$ be an $S$-Fredholm section of a Banach orbibundle $E$ over a
  Banach orbifold $X$ submerging onto $S$, and assume $s$ is transverse
  along $A$. Then in a neighbourhood of $A$, $Z(s)$ has naturally a
  structure of topological orbifold, submerging onto $S$ with finite
  dimensional fibers.
\qed
\end{prop}
For abstract transversality, we want to exhibit $Z(s)$ as zero set of a
section of a finite rank bundle over a finite dimensional orbifold. To this
end we now introduce another technical condition. For a Fredholm
operator $\sigma:T\rightarrow E$ we say at set of vectors $\{v_i\}_{i\in
I}$ {\em span $\coker\sigma$} if $E$ is generated by $\img\sigma$ and
$v_i$, $i\in I$.
\begin{defi}\label{Kuranishi_structure}\rm
Let $p:X\rightarrow S$ be a submersion of Banach orbifolds and $s$ an
$S$-Fredholm section of a Banach orbibundle $E$ over $X$. A {\em
Kuranishi structure} for $s$ relative $S$ is a morphism
\[
  \tau: F\longrightarrow W
\]
from a finite rank orbibundle $F$ defined over an open suborbifold
$X'\subset X$ containing $Z(s)$ such that for any distinguished local
trivialization $\hat E_U= \hat U\times E_0/G^E\downarrow U=\hat U/G$
centered in some $z\in Z(s)$, $\hat U=W\times L\times V$:
\begin{itemize}
\item
  $\img \hat\tau_z$ spans $\coker \sigma$, $\sigma=
  D_V(\Pi\circ \hat s)(0)$
\item
  $\hat \tau$ is continuously differentiable relative $W\times L$.
\end{itemize}
Two Kuranishi structures $\tau:F\rightarrow E$, $\tau': F'\rightarrow
E$ are {\em compatible} iff $\tau+\tau': F\oplus F'\rightarrow E$ is a
Kuranishi structure too. An $S$-Fredholm section together with an
equivalence class of compatible Kuranishi structures is called {\em
$S$-Kuranishi} section.
\qed
\end{defi}
\begin{rem}\rm
1)\ \ For a distinguished local trivialization of an $S$-Fredholm section
$s$ put $\sigma=D_V (\Pi\circ\hat s) (0,0,0)$ and choose $G$-invariant
splittings $T=K\times P$, $E_0=Q\times C$ with $K=\ker\sigma$,
$Q=\img\sigma$. These exist by finite-dimensionality of $K$ and of
$\coker\sigma$. We may then apply the implicit function theorem
Proposition~\ref{diff_kuranishi} to $\Pi\circ\hat s$ relative $W\times L$
to conclude (possibly after shrinking $W$ and $L$) the existence of
maps $\Phi:W\times L\times V' \rightarrow W \times L\times V$,
$\kappa:W\times L\times V'\rightarrow C$ with $V'\subset T$ open and
\[
  \Pi\circ\hat s\circ\Phi\ =\ (\kappa,\sigma):W\times L\times V'
  \longrightarrow C\times Q\, .
\]
We hereby obtain a {\em Kuranishi model} for $s$, or rather $Z(s)$, near
$z$: Up to a change of coordinates via $\Phi$, $\hat s$ splits off an
infinite dimensional invertible linear part $\sigma$ such that on
$\sigma^{-1}(0)=W\times L\times K$, $\hat s$ is given by the
continuous map $\kappa|_{W \times L\times K}$ between open sets in
finite dimensional vector spaces. In particular, it holds $Z(\hat s)
=\Phi((\kappa|_K)^{-1}(0))$ near $z$.

As we will see below the morphism $\tau$ in
Definition~\ref{Kuranishi_structure} enables us to globalize this
construction. This is why we called $\tau$ ``Kuranishi structure''.
\\[1ex]
2)\ \ Another aspect of this notion is that in a sense $\tau$ connects the
a priori unrelated differentiable structures on the various distinguished
local trivializations coming from the Fredholm structure. In particular
it will be possible to introduce index and orientation of Kuranishi
sections.
\qed
\end{rem}

From the definition and Proposition~\ref{transverse_section} we
readily obtain an abstract transversality construction:
\begin{prop}\label{abstract_transversality}
  Let $s$ be an $S$-Kuranishi section of a Banach orbibundle $E$ over a
  Banach orbifold $X$ submerging onto $S$. Let $\tau: F\rightarrow E$
  represent the Kuranishi structure, $q:F \rightarrow X$ the bundle
  projection. Then the section
  \[
    \tilde s\ :=\ q^*s +\tau
  \]
  of $q^*E$ over the total space $F$ is transverse. In particular, $\tilde Z
  := Z(\tilde s)$ is a topological suborbifold of $F$ that submerges onto
  $S$ with finite dimensional fibers.
\end{prop}
\pf
In distinguished local trivializations $\Pi\circ\tilde s$ is uniformized by
the map
\[
  \Pi\circ \hat{\tilde s}:\hat U \times F_0 \longrightarrow E_0,\quad
  (\hat x,f)\longmapsto \hat s(x)+\hat\tau(x)\cdot f\,,
\]
$\hat U= W\times L\times V$. By assumption $\Pi\circ \hat{\tilde s}$ is
differentiable relative $W\times L$. The differential at $0\in \hat U$ is
\[
  T_0 V\times F_0\longrightarrow E_0,\quad
  (v,f)\longmapsto \sigma(v)+\tau(0)\cdot f\,,
\]
$\sigma= D_V(\Pi\circ \hat s)(0)$, and this map is surjective by the
Kuranishi-property of $\tau$.
\qed

\begin{rem}\label{tilde_Z}\rm
Let $s_\can$ be the tautological section of $q^*F\rightarrow F$. So on
local uniformizers $s_\can$ is given by putting
\[
  \Pi\circ\hat s_\can:= \prj_{F_0}: \hat U\times F_0
  \longrightarrow F_0\,.
\]
The zero locus of $s_\can$ is just the zero section of $F$, and
can be identified with $X$. Restricting to $\tilde Z= Z(\tilde s)$ we obtain
\[
  Z(s)=Z(s_\can|_{\tilde Z})\,.
\]
In this way we have exhibited the zero locus of $s$ as zero locus of a
{\em finite rank} orbibundle over a finite dimensional orbifold.
\qed
\end{rem}
Using the finite dimensional orbifold $\tilde Z$ from the last proposition
we can define index and orientation of Kuranishi sections. One more
notation: Let $\tilde Z\rightarrow S$ be a submersion of topological
orbifolds with finite dimensional fibers. Over a connected set in $\tilde
Z$ with fiber dimension $d$ over $S$ the {\em relative orientation
bundle} $\Or_{\tilde Z/S}$ of $\tilde Z$ over $S$ is the real orbi-line
bundle locally uniformized by the real line bundle associated to the
presheaf $(\hat U= W\times V) \mapsto H_d(V)$.
\begin{lemma}
  Let $s$ be an $S$-Kuranishi section of a Banach orbibundle $E$
  over a Banach orbifold $X$ submerging onto $S$ and $\tau:
  F\rightarrow E$, $\tau':F' \rightarrow E$ be compatible Kuranishi
  structures, $q:F \rightarrow X$, $q':F' \rightarrow X$ the bundle
  projections. Let $d$, $d'$ be the fiber dimensions of the topological
  orbifolds $\tilde Z =Z(q^* s+\tau)$ and ${\tilde Z}'= Z({q'}^*s +\tau')$
  at some $z\in Z(s)$.
  
  Then $d-\rk F = d' -\rk F'$ and there is a canonical isomorphism
  \[
    \Or_{\tilde Z/S}\otimes {\det F}^\vee
    \stackrel{\simeq}{\longrightarrow}
    \Or_{{\tilde Z}'/S}\otimes {\det F'}^\vee\,.
  \]
\vspace{-4ex}

\end{lemma}
\pf
Choose a distinguished local trivialization centered at $z\in Z(s)$ over a
local uniformizer of the form $\hat U= W\times L\times V$ with
$\sigma:= D_V(\Pi\circ\hat s)(0,0)$. Then by the implicit function
theorem
\[
  d\ =\ \ind\sigma+\dim L+\rk F,\quad
  \ind\sigma=\dim\ker\sigma -\dim\coker\sigma\,,
\]
and similarly for $d'$. Hence
\[
  d-\rk F\ =\ \ind\sigma+\dim L\ =\ d'-\rk F'\,.
\]
The fiber of $\Or_{\tilde Z/S}\otimes \det F^\vee$ at $z$ is locally
uniformized by
\[
  \det\Big(T_0 L\oplus \ker(\sigma+\tau(0))\Big)\otimes \det F_0^\vee
  \ =\ \det T_0 L\otimes \det \ker\sigma \otimes (\det \coker\sigma)^\vee\,,
\]
(similarly for $\Or_{{\tilde Z}'/S}$) and this isomorphism is canonical.
\qed

By the lemma the following definition makes sense:
\begin{defi}\rm\label{index_defi}
  Let $s$ be an $S$-Kuranishi section of a Banach orbibundle $E$ over a
  Banach orbifold $X$ submerging onto $S$. Let $\tau: F\rightarrow E$
  represent the Kuranishi structure, and $\tilde Z =Z(q^*s+\tau)$.\\[1ex]
  1)\ \ The {\em index} $\ind_z(s)$ of $s$ ({\em relative $S$}) at $z\in
  Z(s)$ is defined as difference of the fiber dimension of $\tilde
  Z\rightarrow S$ at $z$ and $\rk F$.\\[1ex]
  2)\ \ The {\em orientation bundle $\Or(s)$} of $s$ ({\em relative $S$}\,)
  is defined as the real orbi-line bundle $\Or_{\tilde Z/S} \otimes \det
  F^\vee$. We call $s$
  \begin{itemize}
  \item
    {\em locally orientable} if $\Or(s)$ is a real line bundle (i.e.\ iff the
    local groups act trivially on the relative orientation line)
  \item
    {\em orientable} if $\Or(s)$ is a {\em trivial} real line bundle
  \item
    {\em oriented} if a trivialization of $\Or(s)$ has been chosen.
  \end{itemize}
\vspace{-6ex}

\qed
\end{defi}


\subsection{Localized Euler classes on Banach orbifolds}
Here is the main result of this chapter:
\begin{theorem}\label{virfc}
  Let $X$ be a topological Banach orbifold, $p:X\rightarrow S$ a
  submersion onto a finite dimensional, oriented topological orbifold,
  $\dim S=b$, and $s$ an oriented S-Kuranishi section of a Banach
  orbibundle $E$ over $X$ of constant index $d$ (relative $S$).
  We assume $Z=Z(s)$ to be compact.

  Then there exists a \underline{localized Euler class} $[E,s]\in H_{d+b}(Z)$,
  depending only on $E$ and the Kuranishi section $s$, with the
  following properties:
  \begin{enumerate}
  \item
    If $s$ is transverse (so $Z$ is an oriented topological orbifold of
    dimension $d+b$) then $[E,s]=[Z]$.
  \item
    Let $(\bar S,T_{\bar S})\hookrightarrow(S,T_S)$ be a closed oriented
    topological suborbifold, Poincar\'e dual to $\alpha\in H_{\bar S}^*(S)$
    and $\bar X=p^{-1}(\bar S)$, $\bar E=E|_{\bar X}$, $\bar s=s|_{\bar X}$
    the induced orbifold, orbibundle and Kuranishi section (with the
    induced orientation). Then
    \[
      [\bar E,\bar s]\ =\ [E,s]\cap p^*\alpha\, .
    \]
   \end{enumerate}
\end{theorem}
\pf
Let $\tau:F\rightarrow E$ be a representative of the Kuranishi
structure on $s$. By application of
Proposition~\ref{abstract_transversality} we obtain a
$(d+b+\rk F)$-dimensional topological orbifold $\tilde Z\subset F$. Since the
Kuranishi structure and $S$ are oriented, $\tilde Z$ is oriented too.  Let
$q: \tilde Z \rightarrow X$ denote the restriction of the bundle
projection $F\rightarrow X$ and $s_\can\in \Gamma (q^* F)$ the
tautological section, taking the value $f$ at point $f$. Recall also the
Thom class $\Theta_{q^*F} \in H^{\rk F}_{\tilde Z} (q^*F)$ ($\tilde
Z\subset q^*F$ the zero section). Because $Z(s_\can)=Z$ we have
$s_\can^* (\Theta_{q^* F}) \in H^{\rk F}_Z(\tilde Z)$, so we are able to
define
\[
  [E,s]\ :=\ [\tilde Z]\cap s_\can^*(\Theta_{q^* F})\in H_{d+b}(Z)\, .
\]

To establish independence of the choice of $F$ and $\tau$ we first show
that if $\sigma$ was already surjective (i.e.\ $s$ was a transverse
section), in the case of which we could also take $F=0$, then $[E,s]=[Z]$
(this is statement~(1)). Both these classes are of top dimension, so it is
enough to check this on small open sets (for each connected component
of $Z$). In the application of the implicit function theorem that lead to
the local uniformizing system for $\tilde Z$ at $z$ we may now view
${\hat F}_z$ as additional parameter. Using the above notation with
distinguished uniformizers, $Z$ and $\tilde Z$ are locally uniformized by
$W\times K \times L$ and $W\times K \times L\times F_0$ respectively,
with
$K=
\ker\sigma$. So $[E,s]=[Z]$ reduces to the elementary identity
\[
  [W\times L\times K\times F_0]\cap\prj_4^*\delta
  =[W\times L\times K]\,,
\]
where $\delta\in H_{\{0\}}^r(F_0)$ is Poincar\'e dual to $[0]\in
H_0(\{0\})$.

In the general case, if $\tau': F'\rightarrow E$ is a Kuranishi structure
that is compatible with $\tau$, $q':F'\rightarrow X$, ${\tilde
s}'=(q')^*s+\tau'$, $\tilde Z'=Z(\tilde s')$, we consider the sum Kuranishi
structure $\tau+\tau': F\oplus F'\rightarrow E$. We write
$q_{FF'}:F\oplus F' \rightarrow X$ for the bundle projection and
$\tilde{\tilde s\,}=\prj_1^*\tilde s+\prj_2^*\tau' =q_{FF'}^*
s+\prj_1^*\tau+\prj_2^*\tau'$ for the induced section of
$q_{FF'}^*E$, $\prj_1:F\oplus F'\rightarrow F$, $\prj_2:F\oplus F'
\rightarrow F'$ the projections, and $\tilde{\tilde Z\,} =Z(\tilde{\tilde
s\,})$ for the finite dimensional suborbifold of $F\oplus F'$. Note that
${\tilde s}_\can:=\prj_1^* s_\can\oplus\prj_2^* s'_\can$ is the
tautological section of $q_{FF'}^*(F\oplus F')$, $\tilde
Z=\prj_1(Z(\prj_2^* s'_\can)\cap {\tilde{\tilde Z\,}})$ and
\[
  \Theta_{q_{FF'}^*(F\oplus F')}
  =(-1)^{\mbox{{\rm\scriptsize rk}}F\cdot\mbox{{\rm\scriptsize rk}}F'}
  \prj_1^*\Theta_{q^*F}\cup\prj_2^*\Theta_{{q'}^*F'}
  =\prj_2^*\Theta_{{q'}^*F'}\cup\prj_1^*\Theta_{q^*F}\, .
\]
Moreover, as shown in the local description above, we can view
$\tilde s$ as transverse $S$-Fredholm section of $q^*E$ on $F$. This case
has just been treated, so $[\tilde Z]=[\tilde{\tilde Z\,}]\cap (\prj_2^*
s'_\can)^*\Theta_{\mbox{{\rm\scriptsize pr}}_2^*(q')^*F}$ and
\begin{eqnarray*}
  \lefteqn{(q_{FF'})_*\Big[[\tilde{\tilde Z\,}]\cap
  {\tilde s}_\can^*\Theta_{q_{FF'}^*(F\oplus F')}\Big]}\hspace{1cm}&\\
  &=&q_*(\prj_1)_*\Big[[\tilde{\tilde Z\,}]\cap
  \Big((\prj_2^* s'_\can)^*\prj_2^*\Theta_{{q'}^*F'}
  \cup(\prj_1^* s_\can)^*\prj_1^*\Theta_{q^*F}\Big)\Big]\\
  &=&q_*\Big[[\tilde Z]\cap s_\can^*\Theta_{q^*F}\Big]\,,
\end{eqnarray*}
which is nothing but $[E,s]$ when computed with $\tau$. But the first
line is completely symmetric in $F,F'$. This proves well-definedness of
$[E,s]$.

To prove (2) we observe that $\tilde Z$ submerges onto $S$, and
$\bar\tau=\tau|_{\bar X}$ is the Kuranishi structure for $\bar s$
yielding $\widetilde{\bar Z\,} =\tilde Z\cap \bar F$, $\bar F=F|_{\bar X}$.
So $[\widetilde{\bar Z\,}]=[\tilde Z]\cap p^*\alpha$ and
\[
  [E,s]\cap p^*\alpha\ =\ 
  \Big([\tilde Z]\cap s_\can^*\Theta_{q^*F}\Big)\cap p^*\alpha
  \ =\ [\widetilde{\bar Z\,}]\cap{\bar s}_\can^*\Theta_{{\bar q}^*
  \bar F}\ =\ [\bar E,\bar s]\, .
\]
\vspace{-6ex}

\qed
			
\begin{rem}\rm
It is clear by construction that $[E,s]$ depends only on the {\em
germ} of $X$, $E$ and the Kuranishi section $s$ along $Z(s)$.
\qed
\end{rem}


\section{Prestable curves}
\subsection{Definition}\label{2.1}
A prestable curve is a number of Riemann surfaces joined at pairs of
points as in the following picture:\\[3ex]
\nopagebreak
\centerline{\epsffile{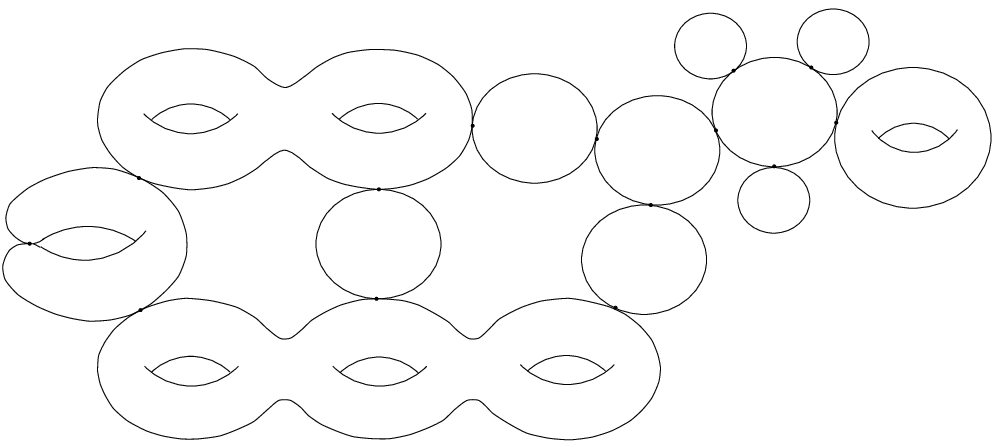}}\\[1ex]
{\bf Figure~2.1:} A prestable curve
\vspace{3ex}

\noindent
A precise definition is best done for our purposes in the language of
complex spaces (of course, an analogous definition can be made in
algebraic geometry). It is convenient to also incorporate additional
points.  Readers not so much acquainted with complex analysis should
consult \cite{grauertremmert} or \cite{kaupkaup} for standard notions of
complex space theory like reduction and normalization. For a complex
space $X$ we write $X_\sing$ and $X_\reg =X\setminus X_\sing$ for the
singular and regular locus respectively. Recall that a one-dimensional
ordinary double point (or node) is the germ of $z\cdot w=0$ in $\cz^2$ at
the origin.
\begin{defi}\rm
A {\em ($k$-marked) prestable curve} is a pair $(C,{\bf x})$ with
\begin{itemize}
\item
  $C$ is a reduced, compact, connected, one-dimensional complex space
  with at most ordinary double points as singularities
\item
  ${\bf x}=(x_1,\ldots,x_k)$ is a $k$-tuple of pairwise distinct points
  $x_\mu\in C_\reg$.
\qed
\end{itemize}
\end{defi}
Prestable curves also go by the names {\em noded Riemann
surfaces} or {\em Riemann surfaces with nodes}. While this is acceptable
terminology, I would not recommend the usage of ``cusp curves''
\cite{gromov}, although this captures well the behaviour of the natural
(Poincar\'e) metrics. The word ``cusp'' is reserved in algebraic geometry
for locally irreducible curve singularities ($y^2-x^3=0$ being the most
simple one), while a node is reducible with smooth irreducible
components.

There is the notion of {\em normalization} of reduced complex spaces.
Applied to a prestable curve $C$ the normalization is a morphism
\[
 \sigma:\hat C=\coprod_{\nu=1}^d D_\nu \longrightarrow C
\]
with the following properties:
\begin{itemize}
\item
  $D_\nu$ are the desingularizations of the irreducible components of
  $C$, hence compact Riemann surfaces
\item
  $\sigma$ is an isomorphism away from $C_\sing$ and two to one over
  $C_\sing$.
\end{itemize}
So normalization provides a decomposition of any prestable curve
$C$ into a union of Riemann surfaces $D_\nu$. The $D_\nu$ come with the
additional datum of finitely many unordered pairs of pairwise distinct
points, the preimages of the points in $C_\sing$.

Conversely, given Riemann surfaces $D_\nu$, a set of pairs of points
$\{y_i, y_{l+i}\}$, $i=1,\ldots,l$, defines a reduced analytic
equivalence relation that identifies points $y_i$ and $y_{l+i}$. The
quotient as ringed space (cf.\ e.g.\ \cite[\S49A]{kaupkaup}) is easily seen
to be a reduced complex space with ordinary double points exactly at
the images of $y_i \sim  y_{l+i}$. Provided there are enough points
identified to make $\coprod D_\nu/ (y_i \sim y_{l+i})$ connected, we
obtain a prestable curve. And by normalization, any prestable curve
arises in this way. This makes precise the intuitive idea of building up a
prestable curve from Riemann surfaces by joining pairs of points.
\vspace{1ex}

For a prestable curve $C$ the dimension of the first coherent
cohomology of the structure sheaf $h^1(C,\calo_C)$ is called
{\em (arithmetic) genus} $g(C)$ of $C$. A more geometric
definition runs as follows: We will see below that any prestable curve
can be deformed to a connected, smooth curve by deforming
neighbourhoods of the nodes $zw=0$ to cylinders $zw=t$ for some
constant $t\in\cz$. Let $C_t$ be such a deformation. Since by
connectedness $h^0(C,\calo_C) =1$ and $\chi(C,\calo_C)= h^0(C,\calo_C)-
h^1(C,\calo_C)$ does not change in flat families, we see
\[
  h^1(C,\calo_C)\ =\ h^1(C_t,\calo_{C_t})\,.
\]
And by Hodge theory (which does not work naively on $C$!)
$h^1(C_t,\calo_{C_t})$ is the ``number of holes'' of $C_t$, which equals
the ``number of holes'' of $C$. So the genus of the prestable curve in
Figure~2.1 is 9. Note that under normalization the genus of an irreducible
component may decrease. For instance, in Figure~2.1 the normalization
of the nodal toric component is the Riemann sphere.
\vspace{1ex}

\begin{defi}\rm
  The {\em automorphism group} $\Aut(C,{\bf x})$ of a prestable curve
  $(C,{\bf x})$ is the group of biholomorphism $\Psi:C \rightarrow C$ with
  $\Psi(x_\mu) = x_\mu$ for all $\mu$.
\qed
\end{defi}
The connected component of the identity $\Aut^0 (C,{\bf x}) \subset \Aut
(C,{\bf x})$ is a normal subgroup  that can be described
explicitely: Any $\Psi \in \Aut^0 (C,{\bf x})$ respects irreducible
components and special (i.e.\ singular or marked) points. Consider
the preimages of special points on $C$ as marked points on the
normalization. The pointed Riemann surfaces with higher dimensional
automorphism groups are
\begin{itemize}
\item
  $\pr^1$ with no marked point
\item
  $\pr^1$ with one marked point: $\Aut^0(\pr^1,\{\infty\}) =
  \cz\rtimes\cz^*$, where $\cz$ acts additively on
  $\pr^1\setminus\{\infty\} \simeq\cz$ and $\cz^*$ acts
  multiplicatively
\item
  $\pr^1$ with two marked points: $\Aut^0(\pr^1,\{0,\infty\}) =\cz^*$
\item
  an elliptic curve without marked points.
\end{itemize}
By connectedness of prestable curves a rational or elliptic curve
without a special point coming from a node can only occur as
single irreducible component of an irreducible prestable curve. Similarly
for $\pr^1$ with two special points that belong to only one node (a
nodal elliptic curve without added points). The remaining cases are:
\begin{itemize}
\item
  a rational component with one special point coming from a node,
  and at most one more marked point
\item
  a rational component with two special points coming from nodes
  and no marked points.
\end{itemize}
We will see in the discussion of deformation theory below (Remark~\ref
{prestable_deform}) that the second case shows different
behaviour depending on whether the irreducible component
$D_\nu$ disconnects $C$ or not, i.e.\ if $C\setminus
D_\nu$ consists of two or one connected component.

Automorphisms not connected to the identity can be described as
follows: The factor group $\Aut(C, {\bf x}) /\Aut^0(C, {\bf x})$ acts on
the irreducible components of $C$ by permutation. And the subgroup
respecting the irreducible components is a direct product of subgroups of
the automorphism groups of the irreducible components.
\vspace{1ex}

Particularly nice are prestable curves without infinitesimal
automorphisms:
\begin{defi}\rm
  A prestable curve $(C,{\bf x})$ is called {\em stable} if $\Aut^0(C,{\bf
  x})$ is trivial, which is if and only if $\Aut(C,{\bf x})$ is finite. An
  irreducible component $D_\nu\subset C$ is called {\em stable
  component} if $\Aut^0(C,{\bf x})$ acts trivially on $D_\nu$, else {\em
  unstable}.
\qed
\end{defi}
The equivalence statement in the definition follows from the finiteness
of automorphism groups of Riemann surfaces without infinitesimal
automorphisms, plus the finiteness of the number of irreducible
components of $C$.

To any $k$-marked prestable curve $(C,{\bf x})$ with $2g(C)+k\ge3$ is
naturally associated a stable $k$-marked curve $(C,{\bf x})^\st$ (the {\em
stabilization} of $(C,{\bf x})$), won by successive contraction of
nonstable components. The stabilization of the prestable curve in
Figure~2.1 looks as follows:\\[3ex]
\nopagebreak
\centerline{\epsffile{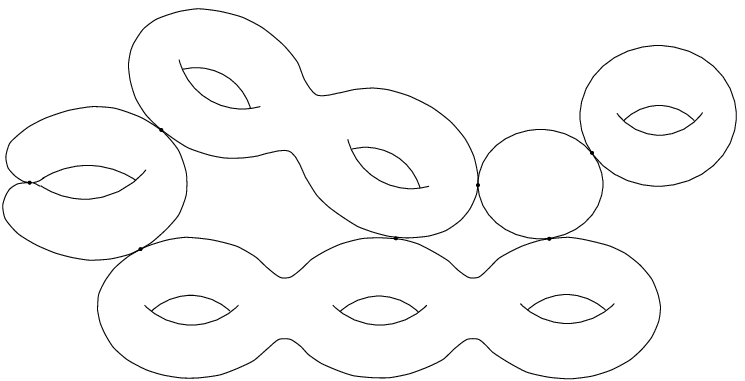}}\\[1ex]
{\bf Figure~2.2:} Stabilization of the curve in Figure~2.1
\vspace{3ex}

\noindent
We see that an originally stable component can become
unstable after contraction of some other
components, but the process stops after finitely many steps by
properness of $C$. We refer to components $D\subset C$ that are
contracted under the stabilization map $(C,{\bf x})\rightarrow(C,{\bf
x})^\st$ as {\em bubble components}. So an unstable component is a
bubble component but not in general conversely. An example is the
rational component with five nodal points in Figure~2.1.

Another, albeit non-canonical, method to make a prestable curve stable
without changing its topology is to add marked points. There
are at least one (two) point(s) needed on each rational component with
only two (one) special points. We refer to this process as {\em
rigidification}. Rigidification is the method of choice for the study of
deformation theory of prestable curves once the deformation theory of
stable curves is understood.


\subsection{Deformation theory of prestable curves}
\label{section_deformation}
Let $(C, {\bf x})$ be a prestable curve. A {\em versal deformation} of
$(C, {\bf x})$ is the germ (at $s_0\in S$) of a proper, flat family $q:\calc
\rightarrow S$  of complex spaces with an identification $C\simeq
q^{-1}(s_0)$, together with $k$ sections $\underline{\bf x}: S \rightarrow
\calc \times_S\ldots \times_S \calc$ with $\underline{\bf x}(0) ={\bf x}$,
having the universal property for germs of deformations of
$(C, {\bf x})$. The latter means that given another such family $p: \calc'
\rightarrow S'$, $C\simeq p^{-1}(s'_0)$ with $k$ sections $\underline{\bf
x'}$, then possibly after shrinking $S'$ to a
smaller neighbourhood of $s'_0$ there is a cartesian diagram
\[\begin{array}{ccc}
  \calc' &\stackrel{F}{\llongrightarrow}& \calc\\[10pt]
  \diagl{q'}&&\diagr{q}\\[10pt]
  S'&\stackrel{f}{\llongrightarrow}& S
\end{array}\]
inducing the identity on the central fibers and which is compatible with
the sections $\underline{\bf x}$ and $\underline{\bf x}'$. The diagram
being cartesian means that $(q',F): \calc'\rightarrow S'
\times_S \calc$ is an isomorphism. If the germs of the morphisms $(f,F)$
over $s_0$ are uniquely determined $q$ is called {\em universal}; if only
the map between the tangent spaces at $s'_0$ and $s_0$ is unique,
$q$ is called {\em semiuniversal}. Similar notions can be formulated for
deformations of complex spaces, of morphisms of complex spaces, of germs
of complex spaces etc.

Now (pointed) compact Riemann surfaces without infinitesimal
automorphisms have universal deformations (with smooth base spaces $S$).
The following methods are available to prove this: Teichm\"uller theory,
geometric invariant theory \cite{mumford}, Banach-analytic methods
\cite{douady}, \cite{grauert}, \cite{palamodov}, a power
series method \cite{forsterknorr}. More generally, stable marked curves
possess universal deformations, while prestable curves that are not stable
have only semiuniversal deformations \cite{delignemumford}
\cite{knudson}. While all these facts are
well-established in algebraic geometry and complex analysis, it
might be worthwile to give an elementary construction of universal
deformation spaces of stable curves assuming only the widely known
deformation theory of Riemann surfaces. Our proof of universality
will however rely on the machinery of formal deformation theory and might
thus be hard to accept for people without proper background. Later
we will also need a statement that is a by-product of the Banach-analytic
proofs in op.\ cit.\ (Proposition~\ref{analytic-dilation_comp}). With some
more effort one could prove both universality of our construction
and this proposition by employing Banach-analytic methods
directly, but we will not do this here.
\vspace{1ex}

Deformations of nodal curves can be effectively decomposed into
deformations of the nodes and deformations of the
(normalizations of the) irreducible components (with added points). We
begin with deformations of the nodes. An obvious deformation of
$Z_0=(zw=0)$ in the product of unit disks $\Delta^2\subset \cz^2$
is given by
\[
  Z_t\ =\ \{(z,w)\in \Delta^2\mid z\cdot w=t\}\,.
\]
Topologically this is a family of cylinders, degenerating to two disks
joined at the origin: Explicitely, write $t=\tau^2 e^{i\theta}\in\Delta$.
By projection onto the coordinate axes, for $t\neq0$, $Z_t$ can also be
described (even complex analytically) as two semi-closed
annuli $A_\tau^+, A_\tau^- =\Delta\setminus \Delta_\tau$
glued along the inner circles $S_\tau^\pm$. Here $\Delta_\tau$
is the open disc of radius $\tau$. The gluing is such that the
circles $S_\tau^+$ and $S_\tau^-$ are shifted against each
other by the angle $\theta$ (Figure~2.3).
Thus circling around the origin in the parametrizing disc
gives rise to a Dehn twist, the family is nontrivial over any
pointed disk $\Delta_r^*$. We write $Z_0= \Delta^+
\amalg_{\{0\}}\Delta^-$, and generally $O\amalg P$ for any objects
defined by $O$ on $\Delta^+$, $P$ on $\Delta^-$.\\[3ex]
\nopagebreak
\centerline{\epsffile{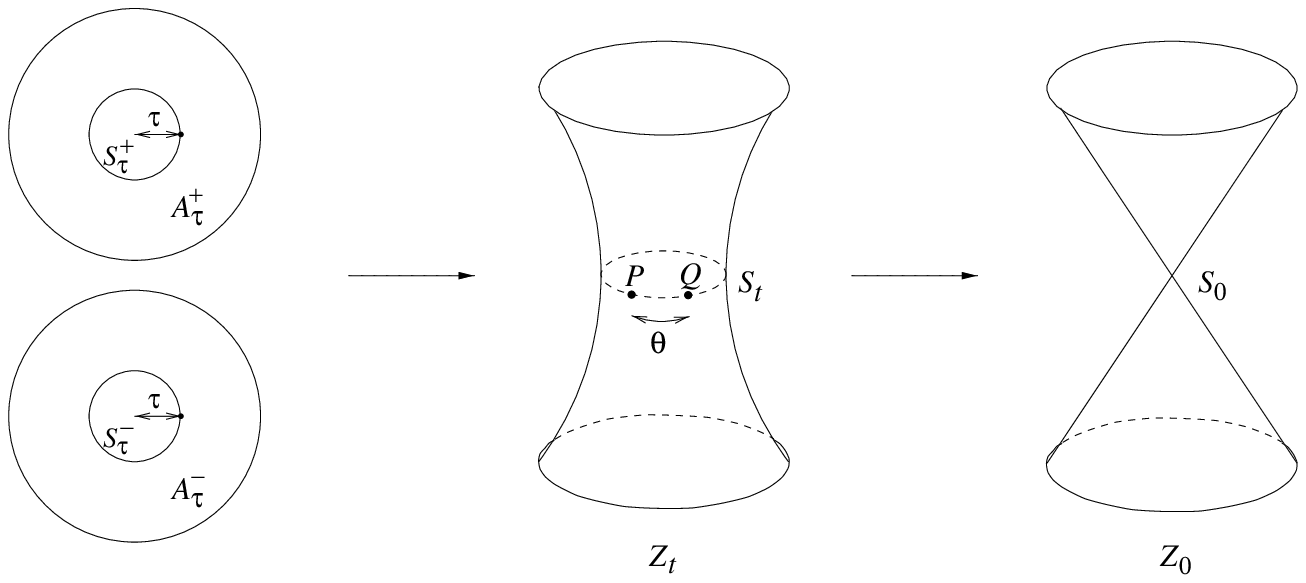}}\\[1ex]
{\bf Fig.~2.3:} Joining the annuli $A_\tau^+$, $A_\tau^-$.
\vspace{3ex}

Now let $(C,{\bf x})$ be a stable $k$-marked curve. To construct a
universal deformation of $(C,{\bf x})$ we consider its normalization
$\hat C= \coprod_\nu D_\nu$. Let $\{y_i,y_{l+i}\}_{ i=1,\ldots,l}$ be
the preimages of the nodes. Write $\hat{\bf x}$ for the lift of ${\bf x}$
to $\hat C$, ${\bf y}= (y_1,\ldots, y_{2l})$ and $\hat{\bf x} \vee {\bf y} =
(\hat x_1,\ldots, \hat x_k, y_1,\ldots, y_{2l})$ for the concatenation of
$\hat{\bf x}$ and ${\bf y}$. Then $(\hat C,\hat{\bf x}\vee {\bf y})$ is a
disjoint union of stable (!) marked Riemann surfaces $(D_\nu,{\bf
x}^\nu)$. As mentioned before, each of these has a universal deformation
$(q_\nu:\calc_\nu\rightarrow S_\nu, \underline{\bf x}^\nu)$. The product
of these deformations
\[
  \hat q:= \prod_\nu q_\nu: \hat \calc:= \prod \calc_\nu\rightarrow
  \hat S:=\prod S_\nu\,,
\]
together with the deformation $\underline{ \hat{\bf x}}\vee \underline{\bf y}$
of $\hat{\bf x}\vee {\bf y}$ obtained from the deformations $\underline{
\bf x}^\nu$ of ${\bf x}^\nu$, is obviously a universal deformation of $(\hat C,
\hat{\bf x} \vee {\bf y})$ (this curve need not be connected, but
the above definitions make of course sense for these slightly more
general curves).

To bring in deformations of the $l$ nodes, let $z_i$ be
local holomorphic functions defined near $\img (\underline{y_i})$, restricting to
local holomorphic coordinates on the fibers of $\hat q$ and with
$z_i (\underline{y_i}) \equiv0$. Possibly by rescaling of $z_i$, and shrinking of
$\hat S$ we may assume that each of the sets
\[
  V_i\ :=\ \{|z_i| <1\}\,,\quad i=1,\ldots,2l\,,
\]
is isomorphic to $\hat S\times\Delta$, contains no marked points
except the deformations of the $y_i$, and that the
$z_i$ are defined in neighbourhoods of the closures of these
sets. Let $t_{l+i}=t_i: \hat\calc \times \Delta^l  \rightarrow \Delta$ be the
projection to the $i$-th factor of $\Delta^l$. Define the total space
$\calc$ of the deformation by removing the closed sets
\[
  \{(z,{\bf t})\in V_i\times \Delta^l\mid |z_i(z)| \le |t_i|\}\,,
\]
from $\hat \calc \times\Delta^l$ and identifying the rest of
$V_i\times \Delta^l$ and $V_{l+i} \times\Delta^l$ via
\[
  z_i\cdot z_{l+i}\ =\ t_i\,.
\]
We have thus joined neighbourhoods of the deformations of $y_i$ and
$y_{l+i}$ by the standard deformation $Z_t$ of the node $Z_0$ with
deformation parameter $t=t_i$. The gluing has ben done in such a way
that for fixed $s, t$ the $2l$ circles $|z_i|=1$ and $|z_i|=|t_i|$ are identified with
$|z_{l+i}|=|t_i|$ respectively $|z_{l+i}|=1$. The connected components of the
resulting prestable curve minus these circles consist of a (not necessarily
connected) open Riemann surface with $2l$ ideal boundary curves and $l$
possibly degenerate cylinders $Z_{t_i(s)}$.

Let us now define $q: \calc \rightarrow S:=\hat S\times \Delta^l$ as the
morphism induced by $\hat q\times\id_{\Delta^l}$ and $\underline{\bf
x}: S\rightarrow \calc^k$ as induced by the first $k$ components of
$\underline{ \hat{\bf x}}\vee \underline{\bf y}$.
\begin{prop}\label{calc_univ}
  $(q:\calc\rightarrow S, \underline{\bf x})$ is a universal deformation
  of $(C,{\bf x})$. 
\end{prop}
\pf
The tangent space at $(C,{\bf x})$ to the functor of families of marked
prestable curves is given by the space of isomorphism classes of such
families over the double point $(x^2=0)= \Spec\, \cz[x]/(x^2)$ (a
zero-dimensional complex space or scheme) with central fiber $(C,{\bf x})$.
It is a standard fact (cf.\ e.g.\ \cite{flenner}, \cite{ran}) that this space can
be identified with
\[
  \Ext^1_C(\Omega_C(|{\bf x}|),\calo_C),
\]
where we write $|{\bf x}|$ for the reduced divisor underlying the tuple $\bf x$.
The beginning of the spectral sequence for the composition of the derived
functors of $\Gamma$ and $\Hom$ on $C$ is
\[
  0\longrightarrow H^1(\Omega_C(|{\bf x}|)^\vee)
  \longrightarrow \Ext^1_C(\Omega_C(|{\bf x}|),\calo_C)
  \longrightarrow H^0(\Ext_{\calo_C}^1( \Omega_C(|{\bf x}|),\calo_C)
  \longrightarrow 0\,.
\]
A direct local computation yields
\[
  \Hom_C(\Omega_C(|{\bf x}|),\calo_C)\ =\ \Theta_C(-|{\bf x}|-C_\sing)\,,
\]
that is, the sheaf of holomorphic vector fields on $C$ with zeros at the
special points. By the Leray spectral sequence $H^1$ of this sheaf can be
computed by pulling back to the normalization $\hat C$ of $C$. So the term
on the left of the short exact sequence is
\[
 H^1(\hat C,\Theta_{\hat C}(-|\hat{\bf x}|-|{\bf y}|))\
 =\ \Ext^1_{\hat C}(\Omega_{\hat C}(|\hat{\bf x}|+|{\bf y}|),\calo_{\hat C})\,,
\]
where as above $\hat{\bf x}$, $\bf y$ are the
preimages of $\bf x$ on $C_\sing$ on $\hat C$ respectively. By the same
reasoning as above this is the tangent space to deformations of $(\hat C,\hat{\bf
x}\vee {\bf y})$.

On the other hand, $\Omega_C(|{\bf x}|)$ is locally free away from
$C_\sing$, while at a singular point $0\in(zw=0)$ the stalk is
$\calo_\cz (dz,dw)/ (wdz+zdw)$. So $\Ext^1_{\calo_C}
(\Omega_C(|{\bf x}|), \calo_C)$ is a skyscraper sheaf with one-di\-men\-sional
stalks concentrated at the singular points. Intrinsically the stalks are of
the form $\Ext^1_A (\Omega_A,A)$ where $A=\cz\{z,w\}/(zw)$ is a local
analytic algebra of complete intersection type. For such algebras the first
$\Ext$-group is the tangent space to deformations. So the right-hand term in
the exact sequence is isomorphic to $\cz^l$ with each factor equal
to the tangent space of deformations of one of the $l$ nodes of $C$.

It follows from the spectral sequence for tangent cohomology that this
sequence is indeed compatible with deformation theory:
Given a deformation of $(C,{\bf x})$, the Kodaira-Spencer map to
$\Ext^1(\Omega(|{\bf x}|),\calo)$ composed with the map to
$H^0(\Ext^1(\Omega(|{\bf x}|),\calo)$ is the Kodaira-Spencer map for the
induced deformation of the local analytic
algebras $\calo_{C,y}\simeq \cz\{z,w\}/(zw)$. This composition vanishes iff the
singularity remains analytically unchanged up to first order. The
infinitesimal deformation is then by (first order-) changing the analytic
patching data away from the singularities.

In our case the base of the deformation is $S=\hat S\times \Delta^l$. By
construction the second factor generates the tangent space to deformations
of the nodes. And restricted to $\hat S \times \{0\}$, the lift of the
Kodaira-Spencer map to $H^1(\Theta_{\hat C}(-|\hat{\bf x}|-|{\bf y}|))$ is an
isomorphism by universality of the family $(\hat\calc \rightarrow \hat S,\hat{\bf
x}\vee{\bf y})$ of marked Riemann surfaces. This shows that the
Kodaira-Spencer map
\[
  T_{s_0}S \rightarrow \Ext^1(\Omega_C(|{\bf x}|),\calo_C)
\]
is an isomorphism. Under the hypothesis of existence of a universal
deformation space this is nothing but the differential of the map from $S$
to the base space of the universal deformation. So this map is a local
isomorphism and we have proven universality of our family.
\qed
\begin{rem}\rm
Alternatively, the fact that the obstruction space $\Ext^2(\Omega_C (|{\bf x}|),
\calo_C)$ vanishes identically implies that the base space of the formal
(uni-) versal deformation space (\cite{schlessinger}) is formally smooth.
In this case bijectivity of the Kodaira-Spencer map shows that our family is formally
universal. And by \cite[Satz7.1]{bingener} the latter implies analytic universality.
\qed
\end{rem}
\begin{rem}\rm\label{prestable_deform}
While we will not explicitely need this it is instructive to also understand
the situation for prestable curves with unstable components. First we remark
that rigidification (introduced at the end of Section~\ref{2.1}) works locally
also in families; that is, given a family of prestable curves
$(q':\calc'\rightarrow S', \underline{{\bf x}'})$ there is a tuple $\underline{\bf
y}$ of local holomorphic sections of $q$ making $(q':\calc'\rightarrow S',
\underline{{\bf x}'} \vee \underline{\bf y})$ into a family of {\em stable} curves.
So letting $(q: \calc\rightarrow S, \underline{\bf x} \vee \underline{\bf y})$
be a universal deformation of a rigidification $(C,{\bf x}\vee {\bf y})$ of a
prestable curve $(C,{\bf x})$, the family $(q:\calc \rightarrow S,
\underline{\bf x})$ is in fact a {\em versal deformation} of $(C,{\bf x})$.
But since we dropped the information about the added points there will
be higher-dimensional subsets of $S$ parametrizing a given isomorphism
class of prestable curves. The largest occurring dimension equals the
number of points added. The equivalence relation on $S$ induced by
isomorphism classes is not closed, but the closure of each equivalence
class is analytic. In fact, there is the germ of a holomorphic action of the
non-compact complex Lie group $\Aut^0 (C,{\bf x})$ on $S$ stabilizing the
center $s_0\in S$, that is induced from an action on $\calc$, and that extends
the action on the central fiber. Being a germ of an action means that the map
$\Aut^0 (C,{\bf x}) \times S\rightarrow S$ is only defined in a neighbourhood
of $\{\id\} \times S$. For finite automorphism groups such action will be
constructed in the next section. The equivalence relation of isomorphism
classes is induced by this (non-proper) action.

The easiest example is probably two copies of $\pr^1$ with two marked
points on one component and one marked point on the other. Adding one
more point gives rise to only one module, the deformation parameter $t$ of
the node. Dropping the additional point makes all curves $C_t$ isomorphic
for $t\neq 0$. And indeed, in appropriate coordinates, the action of
$\Aut^0 (C,{\bf x}) =\cz^*$ is just the (germ of the) linear action on
$\Delta\subset \cz$. The topological quotient consists of two points and
is non-Hausdorff.
\vspace{1ex}

For the study of properties of GW-invariants it is important to relate the
semiuniversal deformation of $(C,{\bf x})$ to the universal deformation of
$(C,{\bf x})^\st$. Let $S_\st$ be the base space of the latter. We claim that
stabilization induces a holomorphic surjection $S\rightarrow S_\st$ of
complete intersection type. To this end we observe that the curve $C_\st$
underlying $(C,{\bf x})^\st$ is a closed subcurve of $C$. The closure of
$C\setminus C_\st$ in $C$ is the union of all bubble components. Its
connected components intersect $C$ in at most two points. Let us call such
a connected component {\em bubble tree} if it intersects $C$ in exactly one
point and {\em bubble chain} in case of two points. In terms of the contraction
map $C\rightarrow C_\st$ a bubble tree lies over a smooth point of $C_\st$
while a bubble chain maps to a node. In our analytic description the
deformation parameters belonging to bubble chains and bubble trees split up
into deformation parameters of nodes and possibly some moduli of stable
bubble components (the position of nodes and at most one more marked
point on the bubble component).

Now since a bubble tree contracts to either a non-special point of $C_\st$ or
a marked point, only the deformation parameter of a possible marked point
survives stabilization. So on these components of the base spaces
$S\rightarrow S_\st$ is just a linear projection.

For a bubble chain there is a chain of rational components $D_1,\ldots, D_s$
connecting the two preimages $y,y'\in C$ of the node, i.e.\ $y\in D_1$, $y'\in
D_s$, $D_\nu\cap D_{\nu+1}\neq\emptyset$. We take holomorphic coordinates
$z_\nu$ on the bubble components to be $z$ or $z^{-1}$ at $0,\infty \in\pr^1=
\cz\cup\{\infty\}$. Let $t_\nu$ (respectively, $t_0$, $t_s$) be the deformation
parameters at $D_\nu\cap D_{\nu+1}$ (respectively, at $y$, $y'$) and $t$ the
deformation parameter of the node in $C_\st$. The relevant part of the
stabilization map $S\rightarrow S_\st$ is then of the form
\[
  (t_0,\ldots,t_s)\ \longmapsto\ t=t_0\cdot\ldots\cdot t_s\,.
\]
This can be checked by explicitely working out an isomorphism of the
deformed bubble chain with $Z_{t_0\cdot\ldots\cdot t_s}$, for any $t_0, \ldots,
t_s\neq 0$. Possible further deformation parameters of the bubble chain are
just dropped under the map. So this map is not smooth in the sense of
analytic geometry (i.e.\ a linear projection in appropriate coordinates) but still
of complete intersection type.
\qed
\end{rem}


\subsection{Automorphisms of universal deformations and $\calm_{g,k}$}
It will also be important for us to understand the role of the
automorphism group in deformation theory.
\begin{lemma}
  Let $(q: \calc\rightarrow S, \underline{\bf x})$ be a universal
  deformation of a stable curve $(C,{\bf x})$. Then possibly after
  shrinking $S$ the action of $\Aut (C,{\bf x})$ on the central fiber
  extends to an action on $\calc$ that descends to $S$ (thus making
  $q$\ \ $\Aut(C,{\bf x})$-equivariant). Moreover, $s$, $s'\in S$ parametrize
  isomorphic marked prestable curves iff there exists $\Psi\in \Aut
  (C,{\bf x})$ with $s'=\Psi(s)$.
\end{lemma}
\pf
For any $\Psi\in \Aut(C,{\bf x})$ the isomorphism $\Psi: C\simeq C_0\subset
\calc$ can be extended by universality to a holomorphic map $\tilde \Psi:
\calc\rightarrow\calc$ (possibly after shrinking $S$). Again by universality
it holds $(\Psi\circ\Psi')\tilde{}= \tilde\Psi\circ\tilde\Psi'$, and
$\tilde\id =\id_\calc$, so this gives the desired action on $\calc$.

For the second statement we need the fact that given two families
of (marked) prestable curves $(p: X\rightarrow T, \underline{\bf x})$,
$(p': X'\rightarrow T, \underline{\bf x}')$ over the same parameter space
$T$, the set of base points $t$ such that $(p^{-1}(t),\underline {\bf x}(t))$
and $({p'}^{-1}(t),\underline {\bf x}'(t))$ are isomorphic, is analytic.
This follows from the representability of the {\em isomorphism functor}
by closed analytic subsets of the base. In the category of
schemes this is proven in \cite[Thm.1.11] {delignemumford}, but the proof
is almost literally valid in the analytic setting too. ("Locally analytic", which
is enough for our purposes, follows also by general deformation theoretic
arguments for any deformations of compact objects.)

Now let $R\subset S\times S$ be the (graph of the) equivalence relation
given by isomorphism of fibers over our universal
deformation. $R$\, is the set of isomorphic fibers of the two analytic
families of marked prestable curves obtained from the universal
object over $S$ by pull-back via the two projections $\prj_1, \prj_2:
S\times S\rightarrow S$. By what we have said above $R$ is thus
an {\em analytic} subset of $S\times S$. It also follows from
\cite[Thm.1.11] {delignemumford} that for any (germ at $(s_0,s_0)$ of an)
irreducible component $\tilde S\subset R$ there exists an isomorphism
between the pull-backs of the universal object over $S$ by
$p_1=\prj_1|_{\tilde S}$ and $p_2=\prj_2|_{\tilde S}$. We may compose
one of the $p_i$ with the action of some $\Psi\in \Aut(C,{\bf x})$ to
achieve the (natural) identification of the central fibers to be the identity.
Then we are dealing with isomorphic deformations of $(C,{\bf x})$, so
universality shows $p_1=p_2$. We conclude that up to the action of
the automorphism group, the germ of any irreducible component of $R$
lies in the diagonal. So $R=\{ (s,\Psi(s))\mid \Psi\in \Aut(C,{\bf x}) \}$,
the equivalence relation induced by the action of the automorphism
group on the base.
\qed
\vspace{1ex}

By the lemma $S/\Aut(C,{\bf x})$ locally parametrizes isomorphism classes of
stable marked curves. By such quotients the set $\calm_{g,k}$ of
isomorphism classes of $k$-marked prestable curves of genus $g$ is
endowed with a topology. It is shown e.g.\ in \cite[Thm.3.3.4]{abikoff} that
this space is compact Hausdorff. And given a universal deformation $(q':
\calc' \rightarrow S', \underline{\bf x'})$ of a fiber $(C',{\bf x'})$ of some
other universal deformation $(q: \calc \rightarrow S, \underline{\bf x})$ it
follows by construction that, possibly after shrinking $S'$, there is an
$\Aut(C',{\bf x'})$-equivariant holomorphic embedding $\calc'
\rightarrow \calc$ respecting the fiber structure. Hence:
\begin{prop}\label{mgk_orbifold}
  $\calm_{g,k}$ has naturally a structure of
  (complex-analytic) orbifold with local uniformizing systems of the form
  \[
    S\ \stackrel{/\Aut(C,{\bf x})}{\llongrightarrow}\ \calm_{g,k}\,.
  \]
Moreover, the morphism of orbifolds
\[
  q:\calm_{g,k+1}\ \longrightarrow \calm_{g,k}\,,\quad
  (C,{\bf x}\vee x_{k+1})\ \longmapsto\ (C,{\bf x})^\st
\]
that set-theoretically forgets the last point and then stabilizes, is
a universal curve over $\calm_{g,k}$.
\qed
\end{prop}
Note that if $2g+k\le2$ (that is, $(g,k)\in \{(0,0), (0,1), (0,2), (1,0) \}$) there
are no stable $k$-marked genus $g$ curves, $\calm_{g,k} =\emptyset$.
Being a universal curve means that $q$ represents the
functor of families of $k$-marked genus $g$ curves, i.e.\ that any such
family $p:X\rightarrow T$ is pull-back of the universal family under some
(orbifold-!) morphism $T\rightarrow S$.

That $\calm_{g,k+1}$ is the universal curve over $\calm_{g,k}$
might seem somewhat astonishing at first glance. There are in fact two
cases when omitting a point and stabilizing changes the topology
of the curve:
\pagebreak

\noindent
\centerline{\epsffile{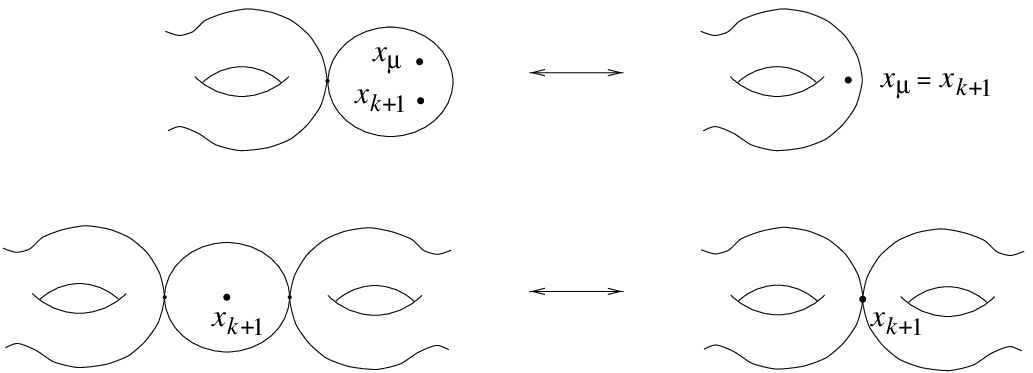}}\\[1ex]
{\bf Fig.~2.4: $(k+1)$-marked curves versus points on $k$-marked
curves} 
\vspace{3ex}

\noindent
However, the arrows pointing in both directions
indicates that this is a one-one
correspondence: The position of $x_{k+1}$ on the rational component
has no moduli. In terms of limiting positions of $x_{k+1}$ the first case
is when the additional point $x_{k+1}$ approaches some
$x_\mu$. The result is a rational component with two marked points
joined with the rest of the curve at the previous position of
$x_\mu$. The second case occurs when $x_{k+1}$ approaches a node.
Then the node is replaced by a rational component
sitting between the previous two irreducible components
and containing $x_{k+1}$. That the corresponding transformation of moduli
of $(k+1)$-marked curves and total spaces of $k$-marked curves is
indeed an isomorphism has been checked in \cite{knudson}.


\subsection{Continuous families and real resolutions}
\label{sect_cont_fam}
We will also deal with non-holomorphic families of prestable
curves. For families of Riemann surfaces $\hat q: \hat\calc \rightarrow
\hat S$ this should mean that locally there is a trivialization $\hat \calc
\simeq \Sigma \times \hat S$, $\Sigma$ an oriented differentiable
surface, such that the holomorphic structure on $\Sigma$ varies
continuously. For prestable curves we make the 
\begin{defi}\rm\label{cont_curves}\sloppy
A {\em continuous family of prestable curves} is a set of prestable
curves $\{C_s\}_{s\in S}$ and a continuous map $q:\calc\rightarrow S$
together with homeomorphisms $h_s: q^{-1}(s) \approx C_s$ for any
$s\in S$ and which locally on $S$ has the following form: For any $s_0 \in
S$ there exists an open neighbourhood $W\subset S$ and
\begin{itemize}
\item
  a continuous family $(\hat q: \hat\calc \rightarrow W,
  \underline{\bf y})$ of compact, $2l$-marked Riemann surfaces
\item
  continuous families of holomorphic functions $z_i$ defined
  on neighbourhoods of the closures of $V_i:= \{|z_i|<1\}$ with
  $z_i(\underline{y}_i) \equiv 0$, $i=1,\ldots,2l$; and continuous maps $t_i:
  W\rightarrow \Delta_{1/2}$, $i=1,\ldots,l$
\item
  letting $\calc' \rightarrow W$ be the family of prestable curves
  constructed from $\hat \calc \setminus  \{|z_i| \le |q^*t_i|\}$
  by gluing through
  \[
    z_i\cdot z_{l+i}\ =\ t_i\,,
  \]
  there exists a continuous family of biholomorphisms $\calc' \simeq
  \calc$ (relative $W$).
\end{itemize}
\vspace{-5ex}

\qed
\end{defi}
\fussy
In other words, $\{C_s\}$ is obtained by joining a continuous family of
compact Riemann surfaces through standard deformations of nodes
in a continuous fashion. Of course, marked points could be painlessly
incorporated in this definition, but we will not need this. If the fibers are
indeed $k$-marked {\em stable} curves, of genus $g$ say, one can
show that such families
are in one-one correspondence to continuous (orbi-) maps from $S$ to the
moduli orbifold $\calm_{g,k}$. So the notion of continuous family fits well
with the analytic discussion of prestable curves given above.

It will be convenient to think of the gluing also in another way: Let $U_i$,
$i=1,\ldots,l$, be the points of $q^{-1}(W)$ coming from $V_i \cup V_{l+i}$
and
\[
  U_0\ :=\ \{z\in q^{-1}(W)\mid z\not\in U_i \mbox{ or }|z_i(z)|> 1/2\
  \forall\,i\}\,.
\]
Then $\calu =\{U_i\}_{i=0,\ldots,d}$ is a covering of $q^{-1}(W)$ with
$U_0$ a continuous family of open Riemann surfaces and for $i>0$, $U_i
\cap C_s\simeq Z_{t_i(s)}$. The only nontrivial intersections among the
$U_i$ are
\[
  U_0\cap U_i\ =\ A_{1/2}\cup A_{1/2}\,,
\]
with $A_r= \{z\in\Delta\mid |z|> r\}$ the annulus with radii $r$ and 1.
We will also use the notation $U_i(s):= U_i\cap q^{-1}(s)$ for any $s\in S$.

The covering $\calu$ and $t_i$, $z_i$ are of course not unique, but locally
any two choices dominate a third one: We say $(\calu,t_i,z_i, z_{l+i})$,
$i=1,\ldots,l$, {\em dominate} $(\calu',t'_i,z'_i,z'_{l+i})$, $i=1,\ldots,l'$ if
$l'\le l$, $U'_i\subset U_i$, $i=1,\ldots,l$, $U'_i\subset U_0$,
$i=l+1,\ldots,l'$. $(z_i,z_{l+i}) \circ(z'_i, z'_{l'+i})^{-1}: Z_{t'_i(s)}
\rightarrow Z_{t_i(s)}$ are then open embeddings intertwining $(z_i,
z_{l+i})$ and $(z'_i, z'_{l'+i})$. A {\em morphism} from $\{C'_{s'}\}_{s'\in
S'}$ to $\{C_s\}_{s\in S}$ is a pair of continuous maps $S'\rightarrow S$,
$\calc'\rightarrow\calc$ commuting with $q': \calc' \rightarrow S'$ and
$q:\calc\rightarrow S$.

We will usually just write $q:\calc\rightarrow S$ for such families, the
holomorphic structure on the fibers being understood. By construction of our
analytically semiuniversal deformations of prestable curves $q:\calc \rightarrow S$
(Proposition~\ref{calc_univ} in conjunction with Remark~\ref{prestable_deform})
are instances of continuous families of prestable curves. And since any flat
analytic family of prestable curves is locally pull-back of such a family by a {\em
holomorphic} map of the base to $S$ the same is true for such families.
Conversely, by the Banach-analytic constructions of the deformation spaces,
any continuous family is locally pull-back (in an appropriate sense) of
$q:\calc \rightarrow S$ by a {\em continuous} map to $S$.
\vspace{1ex}

In a similar way one defines the notion of {\em (continuous families of)
holomorphic vector bundles} $\{E_s\}$ over $q:\calc\rightarrow S$ and
morphisms of such (require $E_s$ to be holomorphically trivialized over
$U_i$ etc.). The short-hand notation will be $E\downarrow \calc$.

Continuous families of prestable curves come with a distinguished
continuous family of holomorphic line bundles, the {\em relative dualizing
bundle} $\omega_{\calc/S}$. To define this note that the sheaf of
holomorphic 1-forms $\Omega_C$ of a prestable curve $C$ is not locally
free at the nodes. The stalk at $P\in C_\sing$, represented as $Z_0=
V(xy) \subset\cz^2$, is
\[
  \Omega_{C,P}\ =\ \calo_{C,P}(dx,dy)/(ydx+xdy)\,.
\]
To make this locally free one considers the double dual
$\omega_C:= \Omega_C^{\vee\vee}$, where we write $\calf^\vee=
\Hom_{\calo_C} (\calf,\calo_C)$ for any coherent sheaf $\calf$ on $C$.
Going over to $\omega_C$ has the effect of admitting simple poles at the
nodes. In fact,
\[
  \omega_{C,P}\ =\ \calo_{C,P}(x^{-1}dx)\ =\ \calo_{C,P} (y^{-1}dy)\,.
\]
In a family situation we obtain the relative dualizing bundle
$\omega_{\calc/S}$ (``dualizing'' because of its role in duality theory):
\begin{lemma}\label{dualizing_bundle}
  Let $q:\calc \rightarrow S$ be a continuous family of prestable
  curves. There exists a continuous family of holomorphic line bundles
  $\omega_{\calc/S} \downarrow \calc$ with sheaf of sections along $C_s$
  equal to $\omega_{C_s}$ for any $s\in S$.
\end{lemma}
\pf
The only thing not immediately clear is the behaviour on
families $U_i$, $i>0$, of possibly degenerating cylinders. These have the
form
\[
  V(zw-q^*t)\ \subset\ W\times\Delta^2
\]
for some continuous function $t:W\rightarrow \Delta_{1/2}$, $(z,w)$ the
linear coordinates on $\Delta^2$. Differentiating in the fiber
directions we get $zdw+wdz=0$. Thus $z^{-1}dz= -w^{-1}dw$ is a local
frame for $\omega_{\calc/S}$.
\qed
\vspace{2ex}

A central topic in this paper is the (local) comparison of maps with domains
varying in a continuous family of prestable curves. On continuous
families $U_i$, $i>0$, of possibly degenerate cylinders $Z_t$ we use the
description in terms of two annuli glued along their inner circles given
before Figure~2.3. We define a linear retraction
\[
  \lambda_\tau:\Delta\setminus\Delta_\tau\longrightarrow\Delta\,,\quad
  r e^{i\ph}\longmapsto\frac{r-\tau}{1-\tau}e^{i\ph}
\]
and declare a Lipschitz map $\kappa_t: Z_t\rightarrow Z_0$ by
$\lambda_\tau$ on each annulus. $\kappa_t$ is diffeomorphic outside the
central circle $S_t=\{(z,w)\in Z_t \mid |z|=|w|\}$ of $Z_t$ and contracts
$S_t$ onto the double point. The resulting retraction $\kappa: \bigcup_
{i>0} U_i \rightarrow \bigcup_ {i>0} U_i(s_0)$ can be extended
by a continuous family of diffeomorphisms $\kappa_0(s): U_0(s)
\rightarrow U_0(s_0)$ to yield a retraction of $q^{-1}(W)$ onto the
``central fiber'' $C_{s_0}$, also denoted
\[
  \kappa:\ q^{-1}(W)\ \longrightarrow C_{s_0}\,.
\]
For the purpose of regularity results we should control the differential of
$\kappa_s$ with respect to the natural holomorphic parameters provided by
uniformization. In fact, by further removing some families of cylinders from the
family of compact Riemann surfaces $\hat q:\hat\calc\rightarrow W$ (in a
continuous fashion) and gluing in holomorphic disks (after enlarging $l$) we may
assume each irreducible component of $\hat\calc$ to be of genus zero. Such
handlebody decompositions of continuous families of Riemann surfaces
into genus zero pieces (in fact, the unit disk minus one or two holes) and gluing
data is well-known in Teichm\"uller theory and gives rise to Fenchel-Nielsen
coordinates on Teichm\"uller space. We refer to \cite{abikoff} for a detailed
discussion.

Following \cite[\S15.7]{conway} we call an open set in $\cz$ all of whose
boundary components are round circles {\em circular region}. A circular region
is a disjoint union of a number of open Riemann surfaces of genus zero. By
uniformization, $U_0\rightarrow W$ can be realized as a family of circular regions
with boundaries varying continuously with $w\in W$.
In other words there is a fiberwise holomorphic homeomorphism
between the total spaces of these families. The pull-back of the linear
coordinate on the unit disk to $U_0$ will be called {\em natural holomorphic
parameter}.

Now let us adapt $\kappa$ to natural holomorphic parameters.
By extending the fiberwise linear retractions of the form $\kappa_t$ on
$U_i$, $i>0$, first linearly to an $\eps$-neighbourhood of $U_i\cap U_0\subset
U_0$, then to a Lipschitz (with respect to natural holomorphic parameters)
homeomorphism, and finally smoothing out on $U_0$ we get:
\begin{lemma}\label{retraction}
  Given a continuous family of prestable curves $q:\calc\rightarrow W$ as in
  Definition~\ref{cont_curves} there is a retraction $\kappa: \calc\rightarrow
  \calc_{s_0}$ to the central fiber which
  \begin{itemize}
  \item
    has continuous differentials relative $W$ with respect to natural
    holomorphic parameters
  \item
    on $U_i$, $i>0$, is fiberwise of the form $\kappa_{t_i(w)}: U_i(w)=Z_{t_i(w)}
    \rightarrow Z_0$.
  \end{itemize}
  \vspace{-4ex}
  
  \qed
\end{lemma}
Similarly, if $E\downarrow\calc$ is a holomorphic vector bundle over
$\calc$, one may also lift $\kappa_s$ to retractions $\tilde\kappa_s: E_s
\rightarrow E_{s_0}=: E_0$. We call retractions $\kappa$, $\tilde\kappa$
with the properties in the lemma {\em admissable}, as retractions of
continuous families of (vector bundles over) prestable curves.
\begin{rem}\rm
We will need this result only for the case of universal deformations. In these cases
we have smooth total spaces and such a retraction could easier be constructed by
gradient flows of $P\mapsto |q(P)|^2$ with respect to certain metrics on the total
space with singularities along the singular locus of $q$. The natural domain for our
results on such families is however the continuous world rather than the
holomorphic one.
\qed
\end{rem}
The restriction $\kappa_s$ of $\kappa$ to some $C_s$ is globally
Lipschitz, contracts a number of circles $S_t\subset \bigcup_{i>0} U_i(s)$
and is a diffeomorphism away from these circles. It decomposes into
cutting $C_s$ along a number of smooth circles and a real
(oriented) blow up of the irreducible components of $C$ at certain singular points
of $C$ (cf.\ \cite{looij}). We introduce the following terminology:
\begin{defi}\rm\label{partres}
Let ${\tilde C}$, $C$ be prestable curves. A {\em partial real resolution} of
$C$  is a map $\kappa:{\tilde C}\rightarrow C$ such that there is an open
neighbourhood $U$ of finitely many embedded circles $S_{(j)}\subset
{\tilde C}\setminus {\tilde C}_\sing$ with $\kappa|_{\tilde
C\setminus\cup S_{(j)}}$ a diffeomorphism onto the image (i.e.\ when
pulled back to the normalization), and such that for any $j$ the restriction
of $\kappa$ to the connected component of $U$ containing $S_{(j)}$ is
biholomorphic to one of the maps $\kappa_t: Z_t
\rightarrow Z_0$ defined above. If ${\tilde C}$ is nonsingular we speak
of a {\em real resolution} of $C$.
\qed
\end{defi}


\subsection{Some lemmata from Teichm\"uller theory}
At some places below we will need a few facts and notions from
Teichm\"uller theory that we collect in this section for the readers
convenience.

Recall that if $\kappa:\Sigma\rightarrow \Sigma'$ is an orientation
preserving diffeomorphism between (possibly non-compact) Riemann
surfaces, the {\em maximal dilation} of $\kappa$ is defined as
\[
  K(\kappa)\ :=\ \sup_{z\in\Sigma} \frac{|\di\kappa(z)|+|\dbar\kappa(z)|}
  {|\di\kappa(z)|-|\dbar\kappa(z)|}\,.
\]
The derivatives are taken with respect to holomorphic coordinates on
$\sigma$ and $\sigma'$, so this is well-defined. The quotient measures
the distortion that happens to a small quadrilateral under $\kappa$. It
equals 1 iff $\kappa$ is holomorphic in $z$. If $K$ is a nonnegative real
number, $\kappa$ is called {\em $K$-quasiconformal} if its maximal
dilation does not exceed $K$. It will be convenient to drop the requirement
of surjectivity of $\kappa$ and to apply the same terminology to
open embeddings.

$K(\kappa)$ is used in Teichm\"uller theory to define
the {\em Teichm\"uller distance} between
two hyperbolic Riemann surfaces as the infimum of $\log K(\kappa)$
over all orientation preserving diffeomorphism $\kappa$. The
Teichm\"uller distance is a complete metric on the
moduli space of nonsingular curves and it coincides with
the Kobayashi metric. Equivalently one could compare
via $\kappa$ the complex structures $j_\Sigma$ and $j_{\Sigma'}$ as
endomorphisms of the tangent bundle using the constant curvature
metrics. In the present context the use of maximal dilations seems
more geometric and natural to me.
\vspace{1ex}

$K$-quasiconformal maps enjoy the following compactness property.
Recall that a (possibly non-compact) Riemann surface is called
{\em hyperbolic} if its universal cover is the unit disk.
\begin{prop}\label{quasiconf_compact}
  Let $\kappa_\nu: \Sigma\rightarrow\Sigma'$ be a sequence of
  $K$-quasiconformal maps between (possibly non-compact) hyperbolic
  Riemann surfaces. Then there exists a subsequence converging locally
  uniformly to a map $\kappa: \Sigma \rightarrow \Sigma'$ that is either
  constant or $K$-quasiconformal.
\end{prop}
\pf
By lifting the $\kappa_\nu$ to the universal cover it suffices to prove the
lemma for maps between unit disks. This case is treated by Lemma~2.1
together with Theorem~2.2 in Chapter~1 of \cite{lehto}.
\qed
\vspace{1ex}

If $(C,{\bf x})$ is a stable marked curve the {\em hyperbolic metric} on
$C$ is defined as the Poincar\'e (constant curvature $-1$) metric on the
connected components of $C\setminus (C_\sing\cup |{\bf x}|)$. In
a holomorphic coordinate $z$ centered at a special point the hyperbolic
metric looks asymptotically ($|z|\rightarrow 0$) like the hyperbolic
metric $dz\wedge d\bar z/ (|z|^2\ln^2|z|)$ on the punctured disc. This shows
that with respect to the hyperbolic metric the length of the circle of radius
$|z|=\eps$ becomes arbitrarily small with $\eps$. We may thus use the
hyperbolic metric to define an intrinsic neighbourhood basis of the set of
special points: Let $V_\eps$ consist of the set of special points together
with all nonspecial points $P$ admitting a homotopically nontrivial
closed curve of length smaller than $\eps$ in $C\setminus
C_\sing\cup |{\bf x}|$, starting and ending at $P$. We call $V_\eps$ the
{\em $\eps$-thin part} of $(C,{\bf x})$. Since in any based homotopy class
of loops there is a unique geodesic loop one may replace the word ``loop''
in this definition by ``geodesic loop''.
\begin{lemma}
  $\eps$-thin parts form a neighbourhood basis
  of the set of special points.
\end{lemma}
\pf
Let $V$ be any neighbourhood of the set of special points.
Using compactness of $C\setminus V$ it is clear that for $\eps$
sufficiently small curves of length less than $\eps$ all enter $V$.
Then one can use the explicit form of the Poincar\'e metric on the
punctured disc to conclude that for even smaller $\eps$ the only
such curves lie entirely in $V$ and encircle a special point.
\qed
\vspace{1ex}

The topology of the moduli space of stable (marked) curves
can also be defined using dilations. Let $(C,{\bf x})$ be a
$k$-marked stable curve. For $\eps>0$
we consider the set $U_\eps(C,{\bf x})$
of stable curves $(\tilde C, \tilde{\bf x})$ admitting partial real
resolutions $\kappa: \tilde C\rightarrow C$ that are standard over
$V_\eps$ (i.e.\ of the form $Z_t\rightarrow Z_0$, cf.\ above) and with
\begin{itemize}
  \item
    $\log K(\kappa|_{\tilde C \setminus\kappa^{-1}(V_\eps)}) < \eps$
  \item
    $\kappa(\tilde x_\mu)\in V_\eps(x_\mu)$, the connected component
    of $V_\eps$ containing $x_\mu$.
\end{itemize}
To verify that the sets $U_{\eps,V}(C,{\bf x})$ form a neighbourhood basis
of $(C,{\bf x})$ in the topology defined by analytic deformation theory we
need the following statement.
\begin{lemma}\label{analytic-dilation_comp}
  Let $(C,{\bf x})$ be a $k$-marked stable curve and
  $(q:\calc\rightarrow S, \underline{\bf x})$ be a universal deformation.
  Then there exists $\eps>0$ such that for any $(\tilde C,\tilde {\bf x}) \in
  U_\eps(C,{\bf x})$ there exists an isomorphism $\Xi: (\tilde C,\tilde{\bf x})
  \simeq (C_s,{\bf x}_s)$ for some $s\in S$.
\end{lemma}
\pf
From the Banach-analytic constructions of universal families \cite{grauert},
\cite{douady}, \cite{palamodov} we know that $(\calc\rightarrow S,{\bf x})$ contains
all marked stable curves obtained from small deformations of the normalization
$(\hat C,\hat{\bf x}\vee {\bf y})$ by gluing with sufficiently small parameters
$\bf t$. To see this we remove $\kappa^{-1}(V_\eps)$ from $\tilde C$ and
fill in unit disks to arrive at a Riemann surface $({\hat C}',\hat{\bf x}'\vee{\bf y}')$,
where we place ${\bf y}'$ in the centers of the disks.
Since we assumed $\kappa$ to be standard over $V_\eps$ the restriction
of $\kappa$ to the complement of $\kappa^{-1}(V_\eps)$ can be extended to
a difffeomorphism $\hat\kappa:{\hat C}'\rightarrow \hat C$ with
dilation still smaller than $\eps$. Here $\hat C$ is the normalization of $C$.
Standard Teichm\"uller theory now shows that for $\eps$ tending to zero,
${\hat C}'$ converges to $\hat C$ in the analytic sense. The fact that the
size of the gluing regions (the preimage of $V$ in $\hat C$) decreases
with $V$ amounts to decreasing the gluing parameters $t$.
\qed
\vspace{1ex}

For an analytic or continuous family $(q:\calc\rightarrow S,
\underline{\bf x})$ of stable curves one expects that
the $\eps$-thin parts of the fibers form a neighbourhood basis
of the singular set of $q$ and the image of $\underline{\bf x}$.
To actually prove this one has to show continuity of the hyperbolic
metrics on the fibers (now viewed as metric on the relative tangent
bundle $T_{\calc/S}= \bigcup_s T_{C_s}$) away from the
special points of the fibers. This apparently classical question
concerning the compatibility of complex-analytic and hyperbolic points of
view has been settled only more recently by Wolpert using a gluing
construction for hyperbolic metrics. On the complex-analytic side he
uses in fact the same gluing as we did in our construction of the
universal family and the definition of continuous families of (pre-)stable
curves. Among more refined issues he showed
\begin{theorem}{\rm\cite{wolpert}}\ \label{conv_wolpert}
  Let $(\calc\rightarrow S, \underline{\bf x})$ be a continuous family of
  marked stable curves. Then the hyperbolic metrics on the fibers
  glue to a metric on $T_{\calc/S}$ that is continuous away from the
  special points.
\end{theorem}
Wolpert gives this result only for stable curves without marked points.
Marked points can however be easily included by going over from
$(C,{\bf x})$ to the singular stable curve obtained from two copies of $C$
by joining at the corresponding marked points.
\vspace{1ex}

Our main application of the notion of $\eps$-thin parts is the following
lemma that will be needed in the proof of the Hausdorff
property in Section~\ref{C(X)_topology}.
\begin{lemma}\label{eps-thins}
  Let $(C,{\bf x})$ and $(C',{\bf x}')$ be stable curves. Let $(C_\nu,{\bf
  x}_\nu)$ be a sequence of marked curves and $\kappa_\nu:
  C_\nu\rightarrow C$, $\kappa'_\nu: C_\nu\rightarrow C'$ be partial real
  resolutions exhibiting the convergence of this sequence to $(C,{\bf x})$
  and $(C',{\bf x}')$ respectively (as in the above definition of the
  topology on the moduli space using dilations).
  
  If $V_\nu$ is a neighbourhood basis of the set of special points
  on $C$ then $\kappa'_\nu\circ\kappa_\nu^{-1}(V_\nu)$ is a
  neighbourhood basis for the set of special points on $C'$.
\end{lemma}
\pf
We may assume $V_\nu$ is the $\eps_\nu$-thin part of $(C,{\bf x})$,
$\eps_\nu\rightarrow 0$. We have already seen that convergence using
dilations is compatible with the point of view of the gluing construction,
i.e.\ convergence in the analytic sense. Wolpert's result thus implies that
the hyperbolic metrics on $(C,{\bf x})$, $(C_\nu, {\bf x}_\nu)$ and $(C',{\bf
x}')$, compared via $\kappa_\nu$, $\kappa'_\nu$, converge to one
another away from the set of special points. This shows that there are
sequences $\delta_\nu, \delta'_\nu \rightarrow 0$ such that
$\kappa_\nu^{-1}$ of the $\eps_\nu$-thin part of $C$ contains the
$\delta_\nu$-thin part of $C_\nu$, which in turn contains
$(\kappa'_\nu)^{-1}$ of the $\delta'_\nu$-thin part of $C'$. Similarly for
the converse inclusion of sets.
\qed

The hypothesis of the lemma together with the Hausdorff property of
the moduli space of stable curves implies of course $(C,{\bf x}) \simeq
(C',{\bf x}')$. The point is that, by setting $X=\{{\rm pt}\}$, the
Hausdorff property of the moduli space will in fact be included in our
proof of the Hausdorff property of the space $\calc(X)$ to be introduced
in the next Chapter.


\section{Stable complex curves in Hausdorff spaces}
\subsection{Definition}
Given a complex (or tamed almost complex) manifold $M$, spaces of
holomorphic maps from
{\em stable} curves to $M$ (of fixed topological type) are in general not compact:
Additional rational components, so-called ``bubbles'', may appear in limits of
sequences of such maps. This is why we called irreducible components contracted
under $(C,{\bf x})\rightarrow(C,{\bf x})^\st$ bubble components. The fact that limit
maps are non-trivial on such components leads to the concept of bubble curve
(\cite{gromov},\cite{parkerwolf},\cite{ruan}). It is convenient to admit also
bubble components $D$ with $\ph|_D$ trivial, but only if $D$ is a
stable component (otherwise one could clearly contract $D$ without
leaving the category of maps with prestable domains). The role of stable bubble
components has been clearly recognized by Kontsevich in his notion of ``stable map''
\cite{kontmanin}. Stable maps are a natural way to treat maps from singular but stable
curves and from curves with bubble components on an equal footing.
To emphasize that we deal with maps from curves I propose to
adopt the term ``stable curve'' (and ``stable complex curve'' in a real setting) from
Deligne-Mumford theory. After all it is indeed an extension of Deligne-Mumford stability
for algebraic curves that goes well with the idea of a curve as a scheme of relative
dimension one over a base space. The only difference is that instead of
$\mbox{\rm Spec}\,k$ the base is now the variety the curve is mapped to.

Our basic notion is an extension of this notion to the topological setting where
the target space is a Hausdorff space $X$. Later on we will restrict ourselves
to manifolds and require maps of Sobolev class.
\begin{defi}\rm\label{cplx_curve}
A {\em (marked, parametrized) complex curve} in $X$ is a
triple $(C,{\bf x},\ph)$ with $(C,{\bf x})$ a marked prestable
curve and $\ph:C\rightarrow X$ a continuous map.

A morphism  $(C,{\bf x},\ph)\rightarrow (C',{\bf x}',\ph')$ is a holomorphic map
$\Psi:C\rightarrow C'$ with $\ph'\circ\Psi=\ph$, $\Psi({\bf x})={\bf x'}$ (as tuples).
This defines the sets $\Hom((C,{\bf x},\ph),$ $(C',{\bf x}',\ph'))$ and
$\Aut(C,{\bf x},\ph)\subset\Aut(C,{\bf x})$.

$(C,{\bf x},\ph)$ is called {\em stable} if the restriction of $\ph$ to
any unstable component $D\subset (C,{\bf x})$ is non-constant. We write
$\calc(X)$ for the space of isomorphism classes of complex curves in $X$ and
---by abuse of notation--- $(C,{\bf x},\ph)\in\calc(X)$.
\qed
\end{defi}
For later reference let us also introduce at this point stable complex curves of
Sobolev class and stable pseudo holomorphic curves. Recall that if $M$ is a
manifold and $C$ is a {\em fixed} prestable curve there are various natural metric
topologies on spaces of (continuous) maps from $C$ to $M$ by either using
local charts or an embedding $\iota: M\hookrightarrow\rz^N$,
together with Riemannian metrics on each irreducible component of $C$.
E.g.\ we obtain $C^0(C;M)$ and $L_1^p(C;M)$, $p>2$ (meaning that $\ph$ pulled
back to the normalization of $C$ has one weak derivative in $L^p$). We also fix a
Riemannian metric $\rho$ on $M$. Then any $\ph\in L_1^p(C;M)$ has
a well-defined area
\[
  {\rm area_\rho}(\ph)\ =\ \int_C |\ph^*\rho|\,d\mu\,.
\]
\begin{defi}\rm\label{sobolev_curves}
  Let $2<p<\infty$. The subset of $\calc(M)$ of stable {\em (marked, para\-met\-rized)
  complex curves of Sobolev class $L_1^p$} in $M$ is defined as set of
  (isomorphism classes of) curves $(C,{\bf x},\ph) \in\calc(X)$ with:
  \begin{enumerate}
  \item
    $\ph\in \check L_1^p(C;\ph^* T_M)$
  \item
    ${\rm area_\rho}(\ph|_D)>0$ for any unstable component $D\subset (C,{\bf x})$.
  \end{enumerate}
  \vspace{-4ex}
  
  \qed
\end{defi}
The condition of positive area is of course independent of the choice of $\rho$.
\begin{rem}\rm
Our stability condition for continuous curves is designed to make the
proof of the Hausdorff property work. It permits the case of bubble components
$D$ with $\ph|_D$ factorizing over the quotient by some compact subgroup
$D/S^1$. In this case $\ph|_D$ would be contractible. 

A slightly stronger, but possibly more natural stability condition, that is good
for both the continuous and the Sobolev world, would require the restriction
of the map to the bubble components to be {\em homotopically nontrivial}.
To relate this to our condition in Definition~\ref{sobolev_curves}
we observe that by \cite{sacks} the areas of maps
from the two-sphere to $M$ in a fixed homotopy
class are larger than the area of a certain harmonic map (that modulo
the action of the fundamental group lies in the specified homotopy class).
So any non-contractible $\ph\in L_1^p (S^2,M)$ indeed has positive area.
\qed
\end{rem}
\vspace{1ex}

Stable complex curves in $M$ of Sobolev class have isolated points in most fibers
of their unstable components.
\begin{lemma}\label{discrete_fibers}
  Let $\ph:S^2\rightarrow M$ be of Sobolev class $L_1^p$ and of positive area.
  Then there exists a set $A\subset \img\ph$ of positive 2-dimensional Hausdorff
  measure such that for any $P\in A$ the fiber $\ph^{-1}(P)$ is finite.
\end{lemma}
\pf
The proof relies on the somewhat deeper fact that any $L^p_1$-function coincides with
a continuously differentiable function away from ``bad'' sets $Z\subset S^2$ of
arbitrarily small measure \cite[Thm.3.11.6] {ziemer}. Let $\calh^i$ be
$i$-dimensional Hausdorff measure. On $U=S^2\setminus Z$ we may now
invoke the area formula for differentiable maps (cf.\ e.g.\ \cite[Thm.3.2.3] {federer}):
\[
  \int_U |D\ph|\,d\mu\ =\
  \int_M \calh^0((\ph|_U)^{-1}(P))\, d\calh^2(P)\,.
\]
So for $Z$ sufficiently small we conclude that the fiber of $\ph$ over $P$ is finite
for $P$ chosen from the set $A=\ph(S^2\setminus Z)$ of positive 2-dimensional
Hausdorff measure.
\qed
\begin{defi}\rm\label{def.ps.hol}
  Let $(M,J)$ be an almost complex manifold, $J$ the almost complex structure.
  $(C,{\bf x},\ph)\in\calc(M)$ is called {\em pseudo holomorphic (with respect to $J$)},
  or {\em $J$-holomorphic} if for any irreducible component $D$ of $C$, $\ph|_D:
  D\rightarrow M$ is a morphism of almost complex manifolds. The subset
  \[
    \calc^\hol(M,J)\ :=\ \{(C,{\bf x},\ph)\in \calc(M)\mbox{ $J$-holomorphic}\}
    /\mbox{isomorphism}
  \]
  of $\calc(M)$ is the space of (marked, parametrized) stable pseudo holomorphic
  curves on $M$ with respect to $J$.
  \qed
\end{defi}
Note that by elliptic regularity pseudo-holomorphic maps $\ph: C\rightarrow M$
are of Sobolev class $L_1^p$ for any $p$. And non-constant pseudo-holomorphic
maps have positive area because they are immersions away from finite sets of points.
This shows that $\calc^\hol(M,J)$ is in fact a subset of $\calc(M;p)$ for any $p$.


\subsection{The $C^0$-topology on $\calc(X)$}
\label{C(X)_topology}
We will now endow $\calc(X)$ with a topology that induces the compact-open
topology on the subsets $C^0(C,X)/\mbox{isomorphism}$ for any fixed prestable
curve $C$, and the Gromov topology on $\calc^\hol(M,J)$ in case
$X=M$ is a manifold with an almost complex structure $J$ (cf.\
Section~\ref{sect_compact}). The main technical result in this section is
a proof of the Hausdorff property for this topology.

The topology on $\calc(X)$ is defined in terms of dilations as follows.
\begin{defi}\label{C0top}
  \rm {\em The $C^0$-topology on $\calc(X)$}.\ \
  Let $(C,{\bf x},\ph)$ be a $k$-marked stable curve in $X$. Let $V$ be an open
  neighbourhood of $C_\sing\cup |{\bf x}|$ and let $N$ be a neighbourhood
  of the graph of $\ph$ in $C\times X$. $V_\mu$ denote the connected
  component of $V$ containing $x_\mu$. We define $U_{V,N}(C,{\bf x},\ph)$ as the
  set of stable curves $(\tilde C,\tilde{\bf x},\tilde\ph)$ admitting partial real resolutions
  $\kappa:\tilde C \rightarrow C$ that are standard over $V$ and with
  \begin{itemize}
  \item
    $\log K(\kappa|_{\tilde C \setminus\kappa^{-1}(V)}) < \eps$
  \item
    $\kappa(\tilde x_\mu)\in V_\mu$
  \item
    the closure of the graph of $\tilde\ph\circ\kappa^{-1}|_{C_{\rm \tiny reg}}$
    in $C\times X$ is contained in $N$.
  \end{itemize}
  The {\em $C^0$-topology on $\calc(X)$} is defined by taking the $U_{V,N} (C,{\bf
  x},\ph)$ as neighbourhood basis of $(C,{\bf x},\ph)$.
\qed
\end{defi}
In case $X$ is metric we can of course take for $N$ the $\eps$-neighbourhood of the
graph of $\ph$ with respect to the metric on $X$ and the chosen metric on $C$.
So in this case (or, more generally, if $X$ is first countable)
the $C^0$-topology on $\calc(X)$ is first countable.
\begin{rem}\rm
It would be slightly more economical for our purposes to define the topology by using
semiuniversal deformations of $(C,{\bf x})$ together with the retraction to the
central fiber $\kappa: \calc \rightarrow C$ as constructed in the last Chapter. However,
the given definition is somewhat more intrinsic and better to compare with
existing definitions for the Gromov topology on spaces of pseudo-holomorphic
curves, cf.\ below.
\qed
\end{rem}
\begin{prop}\label{C^0_hausdorff}
 The $C^0$-topology on $\calc(X)$ is Hausdorff.
\end{prop}
\pf
To avoid messing up the proof with ``nets'' we give the proof only for first
countable $X$. Then the topology on $\calc(X)$ has countable
neighbourhood bases too and it suffices to show uniqueness
of limits of converging sequences. So let $(C_\nu,{\bf x_\nu},
\ph_\nu)$ be a sequence of stable curves in $X$ converging to two stable curves
$(C,{\bf x},\ph)$ and $(C',{\bf x}',\ph')$. We have to show $(C,{\bf x},\ph) =(C',{\bf x'},
\ph')$. Let $\kappa_\nu: C_\nu\rightarrow C$ and $\kappa'_\nu: C_\nu\rightarrow C'$ be
partial real resolutions as in the definition of the topology on $\calc(X)$ exhibiting the
convergence of the sequence to $(C,{\bf x},\ph)$ and $(C',{\bf x'},\ph')$ respectively.

The first, most important step is a reduction to the case of stable domains
$(C,{\bf x})$, $(C',{\bf x}')$. This is where stability will be used in an essential way.
First we observe that by the definition of convergence $\img \ph=\img \ph'$. $(C,{\bf x})$
and $(C',{\bf x'})$ have only finitely many non-stable components and by stability each
of these components map nontrivially to $X$. We may therefore choose finitely many
points $y_j\in C$, $y'_{j'}\in C'$ ($j=1,\ldots,m$, $j'=1,\ldots, m'$)
stabilizing $(C,{\bf x})$ and $(C',{\bf x}')$ and such that
the $\ph(y_j)$, $\ph'(y'_{j'})$ are pairwise different.
We may also achieve the $\ph(y_j)$, $\ph'(y_{j'})$ to be away from
$\ph(C_\sing\cup |{\bf x}|)\cup \ph'(C'_\sing\cup |{\bf x}'|)$. We set
\[
  {\bf z}_\nu\ :=\ (\kappa_\nu^{-1}(y_1),\ldots,\kappa_\nu^{-1}(y_m),
                         {\kappa'_\nu}^{-1}(y'_1),\ldots,{\kappa'_\nu}^{-1}(y'_{m'})\,,
\]
and ${\bf z}^\nu$, ${{\bf z}'}^\nu$ for the image of this tuple under $\kappa_\nu$,
$\kappa'_\nu$ in $C$ and $C'$ respectively. By compactness of the curves
we may assume ${\bf z}^\nu$, ${{\bf z}'}^\nu$ converge to tuples ${\bf z}$ and ${\bf z}'$.
A priori these could have several coinciding entries or entries that are special points of
$(C,{\bf x})$ respectively $(C',{\bf x}')$. However, convergence of the graphs of the
maps (away from the singular set) shows that, for $\nu$ sufficiently large, the images of
the entries under $\ph$ respectively $\ph'$ are pairwise different and they are
also disjoint from the image of the set of special points. This shows that both
$(C,{\bf x}\vee{\bf z},\ph)$ and $(C',{\bf x}'\vee{\bf z}',\ph')$ are stable curves.
And, by construction,
\[
  (C_\nu,{\bf x}_\nu\vee{\bf z}_\nu,\ph)\ \stackrel{\nu\rightarrow\infty}{\longrightarrow}
  (C,{\bf x}\vee{\bf z},\ph)\,,\quad
  (C_\nu,{\bf x}_\nu\vee{\bf z}_\nu,\ph)\ \stackrel{\nu\rightarrow\infty}{\longrightarrow}
  (C',{\bf x}'\vee{\bf z}',\ph')\,.
\]
It thus suffices to treat the case of stable domains $(C,{\bf x})$, $(C',{\bf x}')$.
\vspace{1ex}

We now consider the sequence of maps
\[
  \Psi_\nu:= (\kappa_\nu)^{-1}\circ \kappa_\nu:\ C_\reg\longrightarrow C'\,.
\]
We only have to show that a subsequence converges to a biholomorphism $\Psi:
C_\reg\rightarrow C'_\reg$. In fact, the (uniquely existing) extension of $\Psi$ across the
singular set will intertwine $\ph$ and $\ph'$ by convergence of the graphs. So $(C,{\bf
x},\ph) = (C',{\bf x}',\ph')$ as stable curves in $X$.

Let $D$ be an irreducible component of $C$ and $V$ a
neighbourhood of the set of special points in $D$. By Lemma~\ref{eps-thins},
for $\nu$ sufficiently large, $\Psi_\nu$ induces a diffeomorphism of
$D\setminus V$ onto its image in $C'_\reg$ with dilation tending to zero in
the limit. By passing to a subsequence we may assume $\psi_\nu(D\setminus
C_\sing)$ to lie in a fixed irreducible component $D'$ of $C'$.
Now apply the compactness theorem for quasiconformal maps to a hyperbolic
Riemann surface (Lemma~\ref {quasiconf_compact}) to find a subsequence
converging to a holomorphic map $D\setminus V \rightarrow D'$.
A diagonal argument allows to do this simultaneously for $V$ running through a
neighbourhhod basis of the set of special points and for all irreducible components
of $C$. This defines $\Psi$ away from the set of special points. 

To see that $\Psi$ is indeed a biholomorphism from $C\setminus C_\sing\cup|{\bf x}|$
to $C'\setminus C'_\sing\cup|{\bf x}'|$ (and hence extends to a biholomorphism
$C\rightarrow C'$), we repeat the argument for
$\Psi'_\nu:= \kappa_\nu \circ (\kappa'_\nu)^{-1}$, with $\nu$ running over the
subsequence we have just constructed. Possibly after reducing the set of indices
once more we obtain a holomorphic map $\Psi'$ from $C'$ to $C$ away from the
special points of $C'$. $\Psi'$ is inverse to $\Psi$ because $\Psi'_\nu
\circ \Psi_\nu=\id$ for $\nu$ sufficiently large locally away from the singular set.
For this last step Lemma~\ref{eps-thins} is crucial; it says that large sets in
$C\setminus C_\sing\cup|{\bf x}|$ correspond to large sets in $C'\setminus
C'_\sing\cup|{\bf x}'|$ under $\Psi_\nu$ for $\nu$ sufficiently large. This concludes
the proof of the proposition.
\qed


\subsection{Intrinsic neighbourhood bases for the singular set}
\label{section_intrinsic}
One of the benefits of stability is that it allows approximate measurements
of the size of sets on unstable components of the domain.
This can not be done by just looking at the prestable curve due to
the non-compact nature of the automorphism group on such components.
Implicitely this fact has been used in the proof of the Hausdorff property
of the $C^0$-topology (Proposition~\ref {C^0_hausdorff}).
The following version for families will be convenient for our construction of Kuranishi
structures in Section~\ref{section_kuranishi}. We now assume $X=M$ to be
a Riemannian manifold.

We consider a continuous family of marked prestable curves $(\Gamma\rightarrow
T,\underline{\bf x})$. Let $\Phi:\Gamma\rightarrow M$ be a continuous map
in such a way that $(C_t,{\bf x}_t,\ph_t):= (q^{-1}(t), \underline{\bf x}(t), \Phi|_
{q{-1}(t)})$ is a stable curve in $M$ for any $t\in T$. Assume that $\ph_t$ is
not contractible for all $t\in T$. We also assume that the
number of special points and the number of unstable components of
$(C_t,{\bf x}_t)$ are uniformly bounded by some numbers $c$ and $b$ respectively.
Observe that since $\ph_t$ is non-trivial on any unstable component
the diameter of $\ph_t(D)$ is larger than the infimum ${\rm inj}_M(\ph)$ over
the injectivity radii of $M$ at any $P\in \ph(D)$. If there are no unstable components
the same holds for some stable irreducible component by the assumption of topological
non-triviality. There are at most $2b$ points needed to rigidify
$(C_t,{\bf x}_t)$. So for any $t$ there exist tuples ${\bf y}$ of $2b$
points of $(C_t,{\bf x}_t)$ with
\begin{itemize}
\item
  $(C_t,{\bf x}_t\vee {\bf y})$ is a {\em stable} curve (i.e.\ $\bf y$
  is a rigidifying tuple of points for $(C_t,{\bf x}_t, \ph_t)$
\item
  for any $j$ and any special point $P\neq y_j$ of $(C_t,{\bf x}_t\vee {\bf y}_t)$
  it holds
  \[
    d(y_j, P)\ge \frac{{\rm inj}_M (\ph_t)}{2(c+2b)}\,.
  \]
\end{itemize}

The set of such curves can be identified with a subset $\tilde T$ of the $2b$-fold
fibered product $\Gamma\times _T\times \ldots \times_T\Gamma$.
The pull-back of $\Gamma$, $\underline{\bf x}$ and $\Phi$ under the
projection $\tilde T\rightarrow T$ is a family of stable curves
$(\tilde\Gamma \rightarrow \tilde T, \tilde{\underline{\bf x}}, \tilde\Phi: \tilde\Gamma
\rightarrow M)$ with stable domains. It follows from Wolpert's result
Theorem~\ref{conv_wolpert} that for any $\eps>0$ the union of the $\eps$-thin
parts of the fibers of $(\tilde\Gamma \rightarrow \tilde T, \tilde{\underline{\bf x}})$
form a neighbourhood basis in $\tilde\Gamma$ of the union of the special points.
Now by choice of $\bf y$ the projection $\tilde \Gamma\rightarrow\Gamma$ is
proper. Hence the image of this neighbourhood basis forms a neighbourhood basis
$V_\eps$ of the set of special points of $(\Gamma\rightarrow T, \underline{\bf x})$.
Note that $V_\eps$ only depends on the family and the Riemannian metric on $M$.
We call $V_\eps$ the {\rm $\eps$-thin part of $(\Gamma\rightarrow T, \underline{\bf x},
\Phi)$}.


\subsection{Pseudo-holomorphic curves and compactness}
\label{sect_compact}
The $C^0$-topology on $\calc(X)$ induces Hausdorff topologies on any of its subsets,
hence on the set of pseudo-holomorphic curves $\calc^\hol(X)$. On this subset
Gromov has already defined a topology in a slightly different way as follows
(we follow the definition given in \cite{pansu},\cite{parkerwolf},\cite{ye}): Let $(M,J)$ be
an almost complex manifold and $(C,{\bf x},\ph)$ be a marked prestable curve. We also
fix a Riemannian metric on $M$. Then given $\eps>0$ and a neighbourhood $V$ of
the set of special points of $(C,{\bf x})$ the set $U'_{\eps,V}$ is defined as the set of
curves $(\tilde C,\tilde{\bf x},\tilde\ph)$ allowing partial real resolutions $\kappa:\tilde C
\rightarrow C$ that are standard over $V$ and with the following properties:
\begin{itemize}
\item
    $\log K(\kappa|_{\tilde C \setminus\kappa^{-1}(V)}) < \eps$
  \item
    $\kappa(\tilde x_\mu)\in V_\mu$, the connected component of
    $V$ containing $x_\mu$
  \item
    $d_{C^0(C\setminus V)}(\tilde\ph\circ\kappa^{-1},\ph) <\eps$
  \item
    $\big|{\rm area}(\ph|_V)-{\rm area}(\tilde\ph|_{\tilde\kappa^{-1}(V)})\big|\,<\,\eps$.
\end{itemize}
The {\rm Gromov topology} is defined as the topology generated
by sets of the form $U'_{\eps,V}$. So the requirement of uniform convergence of the
graphs near the singular points is replaced by convergence of areas. Instead
one requires the area of nearby maps to be small in $V$. This prohibits
bubbling off at the singular points.

We have slightly extended the definition found in the literature by
inclusion of marked points. This is not an issue for one can always
replace marked points by nodes, cf.\ the discussion following
Theorem~\ref{conv_wolpert}.

The Gromov topology is the topology used in compactness theorems for
pseudo-holomorphic curves. Recall that an almost complex
manifold $(M,J)$ has {\em bounded geometry} if it admits a Hermitian metric $\rho$
with nonzero injectivity radius and with sectional curvature and $|DJ|$
uniformly bounded.
\begin{theorem}\label{cptness_thm}
  (Gromov compactness theorem.)\ \ 
  {\rm\cite{gromov}, \cite{pansu}, \cite{parkerwolf}, \cite{ruantian1}, \cite{ye}}\ 
  Let $(M,J)$ be an almost complex manifold with bounded geometry. For $C>0$,
  $g$, $k\in\nz$ put
  \[
    \calc^{\mbox{\scriptsize Vol}_\rho<C}_{g,k}(M;p):= \{(C,{\bf x},\ph)\in\calc(M;p)
    \mid{\rm Vol}_\rho(\ph)\le C, g(C)=g, \sharp{\bf x}=k\}\, .
  \]
  Then $\calc^\hol(M,J)\cap\calc^{\mbox{\scriptsize Vol}_\rho<C}_{g,k}(M;p)$
  (this space is independent of $p$) is compact in the Gromov topology. In particular,
  if $(M,\omega)$ is a closed symplectic manifold, $J$ is tame with respect to $\omega$
  and $R\in H_2(M;\gz)$, let
  \[
     \calc_{R,g,k}(M;p):=\{(C,{\bf x},\ph)\in\calc(M;p)\mid\ph_*[C]=R, g(C)=g,
     \sharp{\bf x}=k\}\,,
  \]
  then $\calc^\hol_{R,g,k}(M,J):=\calc^\hol(M,J)\cap\calc_{R,g,k}(M;p)$ is compact.
\qed
\end{theorem}
To use this result in our context we have to show that the Gromov topology and
the $C^0$-topology on $\calc^\hol(M,J)$ coincide. In fact, one implication suffices.
The proof is a simple application of Gromov's monotonicity lemma for
pseudo-holomorphic curves, cf.\ e.g.\ \cite[Prop.4.3.1]{sikorav}. Note that the
assumption of tameness of $J$ that is usually made is superfluous for this
local form of monotonicity. 
\begin{lemma}
  Let $(M,J)$ be an almost complex manifold and $\ph:\Sigma \rightarrow M$ a
  pseudo-holomorphic curve, $\Sigma$ a Riemann surface without boundary.
  For $Q\in \img\ph$ there exists a constant $C$ such that
  \[
    {\rm area}(\ph|_{\ph^{-1}(B_\eps(Q))})\ \ge C\cdot\eps^2\,,
  \]
  provided $\eps$ is sufficiently small, the ball taken with respect to some fixed
  Riemannian metric on $M$. $C$ may be chosen uniformly for $Q$
  varying in a small set in $M$.
\qed
\end{lemma}
\begin{prop}
  The Gromov topology on $\calc^\hol(M,J)$ is finer than the $C^0$-topology.
\end{prop}
\pf
We have to show that convergence of a sequence $(C_\nu,{\bf x}_\nu,\ph_\nu)$
to $(C,{\bf x},\ph)$ in the Gromov topology implies convergence in the
$C^0$-topology. Let $\kappa_\nu: C_\nu \rightarrow C$
be a sequence of partial real resolutions exhibiting the
convergence of $(C_\nu,{\bf x}_\nu,\ph_\nu)$ to $(C,{\bf x},\ph)$ in the
Gromov topology.

If $\ph_\nu\circ\kappa_\nu^{-1}$ does not converge to $\ph$ uniformly
on $C_\reg\setminus|{\bf x}|$ (i.e.\ in an $L^\infty$-sense on all of $C$), then
there exists a neighbourhood basis $W_\nu$ of a special point $P\in C$ and
points $P_\nu \in \kappa_\nu^{-1} (W_\nu)$ with $\ph_\nu(P_\nu)$
not converging to $\ph(P)$. By continuity of $\ph$
and local $C^0$-convergence of $\ph_\nu\circ \kappa_\nu^{-1}$
towards $\ph$ away from the special points, we may assume
$\ph_\nu(P_\nu)$ to stay in a compact subset of $M$. After going over to a
subsequence, $Q_\nu:= \ph_\nu(P_\nu)$ converges to some point $Q\neq
\ph(P)\in M$. The monotonicity lemma applied to $\ph_\nu$ in $Q_\nu$ gives a
strictly positive lower bound for the area of $\ph_\nu$ away from a neighbourhood
of $\img\ph$. But by area convergence this area becomes arbitrarily small
with increasing $\nu$, contradiction!
\qed

\begin{rem}\rm
  We will see later that the converse implication is also true: Let $(C_\nu,{\bf x}_\nu,
  \ph_\nu)$ be a sequence of pseudo-holomorphic curves converging
  in the $C^0$-topology to $(C,{\bf x},\ph)$. Then Proposition~\ref{covtop} shows that
  the area of $\ph_\nu$ over a neighbourhood basis $V_\nu$ of the set of special
  points of $(C,{\bf x})$ tends to zero with $\nu$. So the Gromov topology and the
  $C^0$-topology on the set of pseudo-holomorphic curves coincide.
\qed
\end{rem}

\begin{corollary}
  Let $(M,\omega,J)$ be a closed symplectic manifold with tame almost complex
  structure and $R\in H_2(M;\gz)$, $g,k\ge0$. Then $\calc^\hol_{R,g,k}(M,J)$ is compact
  in the $C^0$-topology.
\qed
\end{corollary}


\section{The \boldmath $\dbar$-operator in families of prestable curves}
\label{chapter_dbar}
We study the behaviour of the $\dbar$ equation on a holomorphic vector
bundle over a continuous family of prestable curves
(Definition~\ref{cont_curves}). Locally, $L^p$ spaces are naturally
identified in such a family. We show that with this trivialization the
graphs $\Gamma_\dbar\subset L^p\times L^p$ are naturally trivialized
too, i.e.\ they fit into a Banach bundle over the parameter space
of the family. A guide provided the
case $p=2$, that has been dealt with by Seeley and Singer using Hilbert
space methods \cite{sesinger}.

\subsection{The Banach bundle $q^p_* E$}
To begin with we want to extend sheaves of germs of Sobolev sections,
a well-defined notion for any differentiable vector bundle over
any differentiable manifold, over the nodes of a prestable curve. The
extension should be locally trivializable in continuous families. We first
treat our standard family $\{Z_t\}$.

Recall the retractions $\lambda_\tau: A_\tau \rightarrow \Delta$ defined
in the text before Definition~\ref{partres}. $\lambda_\tau$ has the
feature of preserving the cylindrical measure (of finite length) $dr\,d\ph$
on $\Delta^*$, up to multiplication by $1/(1-\tau)$. Letting ${\check
L}^p(\Delta):= L^p(\Delta,dr\,d\ph)$ and $\check L^p(Z_0) ={\check
L}^p(\Delta) \amalg \check L^p (\Delta)$ for $1\le p\le \infty$, we thus get
isomorphisms
\begin{eqnarray*}
  \lambda_\tau^*:\  \check L^p(\Delta) &\longrightarrow&
  L^p(A_\tau)\\
  \kappa_t^*:\  \check L^p(Z_0) &\longrightarrow&
  L^p(Z_t)=\check L^p(Z_t)\,,
\end{eqnarray*}
and similarly for the corresponding sheaves of sections ${\check\call}^p
(\Delta)$, ${\check\call}^p (Z_t)$. For 1-forms we use the frame
$dz\amalg dw$ on $Z_0$ to reduce to the case of functions, i.e.\ on
$\Delta^+$, say, we set $\lambda_\tau^\star(f\cdot dz+g\cdot d\bar
z):= (\lambda_\tau^*f) \cdot dz+ (\lambda_\tau^*g)\cdot d\bar z$. By
abuse of notation we write $\check L^p(Z_t;\Omega_{Z_t})$, ${\check
L}^p(Z_t; \bar\Omega_{Z_t})$ for the corresponding function spaces of
$(1,0)$ respectively $(0,1)$ forms, and ${\check\call}^p (Z_t;
\Omega_{Z_t})$, ${\check\call}^p (Z_t;\bar \Omega_{Z_t})$ for the
sheaves.
\begin{rem}\label{retractions}\rm
1)\ \ There are many choices possible for the retraction $\lambda_\tau$.
Ours is distinguished by the property that for $2\le p<\infty$ functions of
class $\check L_1^p$ on $Z_0$ (i.e.\ $L_1^p$ with weights as before
on each irreducible component and continuous) pull back to $L_1^p$
functions on $Z_t$. It is also the
retraction used by Seeley and Singer in their treatment of the case $p=2$.

Another obvious choice would be the one induced by the gradient flow of
$|z|\cdot|w|$ on $\Delta\times\Delta$ composed with a normalizing
stretch ($|z|^2-|w|^2$ is invariant under this flow). This leads to
\[
  re^{i\ph}\longmapsto\sqrt{\frac{1-(\tau/r)^4}{1-\tau^4}}\, re^{i\ph}\,,
\]
This retraction identifies $L^p(\Delta\setminus\Delta_\tau)$ with
$L^p(\Delta)$ but is not Lipschitz. Everything we say in this chapter
works in this setting too, with the usual Calderon-Zygmund inequality.

\noindent
2)\ \ We emphasize with a different star that $\lambda_\tau^\star$ is 
{\em not} the usual pull back $\lambda_\tau^*$, which would introduce
poles along $S_t$ due to the singular nature of $\lambda_\tau$.
In contrast to $\lambda_\tau^*$ our choice of pull-back also preserves
the decomposition of forms into type.

$\lambda_\tau^\star$ gives rise to an essentially non-holomorphic
``trivialization'' of $\coprod_t\Omega_{Z_t}$ (which taken literally is only
locally free away from the double point). In fact, while
$\lambda_\tau^\star(dz\amalg dw)$ {\em is} holomorphic outside the
diagonal $\delta=\bigcup_t S_t$ and bounded, it is not continuous
along $\delta$. Our trivialization should also be contrasted to the natural
holomorphic one given by $dz/z=-dw/w$, which again leads to poles on
$Z_0$ (the corresponding space would be denoted $\check L^p
(\omega_C)$). Alternatively, we could work with the bundle $\omega_C$
and introduce the $p$-dependent measure $r^{-p}\,dr\, d\ph$.
\qed
\end{rem}

The treatment of the global situation is now immediate: Given a
continuous family of holomorphic vector bundles over prestable curves
$E\downarrow q:\calc\rightarrow S$ let
\[
  \tilde\kappa:\ E\ \longrightarrow\ E_0:= E_{s_0}
\]
be an admissable retraction (for definition cf.\ the end of
Section~\ref{sect_cont_fam}). We thus get a family of identifications of
sheaves
\[
  {\tilde\kappa}_s^*:{\check\call}^p(E_0)\longrightarrow
  {\check\call}^p(E_s)\, .
\]
For 1-forms one has to specify a partition of unity $\{\rho_i\}$
subordinate to $\calu=\{U_i\}$ and a continuous choice of holomorphic
1-form $dz$ on $U_0=(\Sigma,j_s)$ (for such a global holomorphic 1-form
to exist one might have to enlarge the set of $U_i$, $i>0$) and set
\[
  {\tilde\kappa}_s^\star:{\check\call}^p(E_0\otimes\Omega_{C_0})
  \longrightarrow {\check\call}^p (E_s\otimes\Omega_{C_s})\,,\quad
  \alpha\longmapsto \sum{\tilde\kappa}_s^\star(\rho_i\cdot\alpha)
\]
with ${\tilde\kappa}_s^\star$ on ${\check\call}^p (E_0\otimes\Omega)
|_{U_i}$ defined using $dz$ ($i=0$) respectively $dz\amalg dw$. Similarly
for $\bar\Omega$.
\begin{prop}\label{Lptriv}
  Let $E\downarrow q:\calc\rightarrow S$ be a continuous family
  of holomorphic vector bundles over prestable curves. Then for
  $1\le p\le\infty$
  \[
    q_*{\check\call}^p(E)\ :=\ \coprod_s\check L^p (C_s;E_s)
  \]
  has naturally a structure of Banach bundle with
  local trivializations ${\tilde\kappa}^*_s$ for any admissable retraction
  $\tilde\kappa$. Similarly for $q_*{\check\call}^p (E\otimes\Omega)$ and
  $q_*{\check\call}^p(E\otimes\bar \Omega)$ (with ${\tilde\kappa}
  ^\star_s$).
\end{prop}
\pf
We have to show that if $\tilde\kappa$, ${\tilde\kappa}'$ are two
admissable retractions of $E$ defined over $W$, $W' \subset S$ with
possibly different centers $0,0'\in S$, then
\[
  (\tilde\kappa^*)^{-1}\circ \tilde\kappa'^*:
  (W\cap W')\times \check L^p (C_{0'};E_{0'})\ 
  \longrightarrow\ \check L^p (C_0;E_0)
\]
is continuous with respect to the Banach space topologies on the ${\check
L}^p$-spaces. So let $(s,v), (s_0,v_0) \in (W\cap W')\times \check L^p
(C_{0'};E_{0'})$. Notice that while the map $\Psi:= \tilde\kappa' \circ
\tilde\kappa^{-1}$ is only well-defined away from a family of contracted
circles (which do not affect arguments concerning $L^p$-spaces), it is
uniformly Lipschitz relative $S$. Let $v_{00}$ be a step function
approximating $v_0$. We have the estimate
\[
  |\!|\Psi_s^*v - \Psi_{s_0}^* v_0|\!|_{\check L^p}
  \le |\!|\Psi_s^*v - \Psi_s^* v_{00}|\!|_{\check L^p} +
  |\!|\Psi_s^*v_{00} - \Psi_{s_0}^* v_{00}|\!|_{\check L^p}+
  |\!|\Psi_{s_0}^*v_{00} - \Psi_{s_0}^* v_0|\!|_{\check L^p}\,.
\]
The first and the last term depend only on the Lipschitz constant of
$\Psi$ relative $S$ and $|\!| v-v_0 |\!|_{\check L^p}$, $|\!| v_{00}-v_0
|\!|_{\check L^p}$ respectively. So these terms can be made arbitrarily
small by choosing $v$ close to $v_0$. Once $v_{00}$ has been chosen, the
middle term can be made arbitrarily small by choosing $s$ close
to $s_0$. In fact, this last step just needs convergence of $\Psi_s
= \Psi| _{C_s}$ towards $\Psi_{s_0}$ in a measure theoretic sense.

The same arguments apply for $1$-form valued sections too, regardless
of the discontinuity along the $S_t$.
\qed


\subsection{Banach bundles of holomorphic \v Cech cochains}
\label{cech}
Trivializations for higher Sobolev spaces are less obvious for there is a
change of topology of the curves when a double point resolves. 
We will solve this problem in the next section for $L_1^p$-spaces by
reduction to $q_*{\check\call}^p(E)$ and (``banachized'') holomorphic \v
Cech cochains via the $\dbar$-operator. The present section deals with \v
Cech cochains.

For a holomorphic vector bundle $E$ over a relatively compact, open set
$U$ in a complex space $X$ we write
\[
  \calo_{(c)}(\overline{U}; E)\ :=\ \calo(U;E) \cap C^0(\overline{U}; E)\,,
\]
for the space of holomorphic sections extending continuously to the boundary.
Choosing a hermitian metric on $E$, the sup-norm makes $C^0_{(c)}(U;E)$
into a Banach space. Now given a covering $\calu:=\{ U_i\}$ of a
prestable curve $C$ and a holomorphic vector bundle $E$ over $C$
we define spaces of (alternating) \v Cech cochains
\[
  \check C^j_{(c)} (\calu;E)\ := \prod_{I\mid \sharp I=j} 
  C^0_{(c)} (U_I,E|_{U_I})\,,
\]
where $U_I:= U_{i_1}\cap\ldots\cap U_{i_j}$ and $I$ runs over all
multi-indices $(i_1,\ldots, i_j)$ with $i_1 <i_2 \ldots <i_j$. And for a
family $E\downarrow q:\calc\rightarrow S$ of holomorphic vector bundles
over prestable curves with covering $\calu$ as in Definition~\ref
{cont_curves} we put
\[
  q_*^{(j)}(\calu,E)\ :=\ \coprod_{s\in S} \check C^j_{(c)} (\calu_s; E_s)\,,
\]
with index $s$ meaning restriction to $C_s$. Since in our covering triple
intersections are empty, $q_*^{(j)}(\calu,E) =0$ unless $j=0,1$. So there is only
one nontrivial  \v Cech differential
\[
  \check d:\ (h_i)_i\ \longmapsto\ (h_j-h_i)_{i<j}\,.
\]
For the following central result of this section we now assume that $U_0$
is a family of bounded {\em circular regions}, with boundaries varying
continuously with the parameter (cf.\ the discussion before Lemma~\ref{retraction}).
\begin{prop}\label{Cech-zero-cycles}
  For $j=0,1$, $q^{(j)}_* (\calu,E)$ have natural structures of Banach bundles
  over $S$.
\end{prop}
\pf
It is well-known that the restriction map
\[
  \calo_{(c)}(\overline{\Delta})\ \longrightarrow\  
  C^0(S^1=\partial\Delta)\,,\quad
  f\ \longmapsto\ f|_{\partial\Delta}
\]
identifies $\calo_{(c)}(\overline{\Delta})$ with the Banach subalgebra
$C^0_{\ge0}(S^1)\subset C^0(S^1)$ characterized by the vanishing of
negative Fourier coefficients (cf.\ e.g.\ \cite[\S\,20.4]{conway}).
Similarly, if $G\subset\cz$ is a domain with outer boundary $N_1=\partial
B_{r_1}(c_1)$ (i.e.\ $G\subset B_{r_1}(c_1)$) and inner boundaries
$N_2=\partial B_{r_2}(c_2),\ldots,N_k=\partial B_{r_k}(c_k)$ we get an
isomorphism
\[
  \calo_{(c)}({\bar\Delta})\ \longrightarrow\ 
  C^0_{\ge0}(N_1)\oplus\bigoplus_{j=2}^k C^0_{<0}(N_i)\,,
  \quad h\ \longmapsto\ (h_1|_{N_1},\ldots,h_k|_{N_k})
\]
with
\begin{eqnarray*}
  h_1(z)&=&\frac{1}{2\pi i}\int_{N_1}\frac{h(w)}{w-z}\, dw\,,
  \quad z\in B_{r_1}(c_1),\\
  h_i(z)&=&\frac{1}{2\pi i}\int_{-N_i'}\frac{h((w-c_i)^{-1})}
  {w-(z-c_i)^{-1}}\, dw
  -h_1(c_1),\ z\in\cz\setminus\overline{B_{r_i}(c_i)},\ i=2,\ldots,l\,,
\end{eqnarray*}
(these extend continuously to $\overline{G}$) where $N'_i= \{z\in\cz
\mid(z-c_i)^{-1} \in N_i\}$ and $C^0_{<0}(S^1)$ is the space of continuous
functions with vanishing non-negative Fourier coefficients. Note that
$h_i$ can be recovered by the above Cauchy integrals from the limit on
$N_i$. The inverse to the above isomorphism is just given by summing
over these $h_i$.

In our family over $S$ we simply choose isomorphisms of $\partial
U_0(s)$ for different $s$ as $U(1)$-homogeneous spaces, varying
continuously with $s$. This identifies $C^0(\partial U_0(s))$ in a way that
respects the subalgebras defined by properties of Fourier series, thus
trivializing $\coprod_s\calo_{(c)}( \overline{U_0(s)}; E_s)$ by application
of a vector valued version of the above to each connected component of
$U_0$. For the possibly singular parts $U_i(s)$, $i>0$, the identification of
the boundary components is already given, and the above trivialization
extends over singular curves by setting
\[
  \calo_{(c)} (\overline{Z_0})\ \longrightarrow\
  \calo_{(c)} (\overline{Z_t})\,,
  \quad h^+\amalg h^-\longmapsto\prj_1^*h^+ +\prj_2^*h^- -h^-(0)\, .
\]
Similarly for $U_{ij}$ (which can be exhibited as circular regions by
projection $U_{ij}\subset Z_t\rightarrow \Delta$ on one of the coordinate
axes) leading to local trivializations for 1-cochains.

Different choices of isomorphisms with circular domains lead to changes
of trivialization that are even norm-continuous families of isomorphisms.
\qed


\subsection{A quasi-isomorphism and ${q_1^p}_*E$}
\label{Sobolev_sections}
For a family $E\downarrow q:\calc\rightarrow S$ of holomorphic vector
bundles over prestable curves we now want to use the $\dbar$-operator to
reduce the case of Sobolev sections $L_1^p$ to $\check L^p$. Viewing
$\dbar$ as closed unbounded operator from $\check L^p (\Delta)$ to $\check
L^p (\Delta,\bar\Omega)$ it is natural to define a (weighted) Sobolev
space ${\check L}^p_1(\Delta) \subset L_1^p(\Delta)$ by the graph norm. We
will see that this is equivalent to take the closure of $C^\infty (\bar\Delta)$
with respect to the norm
\begin{eqnarray*}
  |\!|f|\!|_{\check L_1^p(\Delta)}&:=&|\!|f|\!|_{\check L^p(\Delta)}
  +|\!|\nabla f|\!|_{\check L^p(\Delta)}\\
  &=&\int_\Delta\Big(|f|^p+|\partial_r f|^p+|r^{-1}\partial_\ph f|^p\Big)
  \,dr\,d\ph\, .
\end{eqnarray*}

For $2< p<\infty$ we define $\check L^p_1(Z_0):=({\check
L}_1^p(\Delta)\amalg \check L_1^p(\Delta))\cap C^0(Z_0)$, i.e.\ we
require continuity at the node (this makes sense in view of ${\check
L}^p_1(\Delta)\subset L^p_1(\Delta) \subset C^{1-(2/p)} (\overline{\Delta})$). The
corresponding sheaf of germs extends $\call_1^p (E|_{C\setminus
C_{\mbox{\tiny sing}}})$ to a sheaf denoted ${\check \call}_1^p(E)$, $E$
any holomorphic vector bundle over $C$.
\vspace{1ex}

We now turn to the $\dbar$-operator, as map from $\check L_1^p$ to
$\check L^p$. Recall the right inverse
to $\dbar$ on the disk (cf.\ \cite[Ch.1]{vekua} for an extensive discussion)
\[
  T(g d\bar z)(z):=\frac{1}{2\pi i}\int_\Delta\frac{g(w)}{w-z}\,
  dw\wedge d\bar w\, .
\]
This is a weakly singular integral operator and as such defines a compact
mapping from $\check L^p(\Delta;{\bar\Omega}_\Delta)$ to
$C^0(\bar\Delta)$. In fact, as we will see instantly,
$T$ is a bounded map to $L_1^p(\Delta)\subset
C^{1-(2/p)}(\overline{\Delta})$. It fulfills $\dbar\circ T=\id$ on $L_1(\Delta)$.
To study the regularity properties of $T$ one introduces
\[
  S(g d\bar z)(z):=\frac{1}{2\pi i}\Bigg(\lim_{\eps\rightarrow 0}
  \int_{\Delta\setminus B_\eps(z)} \frac{g(w)}{(w-z)^2}\,
  dw\wedge d\bar w\Bigg)\, dz\,.
\]
$S$ is a singular integral operator and hence extends to a continuous
map from $L^p(\Delta;\bar\Omega_\Delta)$ to $L^p(\Delta;
\Omega_\Delta)$ (Calderon-Zygmund) and satisfies
$\partial\circ T=S$ weakly. The existence of $T$ and $S$ imply
that $\dbar: L_1^p(\Delta)\rightarrow L^p(\Delta)$ is
surjective.

To conclude the analogous statement for our weighted spaces
we have to investigate $S$ on $\check L^p$:
The weight functions $w\in L^1_{\rm\scriptsize loc}$ such that $S$ is a continuous
endomorphism of $L^p$-spaces with weights $w$ are exactly those
fulfilling the so-called Muckenhoupt $A_p$-condition \cite{coifman}. The condition
requires the existence of a constant $C$ such that for any quader $Q$:
\[
  \Big(\int_{Q} w\,d\mu\Big)\cdot
  \Big(\int_{Q} w^{-\frac{1}{p-1}}d\mu\Big)^{p-1}
  \ \le\ C\cdot(\Vol Q)^p\,.
\]
The Muckenhoupt $A_p$-condition is fulfilled for our weights $w=1/r$.
We conclude
\begin{prop}\label{CalderonZygmund}
  $S:\ \check L^p(\Delta;\bar\Omega_\Delta) \rightarrow
  \check L^p(\Delta;\bar\Omega_\Delta)$ is continuous.
\qed
\end{prop}
This proves that $T$ indeed maps $\check L^p(\Delta;{\bar\Omega})
\subset L^p(\Delta; {\bar\Omega})$ to $\check L_1^p(\Delta)\subset L_1^p(\Delta)$.
Alternatively, this follows from elliptic theory on manifolds with conical ends
\cite{lockmac}, \cite{mazya}. Together with interior elliptic estimates and the
imposed continuity at the nodes this shows:
\begin{prop}
  Let $E$ be a holomorphic vector bundle over a prestable curve $C$.
  Then for $2< p<\infty$
  \[
    0\longrightarrow\calo(E)\longrightarrow{\check\call}^p_1(E)
    \stackrel{\dbar}{\longrightarrow}{\check\call}^p
    (E\otimes\bar\Omega_C) \longrightarrow 0
  \]
  is exact.
\qed
\end{prop}
\begin{rem}\rm
Note that while elliptic theory of the $\dbar$ operator holds true for
higher Sobolev spaces ($k>1$) as well, only a subspace of ${\check
L}^p_{k+1}$ is mapped by $\dbar$ to $\check L^p_k$ due to the necessary
continuity at the node. This has the effect of higher tangency conditions,
which are unwanted in our application to singular pseudo holomorphic
curves.
\end{rem}
\begin{corollary}\label{ker_coker}
  There is an exact sequence
  \[
    0\longrightarrow\Gamma(C;\calo(E))\longrightarrow \check L_1^p(C;E)
    \stackrel{\dbar}{\longrightarrow} \check L^p(C;E\otimes\bar\Omega)
    \longrightarrow H^1(C;\calo(E))\longrightarrow 0\,,
  \]
  where $\check L_1^p(C;E):=\Gamma{\check\call}_1^p(C;E)$ and
  $\check L^p(C;E\otimes\bar\Omega):=\Gamma{\check\call}^p
  (C;E\otimes\bar\Omega)$. In particular, the $\dbar$-operator is Fredholm
  as map from $\check L_1^p (C;E)$ to $\check L^p(C;E\otimes\bar
  \Omega)$.
\end{corollary}
\pf
The sequence of the proposition provides a soft resolution of
$\calo(E)$; and cohomology groups of coherent sheaves on compact spaces
are finite dimensional by the finiteness theorem of Cartan-Serre.
\qed
\vspace{2ex}

Now let $E\downarrow q:\calc\rightarrow S$ be a continuous family
of holomorphic vector bundles over prestable curves. To equip
\[
  {q_1^p}_*E\ :=\ \coprod_s \check L_1^p(C_s;E_s)
\]
with the structure of a Banach bundle we introduce the concept of
quasi-isomorphism for morphisms of Banach bundles:
\begin{defi}\label{quasiiso}\rm
A commutative square of Banach bundles
\[\begin{array}{ccc}
  E&\stackrel{\gamma}{\llongrightarrow} &G\\[10pt]
  \diagl{\alpha}&&\diagr{\beta}\\[10pt]
  F&\stackrel{\delta}{\llongrightarrow}& H
\end{array}\]
is called {\em quasi-isomorphism} (between $\alpha$ and $\beta$, and
between $\gamma$ and $\delta$) if the sequence
\[
  0\longrightarrow E\stackrel{(\alpha,\gamma)}{\longrightarrow}
  F\oplus G\stackrel{\delta-\beta}{\longrightarrow} H\longrightarrow 0
\]
is exact and locally split.
\qed
\end{defi}
Exactness of the sequence can be rephrased by saying that the square
is cartesian ($E\simeq F\oplus_H G$) and cocartesian ($H\simeq (F\oplus G)/
E$). The name is motivated by the fact that for coherent sheaves this is
equivalent to the requirement that $\alpha$ and $\beta$ (or equivalently,
$\gamma$ and $\delta$) induce isomorphisms in cohomology, i.e.\ between
kernels and cokernels.
\vspace{1ex}

We will define the bundle structure on ${q_1^p}_*E$ in such a way that
locally, given a covering $\calu$ as in Section~\ref{cech}, there is a
quasi-isomorphism between $\dbar: {q_1^p}_*E \rightarrow
q^p_*(E\otimes \bar\Omega)$ and the \v Cech complex $\check d:
q_*^{(0)} (\calu,E) \rightarrow q_*^{(1)} (\calu,E)$. In fact, the bundle
structure on ${q_1^p}_*E$ will be defined by an exact sequence as in
Definition~\ref {quasiiso}.
\vspace{1ex}

We begin with a family version of the operator $T$. For the construction
we identify $Z_t\subset\Delta\times\Delta$ with $\Delta\setminus
\Delta_{\tau^2}$ via $\prj_1:(z,w)\mapsto z$, $t=\tau^2 e^{i\theta}
\neq0$. For $\gamma\in L^p(Z_t;\bar\Omega_{Z_t})$ set
\[
  (T_t\gamma)(z)\ :=\ \frac{1}{2\pi i}\int_{\Delta\setminus
  \Delta_{\tau^2}} \frac{dw}{w-z}\wedge\gamma(w)\,,\quad z\in\Delta
  \setminus \Delta_{\tau^2}\, .
\]
It is also useful to view $T_t$ as composition of the trivial extension
of $\gamma$ to $\Delta$ with the fixed operator $T$. For $t=0$ we set
\[
  T_0:=T\amalg(T+R)\,,\quad
  \mbox{with }R(\gamma^+\amalg\gamma^-)=T\gamma^+(0)
  -T\gamma^-(0)\, .
\]
$\{T_t\}_{t\in\Delta}$ is a family of operators from ${\check
L}^p(Z_t;\bar\Omega_{Z_t})$ to $\check L_1^p(Z_t)$ with $\dbar
T_t\gamma =\gamma$.

Now let $E \downarrow q:\calc\rightarrow S$ be given as in Section~\ref
{cech}, i.e.\ as gluing of a family of holomorphically trivialized vector
bundles over a continuous family of circular regions $G_s$ and over $d$
(possibly degenerate) standard cylinders $U_i(s)=Z_{t_i(s)}$.  We will
suppress various transition and coordinate functions for the rest of the
construction. For $\gamma\in \check L^p(G_s;{\bar\Omega}^{\oplus r})$
(viewed as $r$-tupel of differential forms) put
\[
  T^0_s \gamma(z)\ :=\ \frac{1}{2\pi i}\int_{G_s}\frac{dw}{w-z}
  \wedge\gamma(w)\,,
\]
while for $i>0$ we define $T^i_s: \check L^p(Z_{t_i(s)};{\bar\Omega}^{\oplus r})
\rightarrow \check L_1^p(Z_{t_i(s)};\cz^r)$ as above, $r= \rk E$.

We are now in position to write down the maps in our quasi-isomorphism.
For notational reasons we put, for a 1-cochain $(f_{ij})_{i<j}$ and $i>j$:
$f_{ij}:=-f_{ji}$, $f_{ii}=0$, and write just $(f_{ij})_{ij}$ for the cochain.
\[\begin{array}{lrclrcl}
  \Theta:&\hspace{-1ex}
  {q_1^p}_*E&\longrightarrow& q_*^{(0)}(\calu,E)\,,\ &
  \Theta_s(f)&\hspace{-1ex}:=\hspace{-1ex}&
  \Big(f|_{U_i(s)}-T_s^i(\dbar f|_{U_i(s)})\Big)_i\\[10pt]
  \Lambda:&\hspace{-1ex}
  q^p_*(E\otimes\bar \Omega)&\longrightarrow&
  q_*^{(1)}(\calu,E)\,,&
  \Lambda_s(\gamma)&\hspace{-1ex}:=\hspace{-1ex}&
  \Big(T_s^j(\gamma|_{U_j(s)})- T_s^i(\gamma|_{U_i(s)}) \Big)_{ij}
\end{array}\]
To check exactness of the corresponding sequence and for expliciteness
of the trivialization of ${q_1^p}_*E$ that we are heading for, it is useful to
write down explicit splittings. Recall the retraction $\kappa: \calc
\rightarrow C_{s_0}$. Let $\calu_s$ be the covering of $C_s$ induced by
$\calu$. If $\{\rho_i\}$ is a partition of unity subordinate to $\calu_{s_0}$
then $\kappa^*\rho_i$ is a partition of unity for\, $\calu$, that we will
also denote $\rho_i$. A crucial role will be played by the inhomogeneous
Cauchy integral formula (cf.\ e.g.\ \cite[Ch.I,\S5]{vekua}), that states, for a
domain $G\subset\cz$ and $f\in L_1^p (\Delta)$
\begin{eqnarray*}
  &f\ =\ T^G (\dbar f) + H^G(f|_{\partial G})\,,\\
  &\mbox{with}\quad\displaystyle (T^G \gamma)(z)\ =\ 
  \frac{1}{2\pi i}\int_G \frac{dw}{w-z}\,\wedge \gamma \,,\quad
  (H^G\ph)(z)\ =\ \frac{1}{2\pi i}\int_{\partial G} \frac{\ph(w)}{w-z}\, dw\,.
\end{eqnarray*}
Here $\gamma$ stands for an $L^p$-valued (0,1)-form, and $\ph$ for a
continuous function on $\partial G$. The formula yields a
decomposition of $L_1^p$-spaces into holomorphic and essentially
non-holomorphic parts. The formula immediately applies to $U_0(s)$ and to
$U_i(s)$, $i>0$, unless $t_i(s)=0$. In the latter case the
domain splits into two irreducible components and the definition of
$T^i_{t=0}$ has been adjusted by the factor $R$ to make the result
continuous at the node. Accordingly, the same factor $R$ has to be
subtracted from $H^{\Delta^-}$. The possibly adjusted boundary integrals
are denoted $H^i_s$. So $\Theta$ could alternatively be defined by
$\Theta_s(f)= (H^i_s(f|_{\partial U_i(s)}))_i$, observing the identification of
$\partial U_i(s)$ with a union of circles. The maps for the splittings are
\[\begin{array}{crclrcl}
  T:&\hspace{-1ex}
  q_*^p(E \otimes\bar \Omega)&\longrightarrow&
  {q_1^p}_*E\,,&
  T_s(\gamma)&\hspace{-1ex}:=\hspace{-1ex}&
  \sum_i\rho_i\cdot T^i_s(\gamma|_{U_i(s)})\\[5pt]
  H:&\hspace{-1ex}
  q_*^{(1)}(\calu,E)&\longrightarrow& q_*^{(0)}(\calu,E)\,,\ &
  H_s(f_{ij})&\hspace{-1ex}:=\hspace{-1ex}&
  \Big(\frac{1}{2}H^i_s(\sum_k f_{ki}\cdot \rho_k)\Big)_i\\[10pt]
  \Phi:&\hspace{-1ex}
  q_*^{(0)}(\calu,E)&\longrightarrow& {q_1^p}_*E\,,&
  \Phi_s(h_i)&\hspace{-1ex}:=\hspace{-1ex}&
  \sum_i \rho_i\cdot h_i\\[10pt]
  \Psi:&\hspace{-1ex}
  q_*^{(1)}(\calu,E)&\longrightarrow&q^p_*(E\otimes \bar\Omega)\,,&
  \Psi_s(f_{ij})&\hspace{-1ex}:=\hspace{-1ex}&
  \frac{1}{4}\sum_{i,j} f_{ij}\cdot\dbar\rho_i
\end{array}\]
A few remarks: On $U_i$, any $f_{ki}\cdot\rho_k$, a priori only defined on
$U_i\cap U_k$, can be continuously extended to all of $U_i$ by zero,
and this extension is meant in the definition of $H$. Similarly, for the
definition of $\Psi$, $\dbar\rho_i|_{U_{ij}}$ has compact support and is
trivially extended to $C_s$.
\begin{lemma}\label{HPhiPsi}
  With respect to the already defined Banach bundle structures on
  $q^{(0)}_*(\calu, E)$, $q^{(1)}_*(\calu, E)$ and $q^p_*(E\otimes\bar\Omega)$
  the maps $\Lambda$, $\check d$, $\Psi$, $H$ are morphism of Banach
  bundles. The maps $\dbar$, $\Theta$, $T$, $\Phi$ are fiberwise continuous.
\end{lemma}
\pf
In the given local trivializations these maps are even given as
norm-continuous families of operators: This is clear for the \v Cech
coboundary operator and for $\Psi$. For $\Lambda$ this will follow from
Lemma~\ref{T_t}. For $H$ we observe that in view of the Cauchy integral
formula one can conclude a $C^0$-bound on $H^i_s(f_{ki}\rho_k)$
from $C^0$-bounds on $f_{ki}\rho_k$ and on $\dbar(f_{ki}\rho_k)
=f_{ki}\dbar \rho_k$.

Fiberwise continuity of the maps involving ${q_1^p}_*E$ is clear in view
of the estimates for Cauchy integral operators $T^i_s$ on plane domains
(for $T$ and $\Theta$), and in view of estimates for derivatives of
holomorphic functions on compact subsets of their domain in terms of
their sup-norm (for $\Phi$).
\qed

\noindent
The maps defined above fit into a family of split exact sequences of Banach
spaces over the base $W$ as follows. 
\begin{lemma}\label{split_sequence}
  The sequence
  \[\begin{array}{ccccccccc}
    0&\llongrightarrow& {q_1^p}_*E
    &\begin{array}{c} \\[-9pt]
    \scriptstyle (\dbar,\Theta)\\[-6pt]
    \llongrightarrow\\[-10pt]
    \longleftarrow\\[-6pt]
    \scriptscriptstyle T+\Phi
    \end{array}&
    q^p_*(E\otimes\bar \Omega)\oplus q_*^{(0)}(\calu,E)
    &\begin{array}{c} \\[-9pt]
    \scriptstyle \Lambda-\check d\\[-6pt]
    \llongrightarrow\\[-10pt]
    \longleftarrow\\[-6pt]
    \scriptscriptstyle (\Psi,-H)
    \end{array}&
    q_*^{(1)}(\calu,E) &\llongrightarrow& 0
  \end{array}\]
  is (fiberwise) exact and splits with the indicated maps.
\end{lemma}
\pf
The following two equations are trivially checked by direct computation:
\begin{eqnarray*}
  &(\Lambda-\check d)\circ(\dbar,\Theta)\ =\ 
  \Lambda\dbar -\check d \Theta\ =\ 0\,,\\
  &(T+\Phi)\circ(\dbar,\Theta)\ =\ T\dbar+\Phi\circ\Theta\ =\ 
  \id_{{q_1^p}_*E}\,.  
\end{eqnarray*}
For the composition of $\Lambda-\check d$ with $(\Psi,-H)$ we observe
that on $U_i$
\[
  \frac{1}{2}\sum_{k,l} f_{kl}\dbar\rho_k\ =\ 
  \frac{1}{2}\sum_k f_{ki}\dbar\rho_k + \frac{1}{2}\sum_l f_{il}\dbar\rho_i
  \ =\ \sum_k f_{ki}\dbar\rho_k\,,
\]
for, whenever $U_i\cap U_l \neq\emptyset$ then $\dbar\rho_i=-\dbar
\rho_l$. The component of $2(\Lambda\Psi+\check d H)(f_{kl})$ on
$U_i\cap U_j$ can thus be computed as
\begin{eqnarray*}&&
  \frac{1}{2}\sum_{k,l}\Big(T^j (f_{kl}\dbar\rho_k) 
  -T^i (f_{kl}\dbar\rho_k)\Big)
  +\sum_k\Big(H^j(f_{kj}\rho_k)- H^i(f_{ki}\rho_k)\Big)\\
  &=&\sum_k\Big({T^j (f_{kj}\dbar\rho_k)- T^i(f_{ki}\dbar\rho_k)
  -H^i(f_{ki}\rho_k})+ H^j(f_{kj}\rho_k)\Big)\\
  &=& \sum_k(f_{kj}\rho_k- f_{ki}\rho_k)\ =\ f_{ij}-f_{ji}\ =\ 2f_{ij}\,.
\end{eqnarray*}
Finally, we have to show the second half of exactness at the middle. So let
$\gamma\in \check L^p (C_s;E\otimes\bar\Omega)$, $(h_i)\in \check C^0
_{(c)} (\calu;E)$ be such that $\Lambda\gamma -\check d(h_i)=0$. Written
out this means
\[
  T^j\gamma +h_j\ =\ T^i\gamma+h_i\,,
\]
on $U_i\cap U_j$. The candidate for $f\in \check L_1^p(C_s;E)$ with
$\gamma=\dbar f$, $(h_i)=\Theta f$ is
\[
  T\gamma+\Phi(h_i)\ =\ \sum_i \rho_i(T^i\gamma+ h_i)\,.
\]
The previous equation shows that on $U_i$ this is nothing but $T^i\gamma
+h_i$. Hence, on $U_i$ it holds
\begin{eqnarray*}
  \dbar f&=& \dbar(T^i\gamma+h_i)\ =\ \gamma\\
  (\Theta f)_i&=& (T^i\gamma +h_i) - T^i \gamma\ =\ h_i\,.
\end{eqnarray*}
The proof is finished.
\qed

\noindent
By the previous two lemmas we may now {\em define} the Banach bundle
structure on ${q_1^p}_*E$ by the exact sequence above, i.e.\ as kernel of
the morphism $T+\Phi$ of Banach bundles. We will show below (Proposition~
\ref{bundle_structure}) that various
choices made in the construction (covering $\calu$, retraction $\tilde
\kappa$ etc.)  lead to compatible structures.
Together with the previous lemma we obtain
\begin{corollary}
The following squares are quasi-isomorphisms:
\[\begin{array}{ccccccc}
  {q_1^p}_*E&\stackrel{\dbar}{\llongrightarrow}&
  q^p_*(E\otimes\bar\Omega)&\hspace{2cm}&
  {q_1^p}_*E&\stackrel{T}{\llongleftarrow}&
  q^p_*(E\otimes\bar\Omega)\\[10pt]
  \diagl{\Theta}&&\diagr{\Lambda}&&
  \updiagl{\Phi}&&\updiagr{\Psi}\\[10pt]
  q_*^{(0)}(\calu, E)&\stackrel{\check d}{\llongrightarrow}&
  q_*^{(1)}(\calu,E)&&
  q_*^{(0)}(\calu, E)&\stackrel{H}{\llongleftarrow}& q_*^{(1)}(\calu,E)
\end{array}\]
\vspace{-4ex}

\qed
\end{corollary}


\subsection{Regularity properties of ${q_1^p}_*E$}
This section is devoted to the proof of well-definedness and some related,
but more refined regularity questions concerning the Banach bundle structure
on ${q_1^p}_*E$ defined in the previous section.
The following continuity result will be used at
several places in the sequel. It shows that the family $T_t$ is compatible
with retraction maps $\kappa_t$.
\begin{lemma}\label{T_t}
  For $2< p <\infty$ the family of bounded operators
  \[
    (\kappa_t^*)^{-1}\circ T_t\circ\kappa_t^\star: L^p(Z_0;\bar\Omega_{Z_0})
    \longrightarrow  L^\infty(Z_0)\,,
  \]
  is continuous in norm.
\end{lemma}
\pf
To save on typing we suppress trivial extensions of $L^p$ forms. Letting
$\iota_t:\Delta\setminus\Delta_\tau\rightarrow\Delta_\tau\setminus\Delta_{\tau^2}$,
$w\mapsto t/w$ ($t=\tau^2 e^{i\theta}$), one computes for $z=\kappa_t(z_t)\in\Delta^+
\setminus\{0\}$:
\begin{eqnarray*}
  T_t(\kappa_t^\star\gamma)(z_t)&=&\frac{1}{2\pi i}\int_{\Delta\setminus\Delta_\tau}
  \frac{dw}{w-z_t}\wedge\lambda_\tau^\star\gamma^+(w)
  +\frac{1}{2\pi i}\int_{\Delta\setminus\Delta_\tau}
  \Big(\iota_t^*\frac{dw}{w-z_t}\Big)\wedge\lambda_\tau^\star\gamma^-(w)\\
  &=&(T\lambda_\tau^\star\gamma^+)(z_t)+(T\lambda_\tau^\star\gamma^-)(t/z_t)
  -(T\lambda_\tau^\star\gamma^-)(0)\, .\hspace{3cm}(*)
\end{eqnarray*}
Here we have used
\[
  -\frac{t}{w^2}\frac{1}{(t/w)-z_t}\ =\ \frac{1}{w-(t/z_t)}-\frac{1}{w}
\]
to split $\displaystyle\iota_t^*\frac{dw}{w-z_t}$ into two parts. We get
\[
  |T_t(\kappa_t^\star\gamma)(z_t)-(T\gamma^+)(z)|
  \ \le\ |(T\lambda_\tau^\star\gamma^+)(z_t)-(T\gamma^+)(z)|
  +|(T\lambda_\tau^\star\gamma^-)(t/z_t)-(T\lambda_\tau^\star\gamma^-)(0)|\,.
\]
The second term can be estimated by $|t/z_t|^{1-(2/p)}$ times the
$C^{1-(2/p)}$-H\"older norm of $T\lambda_\tau^\star \gamma^-$. Now
$t/z_t$ becomes arbitrarily small with $t$ and
\[
  |\!|T\lambda_\tau^\star \gamma^-|\!|_{C^{1-(2/p)}}\ \le\ 
  |\!|T|\!|_{L^p,C^{1-(2/p)}}\cdot |\!|\lambda_\tau^\star\gamma^-|\!|_{L^p}\,.
\]
The operator norm of $T$ to the H\"older space is bounded in view of the Sobolev
embedding theorem $L_1^p\subset C^{1-(2/p)}$ for $p>2$, while
$|\!|\lambda_\tau^\star\gamma^- |\!|_{L^p} = (1+\tau)^{-1/p} |\!|\gamma^-|\!|_{L^p}$
as one easily checks.
\vspace{1ex}

The first term is a little harder to deal with. Write $\mu$ for Lebesgue measure and
$\gamma^+(w)=f(w) d\bar w$. By a change of variables we obtain
\begin{eqnarray*}
  (T\lambda_\tau^\star\gamma^+)(z_t)&=&
  -\frac{1}{\pi}\int_{\Delta\setminus\Delta_\tau}\frac{\lambda_\tau^*f(w)}{w-z_t} d\mu(w)\\
  &=& -\frac{1}{\pi}\int_{\Delta}\frac{f(w)}{w-z}\cdot
  \frac{w-z}{\lambda_\tau^{-1}(w)-\lambda_\tau^{-1}(z)}
  \cdot|D\lambda_\tau|^{-1} d\mu(w)\,.
\end{eqnarray*}
The relevant term thus takes the form
\[
  \int_{\Delta}\frac{f(w)}{w-z}\Big(
  \frac{w-z}{\lambda_\tau^{-1}(w)-\lambda_\tau^{-1}(z)}\cdot|D\lambda_\tau|^{-1}-1
  \Big)d\mu(w)\,.
\]
In view of the H\"older type estimate $a^{1/p}+b^{1/p}\le 2^{1- \frac{1}{p}}
(a+b)^{1/p}$ for any nonnegative numbers $a,b$, we may split the integral into one
over $\Delta_{\tau^{1/2}}$ and one over $\Delta\setminus \Delta_{\tau^{1/2}}$. This
amounts to consider only functions $f$ with supports in these subsets of the unit disk.

To facilitate computations and readability of formulas let us first
convince ourselves that we may remove the normalizing stretch from $\lambda_\tau$.
In fact, $\lambda_\tau=  m_{(1-\tau)^{-1}}\circ\bar\lambda_\tau$ with $m_a:\Delta
\rightarrow \Delta_a$ the multiplication by $a\in\, ]\,0,1\,[$ and $\bar\lambda_\tau
(r e^{i\ph}) = (r-\tau) e^{i\ph}$. A change of variables in the integral defining $T$
shows that
\[
  m_a^*\circ T\ =\ a\cdot T\circ m_a^\star\,.
\]
So $(m_a^*)^{-1}\circ T\circ m_a^\star$ clearly depends continuously on $\tau$,
for $a=(1-\tau)^{-1}$.
\vspace{1ex}

In the second case of support away from the origin we proceed by estimating
\[
  \Big|\frac{w-z}{\bar\lambda_\tau^{-1}(w)-\bar\lambda_\tau^{-1}(z)}
  \cdot|D\bar\lambda_\tau|^{-1}-1\Big|
  \ \le\ |\!| D\bar\lambda_\tau -\id|\!|_\infty\cdot |\!|\,|D\bar\lambda_\tau|^{-1}\,|\!|_\infty
  + |\!|\,|D\bar\lambda_\tau|^{-1} -1|\!|_\infty\,.
\]
Here we treated the term involving $w-z$ by the mean value theorem and inserted
$|D\bar\lambda_\tau|^{-1}$. The norms are sup-norms, with respect to the operator norm
in case of matrices. A direct computation yields $|D\bar\lambda_\tau
((r+\tau) e^{i\ph})|^{-1}= \frac{r+\tau}{r}$ and, with $r e^{i\ph}=x+iy$
\[
  D\bar\lambda_\tau ((r+\tau)e^{i\ph})\ =\ \frac{r}{r+\tau}\cdot
  \left(\begin{array}{cc} 1+\tau\frac{x^2}{r^3}&\tau\frac{xy}{r^3}\\
  \tau\frac{xy}{r^3}&1+\tau\frac{y^2}{r^3} \end{array}\right)\,.
\]
Thus
\[
  D\lambda_\tau-\id\ =\ (\frac{r}{r+\tau}-1)\cdot\id+\frac{\tau}{r+\tau}
  \left(\begin{array}{cc} \frac{x^2}{r^2}&\frac{xy}{r^2}\\
  \frac{xy}{r^2}&\frac{y^2}{r^2} \end{array}\right)
  \ =\ \frac{\tau}{r+\tau}\cdot
  \left(\begin{array}{cc}-1+ \frac{x^2}{r^2}&\frac{xy}{r^2}\\
  \frac{xy}{r^2}&-1+\frac{y^2}{r^2} \end{array}\right)\,.
\]
This is a matrix with entries bounded by $2$ times $\tau/(r+\tau) \le1/(\tau^{-1/2}+1)$.
Moreover, $|D\bar\lambda_\tau|^{-1}=(r+\tau)/r\le 1+\tau^{1/2}$ and
$|D\bar\lambda_\tau|^{-1}-1=\tau/r\le\tau^{1/2}$. Thus
\[
  C(\tau)\ :=\ |\!| D\bar\lambda_\tau -\id|\!|_\infty\cdot |\!|\,|D\bar\lambda_\tau|^{-1}\,|\!|_\infty
  + |\!|\,|D\bar\lambda_\tau|^{-1} -1|\!|_\infty
  \ \le\ \frac{2}{\tau^{-1/2}+1}\cdot(1+\tau^{1/2})+\tau^{1/2}
\]
tends uniformly to zero with $\tau$. Putting things together,
we obtain for $\gamma^+$ with support away from $\Delta_{\tau^{1/2}}$
\[
  |(T\bar\lambda_\tau^\star\gamma^+)(z_t)-T\gamma^+(z)|\ \le\ 
  \frac{C(\tau)}{\pi}|T\gamma^+(z)|\ \le\ \frac{C(\tau)}{\pi}|\!|T|\!|_{L^p,C^0}
  |\!|\gamma^+|\!|_{L^p}\,,
\]
which is the desired estimate in this case.
\vspace{1ex}

The other case of $f$ with support contained in $\Delta_{\tau^{1/2}}$ is done by a rescaling
argument. We want to show that $T\beta(z)/|\!|\beta|\!|_{L^p}$ tends to zero uniformly
with $\tau$ for $\beta=\gamma^+$ or $\beta=\lambda_\tau^\star\gamma^+$, and $z\in \Delta$.
Note that in any case $\beta$ has support in $\Delta_{\tau^{1/2}+\tau}$, which is contained in
$\Delta_{2\tau^{1/2}}$ for small $\tau$. Let us write $\beta(w) = g(w) d\bar w$.

For $z\in \Delta_{\tau^{1/4}}$ we may pull back by multiplication by $a=\tau^{1/4}$
(smaller than $2\tau^{1/2}$ for small $\tau$) to conclude
\begin{eqnarray*}
  |T\beta(z)|&\le& a\cdot|\!|Tm_a^\star\beta|\!|_{C^0}\ \le\ 
  a\cdot |\!|T|\!|_{L^p,C^0}\cdot |\!|m_a^\star \beta|\!|_{L^p}
  \ =\ a^{1-\frac{2}{p}}|\!|T|\!|_{L^p,C^0}\cdot |\!|\beta|\!|_{L^p}\,.
\end{eqnarray*}
And since $p>2$ the exponent of $a=\tau^{1/4}$ is positive, so this tends uniformly
to zero with $\tau$.

For the remaining case $z\in\Delta\setminus \Delta_{\tau^{1/4}}$ we choose polar
coordinates $se^{i\psi}$ centered in $z$ to give a direct estimate of the integral defining
$T$. Note that since the support of $\beta$ is contained in $\Delta_{2\tau^{1/2}}$ the
distance from $z$ to the support of $\beta$ is at least $\tau^{1/4}-2\tau^{1/2}$,
which for small $\tau$ can be bounded by $\tau^{1/4}/2$. We get
\[
  |T\beta(z)|\ =\ \Big|\int g(s e^{i\psi}) e^{-i\psi}\,ds\,d\psi \Big|
  \ \le\ 2\tau^{-1/4}\int|g(se^{i\psi})|s\,ds\,d\psi\ =\ 2\tau^{-1/4}|\!|g|\!|_{L^1}\,.
\]
Now the H\"older estimate $|\!|g|\!|_{L^1}\le |\!|1|\!|_{L^{1-(1/p)}
(\Delta_{2\tau^{1/2}})}|\!|g|\!|_{L^p} = (4\tau\pi)^{1-\frac{1}{p}}|\!|g|\!|_{L^p}$ together
with $p>2$ gives the final result $|T\beta(z)| \le 8\pi^{1-\frac{1}{p}} \tau^{1/4}
|\!|g|\!|_{L^p}$.
\vspace{2ex}

For $z\in\Delta^-\setminus\{0\}$, $z_t$ and $t/z_t$ change their roles, so the three terms
in ($*$) above converge to $T\gamma^+(0)$, $T\gamma^-(z)$ and $T\gamma^-(0)$
respectively. The rest of the argument works mutatis mutandis. Finally, for
$t\in\Delta^*$ we have to deal with a family of diffeomorphisms $\lambda_\tau$.
Then uniform bounds on the derivatives exist and continuity follows from uniform
continuity of
\[
  \frac{w-z}{\lambda_\tau^{-1}(w)-\lambda_\tau^{-1}(z)}\cdot|D\lambda_\tau|^{-1}
\]
as in the case of $\gamma$ with support away from $\Delta_{\tau^{1/2}}$ above.
\qed
\vspace{1ex}

\noindent
Note that we have given the lemma in the slightly stronger form for $L^p$-spaces
without weights.
\vspace{2ex}

For the derivative of $T_t$ one only expects continuity in the compact-open topology.
\begin{lemma}\label{S_t}
  For $1<p<\infty$ the family of bounded operators
  \[
    (\kappa_t^\star)^{-1}\circ \partial\circ T_t\circ\kappa_t^\star:
    \check L^p(Z_0;{\bar\Omega}_{Z_0})
    \longrightarrow\check L^p(Z_0;\Omega_{Z_0})\,,
  \]
  depends continuously on $t\in\Delta$ in the compact-open topology.
\end{lemma}
\pf
For nonzero $t$, $\partial\circ T_t$ is the restriction of the singular
integral operator $S$ introduced in Section~4.3 to
$\Delta\setminus \Delta_{\tau^2}$, while for $t=0$ we get the operator
$S\amalg S$ on $\Delta\amalg\Delta$. The proof now proceeds similarly
to Lemma~4.12. This time we find
\[
  \iota_t^*\frac{dw}{(w-z_t)^2}\ =\ -\frac{t}{w^2}\frac{dw}{((t/w)-z_t)^2}
  \ =\ -\frac{t}{z_t^2}\frac{dw}{(w-(t/z_t))^2}\,,
\]
so
\begin{eqnarray*}
  \lefteqn{(\kappa_t^\star)^{-1}S_t\kappa_t^\star\gamma|_{\Delta^+} -S\gamma^+\ =\ 
  (\lambda_\tau^\star)^{-1}S\lambda_\tau^\star\gamma^+
  +(\lambda_\tau^\star)^{-1} m_t\iota_t^\star S(\lambda_\tau^\star\gamma^-)
  -S\gamma^+}\hspace{2cm}\\
  &=& (\lambda_\tau^\star)^{-1}S(\lambda_\tau^\star\gamma^+-\gamma^+)
  +((\lambda_\tau^\star)^{-1}-{\rm Id})S\gamma^+\\
  &&+(\lambda_\tau^\star)^{-1}m_t\iota_t^\star S
      (\lambda_\tau^\star\gamma^- -\gamma^-)
  +(\lambda_\tau^\star)^{-1}m_t\iota_t^\star S(\gamma^-)\,,
\end{eqnarray*}
where $m_t:\check L^p(\Delta\setminus\Delta_\tau;\Omega)\rightarrow
\check L^p(\Delta\setminus\Delta_\tau;\Omega)$ is multiplication by
$-t/z^2$ and
\[
  \iota_t^\star:\check L^p(\Delta\setminus\Delta_\tau;\Omega)\longrightarrow
  \check L^p(\Delta_\tau\setminus\Delta_{\tau^2};\Omega)\,,\quad
  fdz\longmapsto(\iota_t^*f)dz\,.
\]
Now the trivial extensions of $\lambda_\tau^\star$, $(\lambda_\tau^\star)^{-1}$
to endomorphisms of $L^p$-forms on $\Delta$ converge to $\rm Id$ in the
compact-open topology (cf.\ the proof of Proposition~4.2).
Together with uniform boundedness of $m_t$ this treats the first three terms.
The last term is the composition of the uniformly bounded operator $S$
and $(\lambda_\tau^\star)^{-1}m_t\iota_t^\star$ that we again view as endomorphism
of $\check L^p(\Delta,\Omega)$ by composing with the restriction map to
forms on $\Delta_\tau\setminus \Delta_{\tau^2}$. Explicitely, the latter map sends
$h(z)dz$ to $(-t/(\lambda_\tau^{-1}(z))^2) h(t/\lambda_\tau^{-1}(z))dz$. Now
a change of variables shows that the $\check L^p$-norm of this form
($z\in \Delta$) can be bounded in terms of the $\check L^p$-norm of
$h\, dz$ restricted to the
annulus $\Delta_\tau\setminus\Delta_{\tau^2}$. On any compact set of $\check
L^p$-functions the latter restriction tends uniformly to zero. Hence the claim.
\qed
\vspace{2ex}

We are now in position to prove independence of choices of the bundle structure on
${q_1^p}_*E$.
\begin{prop}\label{bundle_structure}
  Let $E\downarrow q:\calc\rightarrow S$ be a continuous family
  of holomorphic vector bundles over prestable curves. Then
  ${q_1^p}_*E$ has naturally a structure of Banach bundle over $S$.
  With this structure the maps $\dbar$, $\Theta$, $T$, $\Phi$ from above
  become morphisms of Banach bundles.
\end{prop}
\pf
For $s_0, s'_0\in S$ we have to compare local, fiberwise holomorphic
trivializations $\Pi_s:\check L_1^p(C_{s_0};E_{s_0})\rightarrow{\check
L}_1^p(C_s;E_s)$, $\Pi'_s:\check L_1^p(C_{s'_0};E_{s'_0}) \rightarrow {\check
L}_1^p (C_s;E_s)$ induced by open coverings $\calu$, $\calu'$ of (parts of)
$\calc$, subordinate partitions of
unity $\{\rho_i\}$, $\{\rho'_i\}$, local trivializations of $E|_{U_i}$ and
$E|_{U'_i}$ and admissable retractions (of $E$) $\tilde \kappa$,
$\tilde \kappa'$ covering retractions $\kappa,\kappa'$ of $\calc$.
The trivializations will be supppressed in the notation.

$\Pi_s$ was given as follows: To $f\in \check L_1^p(C_{s_0};E_{s_0})$
associate $((h_i),\gamma) \in \check C^0_{(c)}(\calu_{s_0}; \calo(E_{s_0}))
\oplus \check L^p (C_{s_0};E_{s_0}\otimes \bar\Omega)$, with
$h_i=f-T_{s_0}^i (\dbar f) =H_{s_0}^i(f)$, $\gamma=\dbar f$. Transport
$h_i$ to $s$ by a Cauchy type integral of $\tilde \kappa_s^* h_i$ over
$\partial G_s$, as in Lemma~\ref{Cech-zero-cycles}. This defines $h_i(s)$,
while $\gamma(s) := \tilde \kappa_s^\star\gamma$. The pair $((h_i(s)), \gamma(s))
\in \check C^0_{(c)} (\calu_s;\calo(E_s)) \oplus \check L^p(C_s;E_s\otimes
\bar\Omega_{C_s})$ corresponds to
\[
  f_s\ =\ \Pi_s(f)\ =\ \sum\rho_i(s) \Big(h_i(s)+T^i_s\gamma(s)\Big)\ \in
  \check L_1^p(C_s;E_s)
\]
To go back by $(\Pi'_s)^{-1}$ we first have to decompose $f_s$ into
\[
  \gamma'(s) =\dbar_{E_s}f_s= \gamma(s) +\sum\Big(h_i(s)+
  T^i_s\gamma(s)\Big) \,\dbar \rho_i(s)
\]
and the holomorphic cochain $({H'_s}^i (f_s))_i$. Transporting back to $s'_0$
we obtain the $(0,1)$-form $\gamma'_s:=({\mbox{$\tilde \kappa'_s$}^\star})^{-1}
\gamma'(s)$ and the holomorphic \v Cech 0-cochain $h'_{i,s}:= {H'_{s'_0}}^{\!\!\!i}
({\mbox{$\tilde\kappa'$}_{s'}^*}^{-1} {{H'_s}^i} (f_s) )$. Thus
\[
  (\Pi'_s)^{-1}\circ\Pi_s(f)\ =\ \sum\rho'_i(s'_0)\cdot\left
  (h'_{i,s}+{T'_{s'_0}}^{\!\!\!i}\gamma'_s\right)\, .
\]
By Remark~\ref{cp-open-top} we have to show that this family of maps is
continuous in the compact-open topology. It suffices to prove continuity
of the family of maps parametrized by $s$ and sending $((h_i),\gamma)$ to
$((h'_{i,s}),\gamma'_s)$. We begin with $\gamma'_s$, as map from $\check L^p
(C_{s_0};E\otimes\bar\Omega)$ to $\check L^p (C_{s'_0};E\otimes\bar\Omega)$.
Written out this is
\[
  (\mbox{$\tilde\kappa'_s$}^\star)^{-1}\tilde\kappa_s^\star\gamma+
  {\textstyle \sum_i}\Big( (\mbox{$\tilde\kappa'_s$}^*)^{-1} h_i(s)
  +(\mbox{$\tilde\kappa'_s$}^*)^{-1} T^i_s \tilde\kappa_s^\star\gamma\Big)\cdot
   ({\kappa'_s}^\star)^{-1}\dbar\rho_i(s)\,.
\]
Continuity of the first term in the compact-open topology follows as in the proof
of Proposition~\ref{Lptriv}. The term involving $h_i(s)$ can be written
$(\mbox{$\tilde\kappa'_s$}^\star)^{-1} (h_i(s)\cdot\dbar\rho_i(s))$. Continuity of
this expression follows from: uniform estimates on $\dbar\rho_i(s)$ with respect to
natural holomorphic parameters (Lemma~\ref{retraction}); uniform estimates on
$\nabla H^i_s(f)$ on compact subsets of $U_i(s)$ (!) in terms of the sup-norm on
$\partial U_i(s)$ (here: on $\supp \rho_i(s)$); uniform equicontinuity
of $h_i$ varying in a compact subset of $C^0(\partial U_i(s_0))$ by Arzela-Ascoli.
The term involving $T^i_s$ is continuous by Lemma~\ref{T_t}, together with
continuity of $(\mbox{$\tilde\kappa'_s$}^\star)^{-1} \circ \tilde\kappa_s^\star$ (again
as in Lemma~\ref{retraction}), and the already mentioned uniform estimates on
$\dbar \rho_i(s)$. (By reexaming the proof of Lemma~\ref{T_t} one could even
verify uniform continuity here.)

For $h'_{i,s}$ we get
\[
  {H'_{s'_0}}^{\!\!\!i}\Big((\mbox{$\tilde\kappa'_s$}^*)^{-1} {H'_s}^i
  {\textstyle \sum_j}\rho_j(s)(h_j(s)+T^j_s\gamma(s))\Big)\,.
\]
From the inhomogeneous Cauchy integral formula together with the estimates for
$T$ one can bound the $\check L_1^p$-norm of any boundary integral $H^G f$, for
any function $f$ on a plane region $G$ and also on $Z_0$, by the $\check L_1^p$-norm
of $f$ on a neighbourhood $N\subset \overline{G}$ of $\partial G$. We thus
only have to show continuity as map to $\check L_1^p$ away from the nodes
of the term inside the big brackets. We write the term containing $h_j(s)$ as
$((\mbox{$\tilde\kappa'_s$}^*)^{-1}{H'_s}^i \mbox{$\tilde\kappa'_s$}^*)\circ
(\mbox{$\tilde\kappa'_s$}^*)^{-1}(\rho_j(s) h_j(s))$.
Again from interior estimates for the gradient of $h_j(s)$ in terms of the sup-norm
of $h_j$ and from uniform estimates for the derivative of $(\kappa'_s)^{-1}$,
one deduces continuity in the compact-open topology of the term to the right
of the composition symbol. The term to the left of the composition symbol is
continuous in the compact-open topology by what we have just said about
$\check L_1^p$-estimates for $H^G$. An analogous argument together with the
estimates for $T$ cover the term with $T^j_s\gamma(s)$.
\qed
\vspace{1ex}

It is important to understand the regularity properties of the Banach bundle
structure on ${q_1^p}_*E$. We will study two cases: Compatibility with the
Banach bundle structure on $q^p_*E$ for any $p$ (actually, for $p=\infty$);
and continuity of fiberwise differentiation with respect to a continuous
family of {\em holomorphic} connections $\nabla$ on $E$. The latter notion means
that in a local representation of $E$ as continuous family of holomorphically trivialized
vector bundles over $\calu$, $\nabla$ is of the form $d+A_s(z)dz$ or
(for $i>0$) $d+A_s(z,w) dz +B_s(z,w) dw$, with $A_s, B_s$ continuous families
of endomorphisms of $E_s$. So we require the $(0,1)$-part of $\nabla$ to induce
the $\dbar$-operator $\dbar_E$ of $E$ on each curve.

The result is:
\begin{prop}\label{localtriv}
  Let $E\downarrow q:\calc\rightarrow S$ be a continuous family of
  holomorphic vector bundles over prestable curves and $2< p<\infty$.
  \begin{enumerate}
  \item
    Let $W\subset S$ be an open set, $\Pi: W\times \check L_1^p(E|_{C_0})
    \rightarrow {q_1^p}_*E$ be a local trivialization as given in the last section,
    and $\kappa: E|_{q^{-1}(W)} \rightarrow E|_{C_0}$ the corresponding
    admissable retraction. Then
    \[
      (\kappa_s^*)^{-1}\circ\Pi_s:\check L_1^p(E|_{C_0})
      \ \longrightarrow\ L^\infty(E|_{C_0})\,,\ s\in W\,,
    \]
    is uniformly continuous.
  \item
    The inclusion ${q_1^p}_*E\rightarrow q^\infty_*E$ is a morphism of Banach
    bundles.
  \item
    Let $\nabla_s=\dbar_{E_s}+\partial_{E_s}$ be a continuous family
    of holomorphic connections on $E$. Then the corresponding family of operators
    \[
      \nabla:{q_1^p}_*E\ \longrightarrow\
      q^p_*(E\otimes(\Omega\oplus\bar\Omega))
    \]
    is a morphism of vector bundles. Moreover, in the given local trivializations
    \[
      \dbar:\ {q_1^p}_*E\ \longrightarrow\
      q^p_*(E\otimes\bar\Omega)
    \]
    is even a norm-continuous family of operators.
  \qed
  \end{enumerate}
\end{prop}
\pf
For (1) we reexamine the proof of Proposition~\ref{bundle_structure} with
$\kappa'=\kappa$. In the notation of this proposition we have to proof
norm continuity of
\[
  {\kappa^*_s}^{-1}f_s\ =\ \sum\rho_i\Big({\kappa_s^*}^{-1} h_i(s)
  +{\kappa_s^*}^{-1}T^i_s\gamma(s)\Big)\,,
\]
as map from $\check L^p( C;E\otimes\bar\Omega) \oplus \check C^0_{(c)}(\calu;
\calo(E))$ to $L^\infty (C;E)$. Now estimates for higher derivatives
of the Cauchy integrals $H^G$ on compact subsets in $G$ lead
to continuity of $(\kappa_s^*)^{-1} h_i(s)$ in the sup-norm on $\supp
\rho_i$. By the explicit form of $h_i(s)$ near the nodes (cf.\ the proof of
Proposition~\ref{Cech-zero-cycles}) the argument holds for $i>0$ too.
Note that for this step we need $\kappa'=\kappa$ to conclude uniform
continuity in $L^\infty$ on $\partial G$. The term containing
$T^i_s$ is taken care of by Lemma~\ref{T_t}.
Claim (2) is immediate from (1).

As for statement (3), we choose local trivializations of ${q_1^p}_*E$ and $q^p_*E$
built on the same retraction $\kappa$. For $f\in \check L_1^p(C_0;E_0)$,
$\gamma=\gamma(f)$, $h_i=h_i(f)$ we obtain
\[
  (\kappa_s^\star)^{-1}\dbar_{E_s}\Pi_s f\ =\ \gamma+
  (\kappa_s^\star)^{-1}\sum(h_i(s)+T_s^i\kappa_s^\star\gamma)
  \dbar\rho_i(s)\, .
\]
This has the same form as the term discussed for (1), except that $\rho_i$ is
replaced by $\dbar\rho_i$ and a new term $\gamma$ is added.

Doing the same thing with $\partial_E$ turns $\gamma$ into the more
serious term
\[
   (\kappa_s^\star)^{-1}\sum\rho_i(s)\left(\partial_{E_s}h_i(s)
  +\partial_{E_s}T^i_s\kappa_s^\star\gamma\right)\, .
  \quad(*)
\]
This is very similar to the expression for $h'_{i,s}$ in the last paragraph of the
proof of Proposition~ \ref{bundle_structure}. So $(\kappa_s^\star)^{-1} \rho_i(s)
\partial_{E_s} h_i(s)$ is continuous in the compact-open topology by estimates for
the derivatives of $(\kappa_s^*)^{-1}$, $\rho_i(s)$ and $h_i(s)$ away from
arbitrarily small neighbourhoods of the nodes
and on the support of $\rho_i(s)$. The term involving $\partial_{E_s}
\circ T_s^i$ is continuous by Lemma~\ref{S_t}.
\qed


\section{Banach orbifold of maps}
\label{section_maporbifold}
The object of this chapter is to construct the weak Banach orbifold
into which the space of pseudo holomorphic curves embeds as locally compact
subspace. Some of the results are summarized in the following theorem.
For $2g+k<3$ we set $\calm_{g,k}= \{{\rm pt}\}$.
\begin{theorem}\label{orbifold}
  Let $(M,J)$ be an almost complex manifold and $p>2$. Then the
  space $\calc(M;p)$ of stable parametrized complex curves of Sobolev class $L_1^p$
  has naturally a structure of Banach orbifold with local groups
  $G_{(C,{\bf x},\ph)} =\Aut(C,{\bf x},\ph)$. The orbifold topology is finer
  than the $C^0$-topology introduced in Section~\ref{C(X)_topology}.\\[2pt]
  \begin{minipage}{11cm}
  The set of (isomorphism classes of) pairs $[(C,{\bf x},\ph),P]$ with $P\in C$
  has also a structure of Banach orbifold. This Banach orbifold, the
  \underline{universal curve} over $\calc(M;p)$,
  comes with a \underline{universal morphism} $\Phi$ to $M$ and
  with a forgetful morphism $\pi$ to $\calc(M;p)$.
  The underlying set-theoretic maps are evaluation of $\ph$ at $P$
  and forgetting the point $P$ respectively.  Another morphism $p$ maps $\calc(M;p)$
  to $\bigcup \calm_{g,k}$ by sending $(C,{\bf x},\ph)$ to the stabilization of
  $(C,{\bf x})$ as Deligne-Mumford stable curve. Restricted to the open suborbifold
  $\calc_k (M;p)$ of $k$-marked curves there exists a tuple $\underline{\bf x}$ of
  $k$ sections of $\pi$ set-theoretically mapping $(C,{\bf x},\ph)$ to $\bf x$.
  The orbifold structures and morphisms do not depend on $J$.
  \end{minipage}\hspace{2ex}
  \begin{minipage}{3.5cm}
  \vspace{-3ex}
  
  \[\begin{array}{ccc}
    \Gamma(M;p)&\hspace{-2ex}\stackrel{\Phi}{\longrightarrow}& M\\[2pt]
    {\scriptstyle \pi}\!\downarrow\phantom{{\scriptstyle q}\!}\\[2pt]
    \calc(M;p)\\[2pt]
    {\scriptstyle p}\!\downarrow\phantom{{\scriptstyle p}\!}\\[2pt]
    \bigcup_{g,k}\calm_{g,k}
  \end{array}\]
  \vspace{1ex}
  
  \end{minipage}
\end{theorem}
We split the construction of local uniformizers and the necessary verifications
into several steps.


\subsection{Holomorphic structures on $\ph^*T_M$}
\label{hol_structure}
Let $q:\calc\rightarrow S$ be a continuous family of prestable curves, $\ph:C
\rightarrow M$ a map of class $\check L_1^p$ from the central
fiber $C=q^{-1}(0)$ and $\kappa: \calc\rightarrow C= q^{-1}(0)$ an admissable
retraction as in Lemma~\ref{retraction}. Our charts
will be based on an identification of $\check L_1^p(C;\ph^*T_M)$ with
$\check L_1^p(C_s,(\ph\circ\kappa_s)^*T_M)$ for any $s\in S$.
To apply the method of Chapter~\ref{chapter_dbar} it is crucial to equip
$\kappa_s^*\ph^*T_M$ with a continuously varying family of holomorphic
structures. This will be done in the present section.

We choose a (complete for convenience) Riemannian metric $\rho$ on $M$.
Independence of this choice will be shown below.
Then $\ph^*T_M$ inherits the connection $\nabla^\ph:=\ph^*\nabla$ from
the Levi-Civit\`a connection $\nabla$ on $(M,\rho)$. We define the $\dbar$
operator on $\ph^*T_M$ as $\ph^*J$-linear part of the natural linearization of $\ph
\mapsto\dbar_J\ph$ to be introduced in the next chapter. Explicitely, we put,
in real notation, for $v\in \check L_1^p(C;\ph^*T_M)$, $\xi\in T_C$,
$\eta=j(\xi)$, $j$ the complex structure of $C$:
\[
  \dbar_\xi v:=\frac{1}{4}\left(\nabla^\ph_\xi v+J\nabla^\ph_\eta v
  +(\nabla_v J)\partial_J\ph(\eta)-J\nabla^\ph_\xi (Jv)+\nabla^\ph_\eta (Jv)
  -J(\nabla_{Jv}J)\partial_J\ph(\eta)\right)\,,
\]
where we introduce $\partial_J\ph:=\frac{1}{2}(D\ph-J\circ D\ph\circ j)$.
Note that $(\nabla_v J)\circ J=-J\circ\nabla_v J$ such that application of $\nabla_v J$
turns $\ph^*T_M$ valued $(1,0)$ forms into $\ph^*T_M$ valued $(0,1)$ forms.
The only thing important for the moment is the fact that this is the $(0,1)$ part
of a complex linear connection on $\ph^*T_M$ (with complex structure $\ph^*J$)
with $L_1^p$ coefficients, and so makes $\ph^*T_M$ a holomorphic vector bundle
as follows. Our choice is suggested by the fact that we want to study the non-linear
$\dbar_J$-equation later on.
\begin{prop}\label{L1p_dbar}
  Let $(\Sigma,j)$ be a Riemann surface and $E$ a rank $r$ vector bundle of class
  $L_1^p$. Let also be given a complex linear operator
  \[
    \dbar_E:L_1^p(\Sigma;E)\longrightarrow L^p(\Sigma;E\otimes\bar\Omega)
  \]
  satisfying $\dbar_E(f\cdot v)=v\otimes(\dbar f)+f\cdot\dbar_E v$.

  Then the sheaf $\Sigma\supset U\mapsto\calo(E)(U):=\{v\in L_1^p(U;E)
  \mid\dbar_E v=0\}$ is analytic and locally free of rank $r$.
\qed
\end{prop}
The proof is by a nice application of the implicit function theorem to find local holomorphic
frames, cf.\ \cite{holizan} for the line bundle case and \cite{ivshev} for the general case.

It is straightforward to check that the holomorphic structures on $\kappa_s^*\ph^*T_M$
make up a continuous family of holomorphic vector bundles over $C_s$. The theory
developed in Chapter~2 thus does in fact apply. We get a family of isomorphisms
\[
  \Pi_s:\check L_1^p(C_0;\ph^*T_M)\ \longrightarrow\
  \check L_1^p(C_s;\kappa_s^*\ph^*T_M)\, .
\]


\subsection{Local uniformizers at stable domains}
\label{section_uniformizer}
Let $(C,{\bf x},\ph)\in\calc(M;p)$ and assume $(C,{\bf x})$ stable for the time being. Let
$q:\calc\rightarrow S$, $\underline{{\bf x}}:S\rightarrow\calc^k$, be an
analytically universal deformation of $(C,{\bf x})$ as discussed in
Section~\ref{section_deformation}. We write $C_s:=q^{-1}(s)$
and ${\bf x}_s=\underline{{\bf x}}(s)$. We have already observed that such a
family also has a description according to Definition~\ref{cont_curves}. Recall also
the admissable retraction $\kappa:\calc\rightarrow C_0=C$. Moreover, under the
presence of non-trivial automorphisms  of $(C,{\bf x},\ph)$ we choose $q$ and
$\kappa$ equivariant under the action of the finite group $\Aut(C,{\bf x})
\supset \Aut(C,{\bf x},\ph)$. We want to model local
uniformizing systems at $(C,{\bf x},\ph)$ on $S\times\check L_1^p(C;\ph^*T_M)$,
where $S$ is naturally viewed as open neighbourhood of the origin in the tangent
space $\hat T_{\calm_{g,k},(C,{\bf x})}$ of the differentiable
orbifold $\calm_{g,k}$ at $(C,{\bf x})$. 

Next we define the structure maps. The obvious choice at $(C,{\bf x},\ph)$ is
\[
  q_{(C,{\bf x},\ph)}:S\times \check L_1^p(C;\ph^*T_M)\longrightarrow\calc(M;p)\,,
  \quad(s,v)\longmapsto(C_s,{\bf x}_s,\ph(s,v))
\]
with $\ph(s,v)\in \check L_1^p(C_s;M)$ defined by
\[
  \ph(s,v)(z)\ :=\ \exp_{\ph(\kappa_s(z))}^\rho\Pi_s(v)(z)\,,
\]
$\exp_P^\rho:T_P M\rightarrow M$ the exponential map of the Riemannian manifold
$(M,\rho)$ at $P$.

It remains to define an action of the automorphism group
$\Aut(C,{\bf x},\ph)$ on (an open subset of) $S\times \check L_1^p (C,\ph^*T_M)$
such that $q_{(C,{\bf x},\ph)}$ descends to an injection of the quotient space.
$q$ and $\kappa$ being compatible with the action of $\Aut(C,{\bf x},\ph) \subset
\Aut(C,{\bf x})$ the action naturally lifts to $\ph^*T_M$ and hence to
$\check L_1^p(C;\ph^*T_M)$. We force equivariance of $\Pi$ through replacement by
\[
  \check L_1^p(C;\ph^*T_M)\ni v\longmapsto\frac{1}{|\Aut(C,{\bf x},\ph)|}
  \sum_{\Psi\in\mbox{\scriptsize Aut}(C,{\bf x},\ph)}\Psi_s^*\Pi_{\Psi(s)}(\Psi_0^{-1})^*v
\]
at $s\in S$. This has $\Pi_0=\id$ and thus defines a trivialization for $s$ sufficiently
close to $0\in S$, and fulfills $\Pi_{\Psi(s)}\circ(\Psi_0^{-1})^*=(\Psi_s^{-1})^*\circ\Pi_s$
for any $\Psi\in\Aut(C,{\bf x},\ph)$. Then if $\Psi\in\Aut(C,{\bf x},\ph)$ we evidently have
for any $(s,v)$, $\ph(\Psi(s),\Psi(v))\circ\Psi\ =\ \ph(s,v)$
with $\Psi(v):=(\Psi_0^{-1})^*v$. So $\Psi_s:=\Psi|_{C_s}:C_s\rightarrow C_{\Psi(s)}$ induces
the identification of $q_{(C,{\bf x},\ph)}(\Psi(s),\Psi(v))$ with $q_{(C,{\bf x},\ph)}(s,v)$.

The converse statement, that $q_{(C,{\bf x},\ph)}(s,v)= q_{(C,{\bf x},\ph)}(s',v')$
implies the existence of $\Psi\in \Aut(C,{\bf x},\ph)$ with $(s',v')= (\Psi(s), \Psi(v))$
may only be true locally. Namely, for any $\Psi\in\Aut(C,$ ${\bf x})\setminus
\Aut(C,{\bf x},\ph)$ it holds
\[
  d(\Psi):=\sup_{z\in C}d(\ph(z),\ph\circ\Psi(z))>0\, .
\]
By finiteness of $\Aut(C,{\bf x})$ (we still assume $(C,{\bf x})$ stable) the infimum
$d_\ph$ over all such $\Psi$ is positive. In case $\Aut(C,{\bf x})=\Aut(C,{\bf x},\ph)$
we put $d_\ph=\infty$. Let $\mbox{inj}(M,\rho,\ph)$ be the injectivity radius of
$(M,\rho)$ along $\ph$, i.e.\ the minimum of the injectivity radii of $(M,\rho)$
at $\ph(z)$ for any $z\in C$.  Taking $S$ sufficiently small we may restrict to an
open set $S\times V\subset S\times \check L_1^p(C; \ph^*T_M)$ of pairs $(s,v)$ with
\[
  \hspace{2cm}\sup_{z\in C_s}d(\ph\kappa_s(z),\ph(s,v)(z))
  <\frac{1}{2}\min\{d_\ph,\mbox{inj}(M,\rho,\ph)\}\,.\hspace{2cm}(*)
\]
Now assume
$q_{(C,{\bf x},\ph)}(s,v)\simeq q_{(C,{\bf x},\ph)}(s',v')$ as parametrized
complex curves on $M$. This means that there exists an isomorphism
$\Psi_s:(C_s,{\bf x}_s)\rightarrow(C_{s'},{\bf x}_{s'})$ with $\ph(s,v)=\ph(s',v')\circ\Psi_s$.
By ($*$) this only happens if $\Psi\in\Aut(C,{\bf x},\ph)$, where $\Psi$ is the
(unique) automorphism of the central fiber inducing $\Psi_s$. By equivariance of
$\kappa_s$, $\ph(s',v')\circ\Psi$ maps $z\in C_s$ to
\[
  \exp_{\ph\circ\kappa_{s'}\circ\Psi_s(z)}(\Pi_{s'}v')(\Psi_s(z))\ =\ 
  \exp_{\ph\circ\kappa_{s}(z)}(\Pi_{s'}v')(\Psi_s(z))\, .
\]
Comparison with $\ph(s,v)(z)$ shows $\Pi_sv=\Psi_s^*(\Pi_{s'}v')$, which by
equivariance of the $\Pi_s$ equals $\Pi_s(\Psi_s^*v')$. Thus in fact
$(s',v')=(\Psi(s),\Psi(v))$ with $\Psi\in\Aut(C,{\bf x},\ph)$.

Note that in this description the action is manifestly complex linear.


\subsection{Rigidification}\label{sect_rigidification}
So far we have assumed $(C,{\bf x})$ stable. Let us now include nonstable
components. In this case there is (a germ of) a nontrivial holomorphic action
of $\Aut^0(C,{\bf x})$ on $q:\calc\rightarrow S$ inducing (a germ of) an action
$(s,v)\mapsto\Psi(s,v)$ on $S\times V$ characterized by
\[
  \ph(\Psi(s,v))\ =\ (\Psi^{-1})^*\ph(s,v): C_{\Psi(s)}\longrightarrow M\, .
\]
We want to take a smooth slice to this action through $(s,v)=(0,0)$ and show that
the slice is in fact the topological quotient of the equivalence
relation generated by the germ of the action of $\Aut^0(C,{\bf x})$. We first
discuss the case where nevertheless $\Aut(C,{\bf x},\ph)$ is trivial.
\vspace{1ex}

\noindent
{\em The rigidification construction.}\\
To fix notations let $(D_1,x_1),\ldots,(D_{b'},x_{b'})$, $(D_{b'+1},(x_{b'+1},y_{b'+1})),
\ldots,(D_b,(x_b,y_b))$ be the normalizations of the unstable
components of $(C,{\bf x})$ with their special points.
According to Lemma~\ref{discrete_fibers} there exist points $P\neq P'$ in
$D_i\setminus \{x_i\}$, that are isolated in their $\ph$-fibers $\ph^{-1}(\ph(P))$,
$\ph^{-1}(\ph(P'))$. Let $V,V'\subset D_i\setminus\{x_i\}$ be disjoint open neighbourhoods
of $P$, $P'$ with $V\cap \ph^{-1}(\ph(P))= \{P\}$, $V'\cap \ph^{-1}(\ph(P'))
= \{P'\}$. By properness of $\ph$ the preimages of a neighbourhood basis of $\ph(P)$,
$\ph(P')$ form a neighbourhood basis of $\ph^{-1}(\ph(P))$, $\ph^{-1}(\ph(P'))$
respectively. We may thus choose bump functions $\rho_i$, $\rho'_i$ on $M$ with
$\rho_i(\ph(P)) =\rho'_i(\ph(P'))=1$ and such that $\ph^*\rho_i|_V$, ${\ph'}^*\rho'_i|_{V'}$
have compact support in $V$, $V'$ respectively. Let $z_i=z_{b+i}$ be a holomorphic function
defined in a neighbourhood in $\calc$ of the closure of $V\cup V'$ and
restricting to a {\em linear} coordinate on its domain of definition in $D_i\subset C_0$.
So $(a,b)\in \Aut(D_i,{\bf x}_i)\subset \cz^* \ltimes \cz$ acts by
\[
  z_i\longmapsto a z_i+b\,,\quad\mbox{\rm on $D_i$.}
\]
We also require $(q,z_i)$ to yield an isomorphism of the domain of
definition of $z_i$ with $S\times \Delta$. This can be achieved by first
taking the domain of definition of $z_i$ on the central fiber to contain
a sufficiently big ball, then restricting to the level set of this ball by shrinking
$S$ and finally rescaling of $z_i$. Put $W_i= z_i^{-1}(z_i(V))$ and
$W_{b+i}= z_i^{-1}(z_i(V'))$. Again, $S$ can be chosen so small that $W_i\cup W_{b+i}$
are disjoint from the deformation of $x_i$ (the image of $\underline{x_i}: S\rightarrow
\calc$).

Next, by choosing the supports of $\rho_i$, $\rho'_i$
sufficiently small we can achieve that $\ph^*\rho_i|_V$ and $\ph^*\rho'_i|_{V'}$ have
different centers of gravity:
\[
  \int_V z_i\cdot\ph^*\rho_i\,d\mu\ \mbox{\LARGE /} \int_V \ph^*\rho_i\,d\mu\ \neq\
  \int_{V'} z_i\cdot\ph^*\rho'_i\,d\mu\ \mbox{\LARGE /} \int_{V'} \ph^*\rho'_i\,d\mu\,.
\]
In fact, the centers of gravity converge to $z_i(P)$ and  $z_i(P')$ with shrinking supports
respectively.

We will also need soon that the trivialization $\Pi$ equals $\kappa^*$
for functions with support in $W_i\cup W_{b+i}$. To this end we take the domain
of definition of $z_i$ as one connected component of
$U_0$, one of the sets in the covering of $\calc$ used in the \v Cech-construction of
Chapter~\ref{chapter_dbar}. $\kappa$ can be chosen on the domain of definition of $z_i$
via the level sets of $z_i$. Compatibility of $\kappa$ with $\Pi$ follows then by construction,
at least for functions with support away from the intersection locus $U_0\cap U_j$, $j>0$.
By choosing the $U_j$ appropriately we may assume that this is the case for
$W_i\cup W_{b+i}$.

For the bubble components with two special points
$D_{b'+1},\ldots, D_b$ we choose just one bump function $\rho_i$ and one set
$W_i$, but we have to take everything also away from (the deformation of) $y_i$.
By shrinking $S$ and $V$ we may also assume that $(\ph(s,v)|_{W_i})^*\rho_i$ and
$(\ph(s,v)|_{W_{b+i}})^*\rho_i$ have compact supports in $W_i$ respectively $W_{b+i}$.
In fact, in view of Lemma~\ref{localtriv},2 this is an open condition.

The rigidifying data will be obtained as center of gravity of the pull-back by $\ph(s,v)$
of the bump functions, that we conveniently assemble in an
auxiliary map $\Lambda:S\times V\rightarrow\cz^{b+b'}$ defined as follows:
\[
  \Lambda(s,v)\ =\ \left(\int_{C_s\cap W_i}z_i\cdot\ph(s,v)^*
  \rho_i\,d\mu(z_i)\ \mbox{\LARGE /}\ 
  \int_{C_s\cap W_i}\ph(s,v)^*\rho_i\,d\mu(z_i)\right)_{i=1,\ldots,b+b'}\, .
\]
One can easily check that the integrands as elements of $L^p$ spaces depend
differentiably on $s$ and $v$. This shows differentiability of $\Lambda$ in a
neighbourhood of $(0,0)\in S\times V$. We define the slice as $\Lambda^{-1}
(\lambda_0)\subset S\times V$, $\lambda_0:=\Lambda(0,0)$.

To obtain a chart for $\Lambda^{-1}(\lambda_0)$ we apply the implicit function
theorem to $\Lambda$. Note that $\Lambda$ factorizes over the restriction map
$S\times V\rightarrow L_1^p(\bigcup W_i\cap C;\ph^*T_M)$ and hence does not
depend on $S$ (by choice of $\kappa$ and $\Pi$!). In particular, $\ker D\Lambda(0,0)$
is of the form $S\times\bar V$ with $\bar V\subset V$ of codimension $b+b'$. Let
$N\subset V$ be a complementary subspace to $\bar V$ which we identify with
$\cz^{b+b'}$ via $D\Lambda(0,0)$ without further notice. Then there exists a map
\[
  \nu:S\times\bar V\times N\longrightarrow V\quad\mbox{with}\quad
  \Lambda\Big(s,\nu(s,v,\lambda)\Big)=\lambda\quad\mbox{for }(s,v,\lambda)
  \mbox{ near }(0,0,\lambda_0)\,,
\]
and such that $(\prj_S,\nu):(s,v,\lambda)\mapsto(s,\nu(s,v,\lambda))$ is a local
diffeomorphism from $S\times\bar V\times N$ to $S\times V$. $\gamma:=\prj_S
\times\nu|_{S\times\bar V\times\{\lambda_0\}}:S\times\bar V\rightarrow S\times V$
provides a chart for $\Lambda^{-1}(\lambda_0)$.
\vspace{1ex}

\noindent
{\em Slice theorem.}\\
To actually prove that $\Lambda^{-1}(\lambda_0)$ is a local slice it is important that the
action of $\Aut^0(C,{\bf x})$ on $\Lambda$ at $s=0$ can be derived explicitely:
Let $A:\Aut^0(C,{\bf x})\rightarrow\mbox{Aff}(\cz^{b+b'})$ be defined
componentwise by the affine linear action of $\Psi$ on the local coordinates $z_i$
as given above. It holds
\[
  \Lambda(\Psi\cdot(0,v))\ =\ A_\Psi\circ\Lambda(0,v)\,,
\]
for $\Psi$ sufficiently close to the identity and $v$ having small $L^\infty$-norm
(so that the support of $\Psi^*\ph(0,v)^*\rho_i|_{W_i\cap C_0}$ is still compact).

By the explicit form of $A$ and because at $(s,v)=(0,0)$ our two centers of
gravity on the unstable components with only one special point have been
carefully chosen rigidifying, the expression $\Lambda(\Psi(s,v)) -\lambda_0
\in \cz^{b+b'}$ has invertible differential in the direction of $\Psi$. By the implicit
function theorem we thus get a differentiable map $\Psi: S\times V \rightarrow 
\Aut^0(C,{\bf x})$, $(s,v)\mapsto \Psi_{s,v}$ with the characterizing property
\[
  \Lambda(\Psi_{s,v}\cdot s, \Psi_{s,v}\cdot v)\ =\ \lambda_0\,.
\]
Hence
\[
  H:\ (\Psi,s,\bar v)\longmapsto\Psi\cdot \gamma(s,\bar v)
\]
from a neighbourhood of $(\id,0,0)\in\Aut^0(C,{\bf x})\times S\times\bar V$ to
$S\times V$ is a local homeomorphism at $(\id,0,0)$. Its inverse maps
$(s,v)\in S\times V$ to $(\Psi_{s,v}, \gamma^{-1}(\Psi_{s,v}\cdot s,
\Psi_{s,v}\cdot v))$. Note that while $H$ is only continuous, the quotient map
\[
  S\times V\ \longrightarrow\ S\times \bar V\,,\quad
  (s,v)\ \longmapsto\ \gamma^{-1}(\Psi_{s,v}\cdot s,\Psi_{s,v}\cdot v)
\]
is in fact differentiable.
\vspace{1ex}

\noindent
{\em Automorphisms.}\\
If $\Aut(C,{\bf x},\ph)$ is nontrivial a little extra care is needed. First of all we
fix an extension to $\calc$ of the action of the finite subgroup $\Aut( C,{\bf x},\ph)$
of $\Aut(C,{\bf x})$ on the central fiber $C$. Then $\kappa$ and $\Pi$ can be made
$\Aut(C,{\bf x},\ph)$-equivariant as before. The unrigidified chart $S\times V
\rightarrow \calc(M;p)$ is then also invariant under the action of $\Aut(C,{\bf x},\ph)$,
provided $V$ has been chosen appropriately.

Next we adjust our linear coordinates $z_i$. For $(D_i,{\bf x}_i)$
an unstable component of $(C,{\bf x})$ as above, $\Aut(D_i,{\bf x}_i,\ph|_{D_i})$ is
a finite subgroup of $\Aut(D_i,{\bf x}_i)$. For $i\le b'$, i.e.\ if $\sharp{\bf x}_i=1$, we
now require the unique copy of $\cz^*$ in $\Aut(D_i,{\bf x}_i)$ containing this subgroup
to coincide with the second factor in the identification $\Aut(D_i,{\bf x}_i)
=\cz\rtimes\cz^*$. Equivalently, $z_i=\infty$ is required to be a fixed point of the
action of $\Aut(D_i,{\bf x}_i,\ph|_{D_i})$.

Moreover, for fixed $i$, let $\Psi_1,\ldots,\Psi_m\in\Aut(C,{\bf x},\ph)$ be representatives
for elements of
\[
  \Aut(C,{\bf x},\ph)\Big/\Big\{\Psi\in\Aut(C,{\bf x},\ph)\,\Big|\,\Psi(D_i)=D_i\Big\}\, .
\]
Then $\Aut(C,{\bf x},\ph)\cdot D_i=\{\Psi_1(D_i),\ldots,\Psi_m(D_i)\}$ and if
$D_j=\Psi_{\mu(j)}(D_i)$ we choose $z_j=(\Psi_{\mu(j)}^{-1})^*(z_i)$.

Now let $\Psi\in\Aut(C,{\bf x},\ph)$ be arbitrary and $\Psi(D_i)=D_j$. Then necessarily
$\Psi^*z_j=R_{ji}(\Psi)\cdot z_i$ for some root of unity $R_{ji}(\Psi)\in\cz^*$.
Define $R_{j'i}=0$ for $j'\neq j$. We claim that $\Psi\mapsto(R_{ij}(\Psi))_{i,j}$
is a $(b+b')$-dimensional representation of $\Aut(C,{\bf x},\ph)$. In fact, let $\Psi$,
$\Psi'\in\Aut(C,{\bf x},\ph)$ be arbitrary and fix $i\in\{1,\ldots,b+b'\}$. If $\Psi(D_i)
=D_j$ and $\Psi'(D_j)=D_k$ we compute
\begin{eqnarray*}
  R_{ki}(\Psi'\circ\Psi)\cdot z_i&=&(\Psi'\circ\Psi)^*z_k
  \ =\ \Psi^*(\Psi')^*z_k\ =\ R_{kj}(\Psi')\cdot\Psi^*z_j\\
  &=&R_{kj}(\Psi')\cdot R_{ji}(\Psi)\cdot z_i
  \ =\ \sum_{j'}R_{kj'}(\Psi')\cdot R_{j'i}(\Psi)\cdot z_i\,,
\end{eqnarray*}
or $R(\Psi'\circ\Psi)=R(\Psi')\circ R(\Psi)$. As with $A(\Psi)$ above we deduce
\[
  \Lambda(\Psi\cdot(s,v))\ =\ R(\Psi)\cdot\Lambda(s,v)\quad
  \mbox{for any }\Psi\in\Aut(C,{\bf x},\ph)\,,
\]
i.e.\ $\Lambda$ is equivariant with respect to the affine linear action $R$ of
$\Aut(C,{\bf x},\ph)$ on $\cz^{b+b'}$.

The final ingredient is an $\Aut(C,{\bf x},\ph)$-invariant choice of complementary
subspace $N$ to $\bar V=\ker D\Lambda(0,0)$. Then $\nu:S\times\bar V\times N
\rightarrow S\times V$ will be $\Aut(C,{\bf x},\ph)$-equivariant, and we have
found the (complex linear) action of $\Aut(C,{\bf x},\ph)$ on our local models.

Injectivity of the maps $S\times\bar V/\Aut(C,{\bf x},\ph)\rightarrow\calc(M;p)$
induced by the structure maps can be derived as in (4), but one has to use that
$S\times\bar V$ is a local slice to the action of $\Aut(C,{\bf x})$. So $d_\ph$ has to be
replaced by the infimum over $d(\Psi)$ with $\Psi\not\in B_\delta(\Aut(C,{\bf x},\ph))
\subset\Aut(C,{\bf x})$ for sufficiently small $\delta>0$ in any metric on $\Aut(C,{\bf x})$.
\vspace{1ex}


\subsection{Change of coordinates and morphisms}
\label{section_change_coord}
Next we will give open embeddings of our local uniformizing systems,
that is we check compatibility of our local orbifold structures.
At the same time we will show independence of choices, including the Riemannian
and almost complex structure on $M$ claimed in the theorem. The whole computation
will be very similar to the proof of Proposition~\ref{bundle_structure}, with a
non-linear term added.
\vspace{1ex}

\noindent
{\em Change of coordinates: Stable domains.}\\
According to Remark~\ref{construct_orbifold} it suffices to cover the following
situation: Let $q_{(C,{\bf x},\ph)}: S\times V\rightarrow \calc(M)$ be a local uniformizing system
centered in $(C,{\bf x},\ph)$ and constructed from the \mbox{(semi-)} universal deformation
$q:\calc\rightarrow S$ of $(C,{\bf x})$ via a trivialization $\Pi: S\times \check L_1^p
(C;\ph^*T_M)\simeq {q_1^p}_*\kappa^* \ph^*T_M$ and the exponential map $\exp^\rho$,
cf.\ Section~\ref{section_uniformizer}. For the time being we assume $(C,{\bf x})$ stable.
The trivialization $\Pi$ depended on: an admissable retraction $\kappa:
\calc \rightarrow C=q^{-1}(0)$ with lift $\tilde \kappa: \kappa^*\ph^*T_M
\rightarrow \ph^*T_M$; an open cover $\calu$ of $\calc$ together with
a realization as continuous family of circular domains or possibly degenerate cylinders;
a subordinate partition of unity $\rho_i$; local, fiberwise holomorphic trivializations of
$\kappa^*\ph^*T_M|_{U_i}$ with respect to the holomorphic structure on $\kappa^*T_M$
induced by the complex structure $J$ and the Riemannian metric $\rho$ on $T_M$. Then let
$q_{(C',{\bf x}',\ph')}: S'\times V'\rightarrow \calc(M)$ be another such local uniformizing
system centered in $(C'{\bf x}',\ph')= (C_{s'_0},{\bf x}_{s'_0},\ph(s'_0,v'_0))$
for some $(s'_0,v'_0)\in S\times V$. The corresponding data are $q':\calc'
\rightarrow S'$, $\Pi'$, $\exp^{\rho'}$, $\kappa'$, $\calu'$, $\rho_i'$,
and the holomorphic structure on ${\kappa'}^*T_M$ induced by another
almost complex structure $J'$ on $M$. We have to show that possibly after
shrinking $S'$ and $V'$ there is an open embedding $S'\times V'\hookrightarrow
S\times V$ compatible with the action of $\Aut(C',{\bf x'},\ph')$.

The first thing to be remarked is that due to the universal property of $q$ we get
an open (holomorphic) embedding of $\calc'$ into $\calc$ compatible
with $q'$ and $q$, at least after shrinking $S'$. We henceforth identify $\calc'$ and $S'$ with
open subsets of $\calc$ and $S$ without further notice. Again by universality
the automorphism group $\Aut(C',{\bf x}')$ is thereby identified with the stabilizer of
$s'_0\in S$ in $\Aut(C,{\bf x})$, making the embedding $\calc'\subset \calc$ equivariant
with respect to $\Aut(C',{\bf x}')$.

Recall that we restricted ourselves
to local uniformizing systems with the property that $d^\rho(\ph(\kappa_{s}(z)),
\ph(s,v)(z))$ does not exceed half of the injectivity radius of $\rho$ for any
$(s,v)\in S\times V$ and any $z\in C_s$. Moreover, we can then take
$V'$ so small that the $\rho$-distance of $\ph'(q(z),v')(\kappa'(z))$ to
$\ph(\kappa(z))$ is less than the injectivity radius of $\rho$ along $\ph$,
for any $z\in \calc'$ and any $v'\in V'$. We now use the fact that for
$P,P'\in M$ with $d^\rho(P,P')$ less than the injectivity radius of $\rho$ at $P$,
the map
\[
  \Theta_{PP'}:B_r(0)\subset T_{M,P'}\longrightarrow T_{M,P}\,,
  \quad v\longmapsto(\exp^\rho_P)^{-1}(\exp^{\rho'}_{P'}v)
\]
is well-defined for $r$ sufficiently small, and is a diffeomorphism onto the image depending
smoothly on variations of $P,P'$ with respect to local trivializations of $T_M$. If we write
$\Theta(z):= \Theta_{\ph(\kappa(z))\ph'(\kappa'(z))}$ then $\Theta$ is a continuous map
from an open neighbourhood of the zero section of ${\kappa'}^*{\ph'}^* T_M$ to
$\kappa^* \ph^*T_M$. 
The change of coordinates $(s,v')\mapsto\hat \sigma(s,v')=(s,v)$
with $\ph'(s,v')=\ph(s,v)$ thus takes the form
\[
  \hat\sigma(s,v')\ =\ \Big(s,\Pi_{s}^{-1}\Theta ( \Pi'_{s} v')\Big)\, .
\]
By equivariance of $\kappa$ and $\kappa'$ it holds $\Theta_z\circ (\Psi_s^{-1})^*
= (\Psi_s^{-1})^*\circ \Theta_{\Psi^{-1}(z)}$ for any $z\in C_{\Psi(s)}$. This together
with equivariance of $\Pi$ and $\Pi'$ (cf.\ Section~\ref{section_uniformizer}) shows
compatibility of $\hat\sigma$ with the action of $\Aut(C',{\bf x}',\ph')$.

Let $\dbar$, $\dbar'$ be the $\dbar$-operators on $\kappa^* \ph^*T_M$
and ${\kappa'}^* {\ph'}^*T_M$ induced by $J,\rho$ and $J',\rho'$ respectively.
In the notation of the proof of Proposition~\ref{bundle_structure} it remains
to check continuity in the compact-open topology of the map associating to
$(h_i,\gamma)\in \check C^0(\calu'; \calo({\ph'}^*T_M,{\ph'}^*J')) \oplus
\check L^p(C'; {\ph'}^* T_M \otimes\bar\Omega)$ the expressions
\begin{eqnarray*}
  &({\tilde\kappa'_s}{}^\star)^{-1}\dbar'\Theta\sum_j\rho_j(s)
  \Big(h_j(s)+T^j_s\gamma(s)\Big)\,, \\
  &({\tilde\kappa'_s}{}^*)^{-1}{H'_s}^i\Theta
  {\textstyle \sum_j}\rho_j(s)(h_j(s)+T^j_s\gamma(s))\Big)\,,
\end{eqnarray*}
as map to $\check L_1^p(C; {\ph}^*T_M\otimes \bar\Omega) \oplus
\check C^0(\calu; \calo({\ph}^*T_M,{\ph}^*J))$. Here we have already ignored
the obviously continuous averaging over automorphism groups involved in our
definition of $\Pi$ and $\Pi'$ (Section~\ref{section_uniformizer}). 
The only facts we will need about $\Theta$ is that its
derivatives in the vertical directions (i.e.\ along the fibers of ${\kappa'}^*{\ph'}^* T_M$)
are continuous, while derivatives along horizontal directions (with respect to any local
trivializations) are continuous sections of appropriate $\check L^p$-bundles. This
together with the arguments already presented in the proof of
Proposition~\ref{bundle_structure} shows the claimed continuity.

Note that by changing the roles of $S'\times V'$ and $S\times V$, i.e.\ by shrinking
$S\times V$ to a small neighbourhood of some $(s,v)\in S'\times V'$, it follows
that $\hat\sigma$ is indeed an open embedding.
\vspace{1ex}

\noindent
{\em Change of coordinates: The general case.}\\
Now assume $(C,{\bf x})$ unstable and let $S\times \bar V\subset S\times V$ be a
rigidifying slice. Let $q: S\times V \rightarrow \calc (M;p)$ induce the structure map.
As shown in Section~\ref{sect_rigidification}, $S\times \bar V$ is the
topological quotient of $S\times V$ by the equivalence relation $R$ generated by
the germ of action of $\Aut^0 (C,{\bf x})$. Given $(C',{\bf x}',\ph')\in \img q$ and sufficiently
close to $(C,{\bf x},\ph)$ as above (depending on the injectivity radius of $\rho$ etc.),
we choose the local uniformizing system with center $(C',{\bf x}',\ph')$ to be a slice
$S'\times \bar V'$ in $S'\times V'$ with $S'$, $V'$ sufficiently small as before. Again,
$S'\times \bar V'$ is the quotient of $S'\times V'$ by the equivalence relation $R'$,
generated by the germ of action of $\Aut^0 (C',{\bf x}',\ph')$. 

Now $(C',{\bf x}')$ belongs to some $s_0\in S$, and with a choice of $s_0$ the
unstable components of $(C',{\bf x}')$ can be identified with
a subset $D_1,\ldots,D_a$ of the set of unstable components of $(C,{\bf x})$ via $\kappa$.
Let $G:=\{\Psi\in\Aut^0(C,{\bf x})\mid\Psi|_{D_i}=\id,\ i=1,\ldots,a\}$. By our
explicit description of the semiuniversal deformation it is not hard to see
that the local action of $G$ fibers $S$ smoothly near $s_0$, and that the restriction of
$q:\calc\rightarrow S$ to a smooth analytic slice of the action of $G$ at $s_0$
is a semiuniversal deformation of $(C',{\bf x}')$, hence locally isomorphic to $S'$.
Again we can thus identify $\calc'$ with a locally analytic subset of $\calc$, this time
of codimension equal to the dimension of $G$. The map $\hat \sigma: S'\times \bar V'
\hookrightarrow S\times \bar V$ is then defined as composition of $(s,v')\mapsto
(s,\Pi_s^{-1} \Theta\Pi'_s v')$ with the quotient map
$S\times V\rightarrow S\times \bar V$. Equivariance and continuity of $\hat\sigma$
are thus inherited by the corresponding properties of the unrigidified map $S' \times V'
\rightarrow S\times V$. Again we can change the roles of $S'\times \bar V'$ and
$S\times \bar V$ to conclude that $\hat\sigma$ is an open embedding.


\subsection{Comparison of topologies: $\calc(M;p)$ is Hausdorff}
To finish the proof of Theorem~\ref{orbifold} it remains to show that the orbifold topology
(the topology induced by the structure maps $S\times V\rightarrow \calc(M;p)$) is finer than
the $C^0$-topology on the space $\calc (M)$ of continuous stable complex curves
in $M$, cf.\ Definition~\ref{C0top}. Since the latter has already been shown to be
Hausdorff (Proposition~\ref{C^0_hausdorff}) this will also establish the necessary
Hausdorff property for the orbifold topology. If $S\times V$ is a chart let us call
{\em $L^\infty$-topology} on $S\times V$ the topology induced by the embedding
$S\times V\subset S\times L^\infty(C;\ph^*{T_M})$.
\begin{prop}\label{covtop}
  The orbifold topology on $\calc(M;p)$ is finer than the $C^0$-topology. In fact,
  the $C^0$-topology is generated by the $L^\infty$-topology on unrigidified charts
  $S\times V\rightarrow \calc(M;p)$.
\end{prop}
\pf
Using an admissable retraction $\kappa:\calc \rightarrow C$ for comparison of
$(C_s,{\bf x}_s, \ph(s,v))$ with $(C,{\bf x},\ph)$ in the $C^0$-topology one sees
that the $C^0$-topology is coarser than the topology coming from $S\times
L^\infty( C;\ph^*T_M)$ on the local uniformizing systems.

Conversely, let $(C_\nu,{\bf x}_\nu,\ph_\nu)$ be a sequence of stable curves in $M$
of Sobolev class $L_1^p$, converging towards $(C,{\bf x},\ph)$ in the $C^0$-topology.
Let $\kappa_\nu: C_\nu \rightarrow C$ be a sequence of retractions exhibiting
this convergence. We first reduce to stable domains $(C,{\bf x})$. To this
end we choose a tuple ${\bf y}$ of points on $C$ such that
$(C,{\bf x}\vee {\bf y})$ is stable. Moreover, we require that for any
unstable component $D$ of $(C,{\bf x})$ the distance of the image under $\ph$
of any two of the special points, involving at least one newly added point $y_j$, is
larger than some constant $d>0$. This is possible by stability. Then for $\nu$ sufficiently
large $(C_\nu,{\bf x}_\nu\vee \kappa_\nu^{-1}({\bf y}))$ is also stable. It thus suffices
to treat the case of stable domains.

If $(C,{\bf x})$ is stable we can apply Lemma~\ref{analytic-dilation_comp} to
deduce, for sufficiently large $\nu$, the existence of biholomorphisms
$\Xi_\nu: C_\nu \rightarrow C_{s_\nu}$ with
$s_\nu\rightarrow 0\in S$, $\Xi_\nu {\bf x}_\nu = {\bf x}_{s_\nu}$ and
$d(\kappa_s\circ\Xi_\nu,\kappa_\nu) \rightarrow 0$ in the limit $\nu\rightarrow\infty$.
Thus
\[
  d(\ph\circ\kappa_{s_\nu},\ph_\nu\circ\Xi_\nu^{-1})\ =\
  d(\ph\circ\kappa_{s_\nu} \circ\Xi_\nu,\ph_\nu)
  \ \le\ d(\ph\circ\kappa_{s_\nu} \circ\Xi_\nu,\ph\circ\kappa_\nu)
  +d(\ph\circ\kappa_\nu,\ph_\nu)\,.
\]
Of these the first term can be made arbitrarily small by equicontinuity of $\ph$
and the fact that $\kappa_{s_\nu}\circ\Xi_\nu$ tends to $\kappa_\nu$ in $C^0$.
The second term tends to zero with $\nu$ because the $\kappa_\nu$ exhibit
$C^0$-convergence. But any of our unrigidified charts for $\calc(M;p)$ is of the
form $S\times V$ with $V\subset \check L_1^p(C;\ph^*T_M)$ in fact open in the
$L^\infty$-topology. This shows that for sufficiently large $\nu$, the curves
$(C_\nu,{\bf x_\nu}, \ph_\nu)$ are of the form $\ph(s_\nu,\ph(s_\nu,v_\nu))$ for some
$v_\nu\in V$.
\qed


\subsection{The universal curve and morphisms}
On a set-theoretic level the universal curve $\Gamma(M;p)$ and the
morphisms $\Phi,\pi, p$ are well-defined. As for
uniformization near the pair $[(C,{\bf x},\ph), P]\in \Gamma(M;p)$ let
$S\times\bar V$ be a chart for $\calc(M;p)$ centered in $(C,{\bf x},\ph)$
and modelled on the semi-universal deformation $(q:\calc\rightarrow S,
\underline{\bf x})$ of $(C,{\bf x})$. Then we may just take $\calc\times
\bar V$ as local uniformizing systems for $\Gamma(M;p)$. These are obviously
compatible with the set-theoretic structure of $\Gamma(M;p)$. The
orbifold topology on $\Gamma (M;p)$ defined by such local uniformizing systems
is Hausdorff in view of the Hausdorff property of the orbifold topology
on $\calc(M;p)$. The map $\pi: \Gamma(M;p)\rightarrow \calc(M;p)$ is readily
uniformized by
\[
  q\times \id_{\bar V}:\ \calc\times\bar V\ \longrightarrow S\times \bar V\,.
\]
Similarly, the tuple of sections of $\pi$ is locally pull-back of $\underline{\bf x}$.
The universal morphism $\Phi$ factorizes over the map to the topological space
$\Gamma(M;p)$ and is easily seen to be continuous. For the morphism to
$\bigcup \calm_{g,k}$ note that the base space $S$ of the semi-universal
deformation of $(C,{\bf x})$ fibers over the base space $\bar S$ of the
universal deformation of $(C,{\bf x})^\st$ (Remark~\ref{prestable_deform}).
The map $\pi$ is locally uniformized by the projection
\[
  S\times \bar V\ \longrightarrow\ S\ \longrightarrow\ \bar S\,,
\]
and these morphisms are clearly compatible with open embeddings
of local uniformizing systems.


\section{The bundle and section}
To apply the theory of localized Euler classes developed in Chapter~1 we still have to
construct a weak Banach orbibundle $E$ over $\calc(M;p)$ and a Kuranishi
section of $E$ with zero locus the set of pseudo holomorphic
curves $\calc^\hol(M,J)\subset\calc(M;p)$. We take ${\hat E}_{(C,{\bf x},\ph)}
=\check L^p(C;\ph^*T_M\otimes_\cz{\bar\Omega}_C)$ and ${\hat s}_{\dbar,J}
(C,{\bf x},\ph)=\dbar_J\ph:=\frac{1}{2}(D\ph+J(\ph)\circ D\ph\circ j)$, the latter
being easily checked to be $(J,j)$-antilinear and so to lie in ${\hat E}_{(C,
{\bf x},\ph)}$. Note that $\hat E$ also depends on the choice of almost complex
structure $J$ on $M$ by the tensor product involved in its definition. Clearly,
${\hat s}_{\dbar,J}(C,{\bf x},\ph)=0$ iff $(C,{\bf x},\ph)$ is a pseudo holomorphic
curve on $(M,J)$.

We are going to verify in this chapter: 1) there is such a Banach orbibundle $E$
over $X$\ \
2) ${\hat s}_{\dbar,J}$ defines a continuous section $s_{\dbar,J}$ of $E$\ \
3) ${\hat s}_{\dbar,J}$ is continuously differentiable for fixed domains with
derivatives that are continuous at the centers of our local uniformizing systems\ \
4) this derivative is Fredholm of locally constant index and naturally oriented if
$(C,{\bf x},\ph)$ is pseudo holomorphic with respect to $J$. Together with the last
point this defines the Fredholm structure for $s_{\dbar,J}$\ \ 
5) there exists an orbibundle $F$ over $\calc(M;p)$ that is locally of finite rank,
together with a morphism $\tau:F \rightarrow E$ spanning the cokernel of the differentials
of $s_{\dbar,J}$ for fixed domains; that is, we define a Kuranishi structure for
the section.

We will then be in position to apply the main theorem from the first chapter
(Theorem~\ref{virfc}) to produce a class $[E,s_{\dbar,J}]\in H_*(\calc^\hol(M,J))$,
the {\em virtual fundamental class} of $\calc^\hol(M,J)$, which will be further investigated
in Chapter~\ref{chapter_GW}. We now begin with the verification of points 1--5.


\subsection{The bundle $E$}\label{bundleE}
The construction of $E$ is rather straightforward, but there is one point of caution:
With our choice of connection, parallel transport does not in general commute with $J$,
so local trivializations of $E$ obtained by parallel transport would not be well-defined.
One could remedy this by working with almost complex connections ($\nabla J=0$),
but then changes of $J$ are not discussed so easily. Instead we chose to stay with
torsion free connections and to force parallel transport to respect the complex structure
by setting
\[
  \Gamma_{P,Q}^\cz:=\frac{1}{2}\Big(\Gamma_{P,Q}-J(P)\circ\Gamma_{P,Q}\circ J(Q)\Big)\,,
\]
where $\Gamma_{P,Q}:T_QM\rightarrow T_PM$ for $P$, $Q$ sufficiently close is
parallel transport along the unique geodesic joining $Q$ with $P$.
Note that $\Gamma^\cz_{P,Q}$ is still an isomorphism for $P$ sufficiently
close to $Q$ for $\Gamma_{P,P}=\id$.

Now let be given one of the local uniformizing systems
with center $(C,{\bf x},\ph)$ from the proof of Theorem~\ref{orbifold}, let $z_0=z_0(s)$ be
a holomorphic coordinate on $U_0(s)\subset C_s$ varying continuously with $s$,
and $z_i=z_i(s)$ for $i>0$ be $z$ resp.\ $w$ on the two connected components
of $Z_{t_i(s)}\setminus S_{t_i(s)}$ (so this is a holomorphic coordinate with
discontinuities along $S_{t_i(s)}$). Let furthermore $\{\rho_i\}$ be a partition
of unity subordinate to $\{U_i(0)\}$  and $\kappa:\calc\rightarrow C$ an
admissable retraction of continuous families of prestable curves. A local trivialization
of $E$ is now given by
\begin{eqnarray*}
  \Pi^E_{s,v}:\check L^p(C;\ph^*T_M\otimes{\bar\Omega}_C)
  &\longrightarrow&\check L^p(C_s;\ph(s,v)^*T_M\otimes{\bar\Omega}_{C_s})\\
  \alpha=\sum\alpha_i d{\bar z}_i&\longmapsto&\sum\Big(\Gamma^\cz_{\ph(s,v),\ph(s,0)}
  (\kappa_s^*\alpha_i)\Big)d{\bar z}_i\,,
\end{eqnarray*}
where $\alpha_i d{\bar z}_i=\rho_i\cdot\alpha$.
Changes of the local uniformizing system, say coming from $\ph':C'\rightarrow M$,
$\kappa'_s$, $\{\rho'_i\}$, $v'=v'(s,v)$, $\Pi'_s$ with holomorphic coordinates
$z'_i=z'_i(s)$, lead to
\[
  \Big({\Pi^E_{s,v'}\!\!\!}'\Big)^{-1}\circ\Pi^E_{s,v}\Big(\sum\alpha_i d{\bar z}_i\Big)
  =\sum_{i,j}\rho'_j\cdot\Big(({\kappa'_s}^*)^{-1}\Gamma^\cz_{\ph'(s,0),\ph(s,v)}
  \Gamma^\cz_{\ph(s,v),\ph(s,0)}(\kappa^*_s\alpha_i)\Big)\, d{\bar z}'_i\,,
\]
which is evidently continuous in the compact-open topology as family of maps from
$\check L^p(C;\ph^*T_M\otimes{\bar\Omega}_C)$ to $\check L^p(C';
{\ph'}^*T_M\otimes{\bar\Omega}_{C'})$.

If $\Aut(C,{\bf x},\ph)$ is not trivial we make $\Pi_{s,v}^E$ equivariant by averaging
as in Section~\ref{section_uniformizer}.
The case of unstable domains is treated by restriction of the
above construction to $S\times\bar V\hookrightarrow S\times
\check L_1^p(C;\ph^*T_M)$.


\subsection{The section $s_{\dbar,J}$}
Using the local trivialization $\Pi^E$ of $E$ and our local uniformizing systems
$S\times\bar V$ we may now define $s_{\dbar,J}$ locally
by
\[
  S\times\check L_1^p(C;\ph^*T_M)\supset S\times \bar V\longrightarrow
  \check L^p(C;\ph^*T_M\otimes{\bar\Omega}_C)\,,\quad
  (s,v)\longmapsto\Big(\Pi^E_{s,v}\Big)^{-1} \dbar_J \ph(s,v)\,.
\]
Obviously, this map is equivariant and compatible with open embeddings of local
uniformizing systems. It remains to show continuity.
To avoid arguments with Jacobi fields
we take up a strictly local point of view. It clearly suffices to check for continuity
after pairing with finitely many generating sections of $T_C$.
Choosing the support of these vector fields sufficiently small we
are reduced to the case of $C_s=G_s\subset\cz$,
where we may take $\partial/\partial_x$ as vector field ($z=x+iy$, so $j_s(\partial/
\partial_x)=\partial/\partial_y$ for any $s$), and $M\subset\rz^{2n}$. The
identification of $M$ with a subset in $\rz^{2n}$ gives a trivialization of $\ph(s,v)^*T_M$
(as real vector bundle) for $(s,v)$ close to $0$. By taking the complex linear part
we get from this another local trivialization
of $E$ compatible with $\Pi^E_{s,v}$ through a continuous family of toplinear
isomorphism in the compact-open topology. Since projection onto a complex linear
part is clearly continuous, it suffices to consider (as map $S\times \bar V\rightarrow
\check L^p(G_0;\rz^{2n})$)
\[
  (s,v)\longmapsto\Big[\partial_x(\exp_{\ph\circ\kappa_s(\, .\, )}(\Pi_s v)(\, .\, )
  +J(\ph(s,v)(\, .\, ))\partial_y(\exp_{\ph\circ\kappa_s(\, .\, )}(\Pi_s v)(\, .\, )\Big]
  \circ\kappa_s^{-1}\,,
\]
where $\partial_x$, $\partial_y$ act on the ``$.$'' that stands for the $z$-variable.
Now as map to $L^\infty(G_0;$ $\End(\rz^{2n}))$, $(s,v)\mapsto
J(\ph(s,v)\kappa_s^{-1}(\, .\, ))$
is continuous, so it is enough to prove contin\-uity
for each of the terms involving the derivatives. We carry this out for the first one,
the other one being completely analogous. Writing $\exp(P,v):=\exp_P v$ (now viewed
as map from an open set in $\rz^{2n}\times\rz^{2n}$ to $\rz^{2n}$) this
becomes at $(s,v)$
\[
  D_1\exp\Big(\ph,(\Pi_s v)\circ\kappa_s^{-1}\Big)
  \Big(\partial_x(\ph\kappa_s)\circ\kappa_s^{-1}\Big)
  +D_2\exp\Big(\ph,(\Pi_s v)\circ\kappa_s^{-1}\Big)
  \Big((\partial_x\Pi_s v)\circ\kappa_s^{-1}\Big)\, .
\]
As maps from $\hat U$ to $L^\infty(G_0;\End(\rz^{2n}))$ the derivatives of the
exponential map in either entry are continuous, while $(s,v)\mapsto\partial_x\Pi_s v\circ
\kappa_s^{-1}$ is continuous by Proposition~\ref{localtriv},3.
Using the Lipschitz property of $\kappa_s$ it is also easy to check the continuity of
$s\mapsto\partial_x(\ph\circ\kappa_s)\circ\kappa_s^{-1}\in L^\infty(G_0;\rz^{2n})$.
This proves continuity of $s_{\dbar,J}$.


\subsection{Fredholm structure}
We are going to show now that our charts $S\times \bar V$ for $\calc(M;p)$
together with its naturally associated trivialization of $E$ (i.e.\  using the
same retraction $\kappa$) endow
$s_{\dbar,J}$ with a Fredholm structure, with $(S=W,L,V)$ in
Definition~\ref{Fredholm_structure} equal to $({\rm pt},S,\bar V)$.
\vspace{2ex}

\noindent
{\em Differentials of $s_{\dbar,J}$ for fixed domains.}\\
For fixed $s$ we want to compute the derivative of
\[
  \check L_1^p(C;\ph^*T_M)\supset \bar V\longrightarrow
  \check L^p(C;\ph^*T_M\otimes{\bar\Omega}_{C_s})\,,\quad
  v\longmapsto\Big(\Pi^E_{s,v}\Big)^{-1}{\hat s}_{\dbar,J}(C_s,{\bf x}_s,\ph(s,v))
\]
at any $v\in \bar V\subset \check L_1^p(C;\ph^*T_M))$. Writing
$\ph_s=\ph\circ\kappa_s$, $v_s:=\Pi_s(v)\in\check L_1^p(C_s;\ph_s^*T_M)$, then for
$w\in\check L_1^p(C;\ph^*T_M)$ we put $w_s:=\Pi_s w$, $\psi(z,t):=\ph(s,v+tw)(z)
=\exp_{\ph_s(z)}(v_s+tw_s)(z)$. We have to look at
\[\hspace{1.5cm}
  \left.\frac{d}{dt}\right|_{t=0}\Gamma^\cz_{\ph_s(z),\psi(z,t)}\Big(D_\xi\psi(z,t)
  +J(\psi(z,t))\cdot D_{j\xi}\psi(z,t)\Big)\,,\hspace{1.5cm}(*)
\]
where $\xi$ is a vector field on $C$. Note that
$D_\xi\psi(.,t)$ does not make sense pointwise, but as $\check L^p$ vector field along
$z\mapsto\psi(z,t)$. Now observe that as $\End(\rz^{2n})$ valued functions on
$C_s\times[0,\delta]$ ($\delta\ll1$), $\Gamma_{\ph_s(z),\psi_t(z)}$ differs from
$\Gamma_{\ph_s(z),\psi(z)}\circ\Gamma_{\psi(z),\psi_t(z)}$ only by a term
coming from holonomy. Here we have put $\psi(z):=\psi(z,0)=\ph(s,v)(z)$
and $\psi_t(z):=\psi(z,t)$. In the limit $t\rightarrow 0$ the holonomy term
leads to a zero order operator vanishing at $(s,v)=(0,0)$
(an integral along $\lambda\mapsto\ph(s,\lambda\cdot v)$
involving the curvature). So this term is uniformly continuous for fixed $s$ or at
$(s,v)=(0,0)$. On the other hand,
\[
  \left.\frac{d}{dt}\right|_{t=0}\Gamma_{\psi,\psi_t}D_\xi\psi_t
  \ =\ \nabla^\psi_{\frac{d}{dt}}D_\xi\psi_t\ =\ \nabla^\psi_\xi\dot{\psi}_0\,,
\]
since $\nabla$ was torsion free and
\begin{eqnarray*}
  \left.\frac{d}{dt}\right|_{t=0}\Big(J(\psi)\Gamma_{\psi,\psi_t}
  J(\psi_t)D_\xi\psi_t\Big)&=&J(\psi)\nabla_{\frac{\partial}{\partial t}}^\psi
  \Big(J(\psi_t)D_\xi\psi_t\Big)\\
  &=&-\nabla_\xi^\psi\dot{\psi}_0+J(\psi)(\nabla_{\dot{\psi}_0}J)D_\xi\psi\, .
\end{eqnarray*}
For the term coming from $J(\psi_t)\partial_y\psi_t$ we obtain
\[
  2J(\psi)\nabla_{\frac{\partial}{\partial y}}^\psi\dot{\psi}_0
  +(\nabla_{\dot{\psi}_0}J)D_{j\xi}\psi\, .
\]
So with $\dot{\psi}_0=\Gamma_{\psi,\ph_s}w_s=:w_{s,v}$ and up to the
zero order operator, expression $(*)$ equals
\[
  \Gamma_{\ph_s,\psi}\Big[\nabla_\xi^\psi w_{s,v}
  +J(\psi)\nabla_{j\xi}^\psi w_{s,v}
  +(\nabla_{w_{s,v}}J)(\psi)\circ\partial_J\psi(j\xi)\Big]\,.
\]
Intrinsically this reads
\begin{eqnarray*}
  \lefteqn{\Big(\Pi^E_{s,v}\Big)^{-1}\Big[\nabla^{\ph(s,v)}(\Pi_{s,v} w)
  +J(\ph(s,v))\circ\nabla^{\ph(s,v)}(\Pi_{s,v} w)\circ j_s}\hspace{5cm}\\
  &&+(\nabla_{\Pi_{s,v} w}J)(\phi(s,v))\circ\partial_J\phi(s,v)\circ j_s\Big]\,,
\end{eqnarray*}
where we wrote $\Pi_{s,v}$ for $\Pi_s(v)$ to emphasize the interpretation of
this map as local trivialization of the Banach orbibundle with fiber $\check L_1^p
(C;\ph^*T_M)$ at $(C,{\bf x},\ph)$. This should be thought of as relative tangent
bundle $T_{\calc(M;p)/\calm_{g,k}}$. Then the above part of the relative differential
can be interpreted as morphism $\sigma_J:T_{\calc(M;p)/\calm_{g,k}}\rightarrow E$
with
\[
  {\hat\sigma}_J(C,{\bf x},\ph):w\longmapsto\frac{1}{2}\Big[\nabla^\ph w
  +(J\circ\ph)\circ\nabla^\ph w\circ j+(\nabla_w J\circ\ph)\circ\partial_J\ph\circ j\Big]\, .
\]
This is norm continuous as map from $\check L_1^p(C;\ph^*T_M)$ to $\check L^p
(C;\ph^*T_M\otimes{\bar\Omega}_C)$ ($s$ fixed) and continuous in the compact-open
topology under changes of $s$.

To show norm continuity at $(s,v)=(0,0)$ we need to take
a closer look at $\sigma_J$ (at fixed $(C,{\bf x},\ph)$ for the time being): First we notice
that ${\hat\sigma}_J=\dbar_{\ph^*T_M,J}+R$ with $R:\check L_1^p(C;\ph^*T_M)
\rightarrow\check L^p(C;\ph^*T_M\otimes{\bar\Omega}_C)$ complex antilinear,
i.e.\ anticommuting with $J$. In fact, by definition of the $\dbar$ operator on $\ph^*T_M$
this is nothing but the decomposition of ${\hat\sigma}_J$ into $J$-linear and antilinear
parts. Let us evaluate $R$ at some $w\in\check L_1^p(C;\ph^*T_M)$ and $\xi$
a vector field on $C$, using $\partial_J\ph(j\xi)=JD_\xi\ph+\dbar_J\ph(j\xi)$ and setting
$u=D_\xi\ph$:
\begin{eqnarray*}
  \lefteqn{4R(w,\xi)=2{\hat\sigma}_J(w)(\xi)+2J{\hat\sigma}_J(Jw)(\xi)}\hspace{0.5cm}\\
  &=&\nabla_\xi w+J\nabla_{j\xi}w+(\nabla_w J)\partial_J\ph(j\xi)\ph
  +J\nabla_\xi(Jw)+J(\nabla_{Jw}J)\partial_J\ph(j\xi)\ph\\
  &=&(\nabla_w J)\partial_J\ph(j\xi)\ph+J(\nabla_\xi J)w-(\nabla_{j\xi}J)w
  -\nabla_{j\xi}(Jw)+J(\nabla_{Jw}J)\partial_J\ph(j\xi)\ph\\
  &=&(\nabla_w J)(Ju)+(\nabla_w J)(\dbar_J\ph(j\xi))+J(\nabla_u J)w-(\nabla_{Ju}J)w
  -2(\nabla_{\dbar_J\ph(j\xi)}J)w\\
  &&+J(\nabla_{Jw}J)(Ju)+J(\nabla_{Jw}J)(\dbar_J\ph(j\xi))\, .
\end{eqnarray*}
The first, third, fourth and sixth term combine to $4N_J(w,u)$, the almost complex torsion
of $J$ (Nijenhuis tensor). We obtain
\[
  4N_J(w,u)+\Big[(\nabla_w J)+J(\nabla_{Jw}J)\Big]\dbar_J\ph(j\xi)
  -2(\nabla_{\dbar_J\ph(j\xi)}J)w\, .
\]
We are now in a position to control the operator norm of
\[
  \Big(\Pi_{s,v}^E\Big)^{-1}\circ{\hat\sigma}_J(s,v)\circ\Pi_{s,v}-{\hat\sigma}_J(0,0):
  \check L_1^p(C;\ph^*T_M)\longrightarrow\check L^p(C;\ph^*T_M\otimes{\bar\Omega}_C)
\]
in terms of $(s,v)$, which will then prove continuity of the derivative at $(s,v)=(0,0)$:
Since we may bound $(\Pi_{s,v}^E)^{-1}\circ{\hat\sigma}_J(s,v)\circ\Pi_{s,v}
-(\Pi_{s,0}^E)^{-1}\circ{\hat\sigma}_J(s,0)\circ\Pi_{s,0}$ by $|\!|v|\!|_{\check L_1^p}$
independently of $s$ (cf.\ the formula that established continuity in the compact-open
topology above) it suffices to treat the case $v=0$, so $\Pi_{s,v}^E=\kappa_s^\star$
and $\Pi_{s,v}=\Pi_s$. The $J$-linear part of ${\hat\sigma}_J(s,0)$ is exactly the
$\dbar$ operator on $\kappa_s^*\ph^*T_M$, which is uniformly continuous in the
chosen trivialization for any $s$.

It remains to treat the various terms coming from the antilinear part $R$
(after pairing with a continuously varying vector field $\xi=\xi(s)$).
The first of these maps (at $(s,0)$) $w\in\check L_1^p(C;\ph^*T_M)$ to $(N_J\circ\ph_s)
(\Pi_s w,D_{\xi(s)}\ph_s)$ and
\begin{eqnarray*}
  &(N_J\circ\ph)\Big((\kappa_s^*)^{-1}\Pi_s w,(\kappa_s^*)^{-1}
  D_{\xi(s)}\ph_s\Big)-(N_J\circ\ph)(w,D_{\xi(0)}\ph)\\
  &=\ (N_J\circ\ph)\Big((\kappa_s^*)^{-1}\Pi_s w-w,(\kappa_s^*)^{-1}
  D_{\xi(s)}\ph_s\Big)\\
  &\ \ +(N_J\circ\ph)(w,(\kappa_s^*)^{-1}D_{\xi(s)}\ph_s-D_{\xi(0)}\ph)\, .
\end{eqnarray*}
So this gets arbitrarily small in $\check L^p$ norm with $s$ by uniform boundedness
of $N_J\circ\ph$, $L^\infty$ convergence $(\kappa_s^*)^{-1}\circ\Pi_s\rightarrow\id$
and $\check L^p$ convergence $(\kappa_s^*)^{-1}D_{\xi(s)}\ph
\rightarrow D_{\xi(0)}\ph$. The next term involves $\dbar_J\ph_s$, which
tends to zero with $s$ and $v$ by continuity of $s_{\dbar,J}$,
composed with an endomorphism that depends continuously in $L^\infty$ on $s$.
Similarly for the last term $(\nabla_{\dbar_J \ph(j\xi(s))}J) \cdot\Pi_s w$.
\vspace{2ex}

\noindent
{\em Fredholm property.}\\
If $(C,{\bf x},\ph)$ is pseudo holomorphic with respect to $J$ the formula for
${\hat\sigma}_J (C,{\bf x},\ph)$ reduces to
\[
  {\hat\sigma}_J(C,{\bf x},\ph)(w)\ =\ \dbar_{\ph^*T_M,J}w+N_J(w,D\ph)\, .
\]
So up to a term of order zero this is just the $\dbar$ operator on the holomorphic
vector bundle $\ph^*T_M$, which is Fredholm by Corollary~\ref{ker_coker}:
$\ker(\dbar_{\ph^*T_M})=H^0(C;\ph^*T_M)$ and $\coker(\dbar_{\ph^*T_M})
\simeq H^1(C;\ph^*T_M)$. If $(C,{\bf x})$ is unstable we have to restrict to the tangent
space of the slice $\Lambda^{-1}(\lambda_0)\subset S\times V$, which is
the subspace of finite codimension $\bar V = \ker D\Lambda(0,0) \subset
\check L_1^p (C;\ph^*T_M)$. So the Fredholm property is also verified in this case.
This finishes the proof that our charts indeed provide a Fredholm structure for
$s_{\dbar,J}$.


\subsection{The finite rank bundle $F$}
\label{section_F}
Since the automorphism groups may act non-trivially on the cokernel of
the linearization, the domain of the morphism $\tau: F\rightarrow E$ spanning
the cokernel of the differentials of $s_{\dbar,J}$ for fixed domains must be
an orbibundle with non-trivial actions of the local groups. The existence
of such finite rank orbibundles over general Banach orbifolds is a fairly
non-trivial issue, even if one assumes the Banach orbifold to be covered by
only finitely many local uniformizers. A general solution to this problem, which
the author gave in a previous version of this paper by constructing a big finite
group into which all local groups embed in a certain compatible
way, was essentially flawed. After several more attempts the author's impression
is that the existence of such finite rank bundles may not always be true.
Fortunately, in the case of $\calc(M;p)$, the bundle $F$ can be constructed
by mimicking the standard construction of such bundles in algebraic gemetry,
cf.\ e.g.\ \cite[Prop.5]{behrend}.

To this end we now assume $J$ tamed by some symplectic form $\omega$.
By slightly deforming $\omega$ and taking a large multiple, we may assume
$\omega$ to represent an integral de Rham class. Then there exists a
$U(1)$-bundle $L$ over $M$ with $[\omega]=c_1(L)$. $L$ is the substitute for
an ample line bundle in the algebraic setting. Let $\nabla$ be a
$U(1)$-connection on $L$. Let $\pi:\Gamma\rightarrow \calc(M;p)$ be the
universal curve and $\Phi: \Gamma \rightarrow M$ be the evaluation map
sending $p\in C$ over $(C, {\bf x},\ph)\in \calc(M;p)$ to $\ph(p)$. According
to Theorem~\ref{orbifold}, $\pi$ is a morphism of orbifolds. As for the
pull-back of the tangent bundle in
Section~\ref {hol_structure} we see that via $\nabla$, $\Phi^*L$ has naturally
the structure of a continuously varying family of holomorphic line bundles
over $\pi$. And since $[\omega]$ evaluates positively on any non-constant
$J$-holomorphic curve, $\ph^*L$ is ample on any unstable component.
To achieve ampleness on the other components we tensor with the
relative dualizing bundle $\omega_\pi= \omega_{\Gamma/\calc(M;p)}$
(Lemma~\ref{dualizing_bundle}), twisted by the bundle $\calo(\underline{x_1}+\ldots+
\underline{x_k})$ of meromorphic 1-forms on the fibers of $q$ with at most simple poles
at the marked points. Note that the latter is also locally uniformized by the product of
$\bar V$ with a bundle on the semi-universal deformation $q:\calc\rightarrow S$.
Now $\Phi^*L\otimes \omega_\pi (\underline{\bf x})$ is a continuous family
of holomorphic line bundles over $\pi$ that on each irreducible component
of any fiber $C_s$ of $q$ has positive degree. So restricted to any fiber this
line bundle is ample.

Now let $(C_0,{\bf x}_0,\ph_0)$ be a stable curve in $M$. Set $G=
\Aut(C_0,{\bf x}_0,\ph_0)$. Recall the standard representation of $G$ on its
group ring $\rz[G]\simeq\rz^G$ via $g\cdot {\bf 1}_h:={\bf 1}_{gh}$.
\begin{lemma}\label{bundle_F}
  Let $K\subset \calc(M;p)$ be a compactum containing $(C_0,{\bf x}_0,\ph_0)$
  and assume that for any unstable component $D$ of any
  $(C,{\bf x},\ph)\in K$ it holds $\omega(\ph_*[D])>0$. Then over a
  neighbourhood of $K$ there exists a finite rank orbibundle $F$ with
  fiber at $(C,{\bf x},\ph)$ uniformized by a $G$-space containing a copy of
  the standard representation on the group ring $\rz[G]$.
\end{lemma}
\pf
Let $N$ be a sufficiently big natural number such that for any $|G|$ points on
$C_0$ there exist a section of $(\ph_0^*L\otimes \omega_{C_0}({\bf
x}_0))^{\otimes N}$ vanishing at all but one point. As a set we put
\[
  F\ :=\ \pi_*\Big((\Phi^*L\otimes \omega_\pi(\underline{\bf x}))^{\otimes N}\Big)\ :=\ 
  \coprod_{(C,{\bf x},\ph)\in \calc(M;p)} \Gamma(C;\Phi^*L^{\otimes N}\otimes
  \omega_\pi(\underline{\bf x})^{\otimes N})\,.
\]
Locally we may use Banach bundles of \v Cech-cocycles as in Section~\ref{cech}
to uniformize $\pi_*(\Phi^*L\otimes \omega_\pi(\underline{\bf x}))^{\otimes N}$
by the kernel of the \v Cech-coboundary operator $\check d$. But the cokernel of
$\check d$ is fiberwise isomorphic to $H^1(C,\ph^*L^{\otimes N}
\otimes \omega_C({\bf x})^{\otimes N})$, which vanishes by ampleness
for $N$ sufficiently large and as long as $(C,{\bf x},\ph)$ lies in a neighbourhood
$V$ of $K$. So over $V$, $F$ is an orbibundle of finite rank.

The lemma will follow once we produce an element $v$ of
the vector space $\Gamma (C_0,\varphi_0^*L^{\otimes N} \otimes
\omega_{C_0}({\bf x}_0)^{\otimes N})$ with $G$-orbit spanning
a $|G|$-dimensional linear subspace. In fact, an isomorphism
of $\rz[G]$ with $\langle G\cdot v\rangle$ is obtained by sending
${\bf 1}_e$ to $v$ and extending $G$-equivariantly.

The automorphism group $G$ acts by permutation on the
set of irreducible components of $C$. We decompose $C=D_1
\cup\ldots\cup D_m$ with each $D_i$ a union of irreducible
components in one $G$-orbit. Let $G_i\subset G$ be the subgroup
of automorphisms that restrict to the identity on the complement
of $D_i$. Then $G=G_1\times\ldots\times G_m$ and generic points
$y_i\in D_i$ have $G_i$-orbits (or, what is the same here,
$G$-orbits) of cardinality $|G_i|$. Put $A=G\cdot \{y_1,\ldots, y_m\}$.
By choice of $N$ there exists a section $v\in\Gamma
(C_0,\varphi_0^*L^{\otimes N} \otimes \omega_{C_0}({\bf x}_0)^{\otimes N})$
with $v(y_i)\neq 0$ for any $i$ but vanishing on the rest of $A$. 
Then $g=g_1 \cdot\ldots\cdot g_m$ with $g_i\in G_i$ produces a
section $g\cdot v$ whose non-zero values on $A$ are only
at $g_i\cdot y_i$. This shows that
\[
  \dim\,\langle G\cdot v \rangle\ =\ \prod_i  \dim\,\langle G_i\cdot v \rangle
  \ =\ \prod_i |G_i|\ =\ |G|\,,
\]
as desired.
\qed


\subsection{Kuranishi structure}\label{section_kuranishi}
We now turn to the construction of a Kuranishi structure (Definition~\ref
{Kuranishi_structure}) for the section $s_{\dbar,J}$. The morphism
$\tau:F\rightarrow E$ spanning the cokernel of the relative linearization
will only be defined over a neighbourhood of a compact subset $K\subset \calc(M;p)$,
that we henceforth assume chosen. In the application in Chapter~\ref{chapter_GW}
we will put $K=\calc^\hol_{R,g,k}(M,J)$, the pseudo-holomorphic curves of
genus $g$ with $k$ marked points and of homology class $R$.
\vspace{2ex}

\noindent
{\em Choice of local uniformizing systems.}\\
We now restrict attention to local uniformizing systems centered in pseudo-holomorphic
$(C,{\bf x},\ph)\in K$. This is enough to define the orbifold structure in a neighbourhood
of $K$ in view of Remark~\ref{construct_orbifold} and our discussion of change of
coordinates in Section~\ref{section_change_coord}.

The operator norm of $w\mapsto N_J(w,D\ph)$ as map from $\check L_1^p(\Delta;\ph^*
T_M)$ to $\check L^p(\Delta;\ph^*T_M\otimes\bar\Omega)$, $\Delta\subset C$, compared
to the operator norm of $\dbar_{\ph^*T_M,J}$, gets arbitrarily small with
$\diam(\Delta)$ (measured in any metric on $C$). From this
and the results of Chapter~\ref{chapter_dbar} one easily concludes the surjectivity
of the relative differential ${\hat\sigma}_J(C,{\bf x},\ph)=\dbar_{\ph^*T_M,J}
+N_J(\, .\,,D\ph)$ of the $\dbar_J$ operator at pseudo holomorphic
$(C,{\bf x},\ph)$, when restricted to sufficiently small open sets in $C$. The cokernel
of ${\hat\sigma}_J(C,{\bf x},\ph)$ can thus be spanned by 1-form valued sections
of $\ph^*T_M$ with (arbitrarily small) support outside some neighbourhood of
the set of special points of $(C,{\bf x},\ph)$. For unstable domains this argument
has to be supplemented by the observation that by the form of the rigidifying map
$\Lambda$, $D_{\bar V}\gamma(0)$ is surjective away from $\bigcup_i\supp\ph^*\rho_i$.

To make this uniform we now use $\eps$-thin parts for stable curves in $M$ as
developed in Section~\ref{section_intrinsic}. By compactness the number of
special points and bubble components is indeed uniformly bounded on $K$,
hence for all stable curves covered by local uniformizing systems centered
in some point of $K$. Combining this with the previous reasoning and uniform continuity
of ${\hat\sigma}_J$ at centers of local uniformizing systems, we infer the
existence of $\eps>0$ such that the cokernel of ${\hat\sigma}_J (C,{\bf x},\ph)$ can
be spanned by sections with support outside the $\eps$-thin part of $(C,{\bf x},\ph)$.

Since we may take the sets $U_i$ for $i>0$ in our definition of charts to be arbitrarily
small neighbourhoods of the singular set we may restrict to charts such that the
pull-back of $U_i$, $i>0$, to the universal curve $\Gamma$ is contained in the union of
the $\eps/2$-thin parts of $(\Gamma, \underline{\bf x}, \Phi)$. We may also assume
that the preimages of the $\eps/2$-thin part of the central fiber under
local admissable retractions $\kappa: \calc\rightarrow C$ used in the
construction of local charts are contained in the $\eps$-thin parts. Then for any of
these local uniformizing systems, centered in $(C,{\bf x},\ph)$ say, away from the
$\eps$-thin parts the trivialization $\Pi$ of ${q_1^p}_* \kappa^*\ph^*T_M$ is
given by $T^0_s\kappa_s^\star \dbar$. We use the same $\kappa$ for the
construction of local trivializations of the bundle $E$.
\vspace{2ex}

\noindent
{\em Construction of the morphism $\tau$.}\\
Let us now come to the construction of $\tau$. By uniform continuity of the relative
differential of $s_{\dbar, J}$ at the centers of local uniformising systems and
by compactness of $K$ it is enough to construct a morphism $\tau: F\rightarrow E$
with the required regularity properties and with image not contained in the image of
the relative differential $\sigma_J$ of $s_{\dbar,J}$ at the center $(C,{\bf x},\ph)$. Now
we assumed $(C,{\bf x},\ph)$ pseudo-holomorphic, so $\ph$ is in fact smooth.
Let $\alpha\in C^1(C;\ph^*T_M \otimes\bar\Omega) \setminus \img
\hat\sigma_J$ (continuously differentiable sections) have support away from the
$\eps$-thin part of $(C,{\bf x},\ph)$. Let $F$ be the finite rank bundle of
Lemma~\ref{bundle_F}, constructed for $(C,{\bf x},\ph)$, and put $G=
\Aut(C,{\bf x},\ph)$. We define $\tau$ locally as map from $F_0\times S\times
\bar V$ to $E$ by
\[
  \hat\tau:\ ({\bf 1}_\Psi, s,v)\ \longmapsto\ \eta(C_s,{\bf x}_s, \ph(s,v))
  \cdot \Pi^E_{s,v} (\Psi^{-1})^*\alpha
\]
on the copy of $\rz[G]\subset F_0$ and $0$ on the $G$-complementary subspace
of $F_0$. For $\eta$ we would like to take a bump function with support in $q(S\times
\bar V)$. Unfortunately, such bump functions, differentiable
relative $\calm_{g,k}$, do not seem to exist. Instead we will restrict to an open
suborbifold of $\calc(M;p)$ containing $K\cap \calc^\hol(M,J)$ that we are going to
specify. Let us first check for equivariance of $\hat\tau$: For $\Psi, \Psi'\in
\Aut(C,{\bf x},\ph)$ we obtain
\begin{eqnarray*}
  \hat\tau(\Psi({\bf 1}_{\Psi'},s,v))&=& \hat\tau({\bf 1}_{\Psi\Psi'},\Psi(s),\Psi(v))\ =\
  \eta(q(s,v)) \cdot \Pi^E_{\Psi(s),\Psi(v)} ({\Psi'}^{-1})^*(\Psi^{-1})^*\alpha\\
  &=&\eta(q(s,v)) \cdot (\Psi^{-1})^*\Pi^E_{s,v} ({\Psi'}^{-1})^*\alpha
  \ =\ (\Psi^{-1})^* \hat\tau({\bf 1}_{\Psi'},s,v)\,.
\end{eqnarray*}
For the definition of $\eta$ let $\chi\in C^\infty(\rz;[0,1])$ be a bump function,
say $\chi|_{[-1,1]}\equiv1$, $\supp\chi\subset[-2,2]$. For $\delta>0$ we set
$\chi_\delta(t)=\chi(\delta^{-1}t)$, which is a bump function with support
in $[-2\delta, 2\delta]$. Choose $G$-invariant scalar products on $S$ and
on $L^2(U_0\cap C;\ph^*T_M)$. The corresponding norms will be denoted
$|\,.\,|_S$ and $|\!|\,.\,|\!|_{U_0}$. Recall also the slice map $\gamma=\prj_S\times\nu:
S\times \bar V\hookrightarrow S\times V$ (Section~\ref{sect_rigidification}).
On $U:= q(S\times\bar V)$, $\eta$ is then defined by
\[
  \eta_\delta(C_s,{\bf x_s},\ph(s,v))\ :=\ \chi_\delta (|s|^2+|\!|\nu(s,v)|_{U_0}|\!|^2_{U_0})\,,
\]
with $\delta>0$ to be specified. By $G$-invariance of the norms this is a
well-defined function on the support of the local uniformizing system. Note that
$\eta$ does {\em not} have bounded support (not even in the $L^2$ topology)
because it disregards the behaviour of $v$ near the singular set of $C$. 
\begin{lemma}\label{support}
  For $\delta\ll1$: $\partial U\cap K\cap \supp \eta_\delta =\emptyset$.
\end{lemma}
\pf
Otherwise we find sequences $\delta_\nu\rightarrow0$ and $(C_\nu,{\bf x}_\nu,\ph_\nu)
\in \partial U\cap K\cap \supp \eta_{\delta_\nu}$. By compactness of $K$
we may assume a limit $(C_\infty,{\bf x}_\infty,\ph_\infty)\in \partial U\cap K$
to exist. By the definition of $\eta_\delta$, $(C_\infty,{\bf x}_\infty)=
(C,{\bf x})$, and $\ph_\infty=\ph$ away from the double points ($L^2$ convergence
suffices for that). Now we apply the identity theorem for pseudo holomorphic curves
to conclude $(C_\infty,{\bf x}_\infty,\ph_\infty)=(C,{\bf x},\ph)$, in contradiction to
$(C_\infty,{\bf x}_\infty,\ph_\infty)\in\partial U$.
\qed

We take a sufficiently small $\delta$ as in the lemma. Then by compactness of $K$
and since the $C^0$-toplogy is Hausdorff there exists a neighbourhood $X'$
of $K$ in the $C^0$-topology with the property $\partial U\cap X'
\cap \supp \eta_\delta =\emptyset$. We define $\eta$ on $X'$ as trivial extension
of $\eta_\delta$.
\vspace{2ex}

\noindent
{\em Regularity properties of $\tau$.}\\
Now let $(C',{\bf x}',\ph')$ be another pseudo-holomorphic curve and
$S'\times\bar V'$ one of our local uniformizing systems with center
$(C',{\bf x}',\ph')$ and structure map $q': S'\times\bar V' \rightarrow
\calc(M;p)$. We claim continuous differentiability of the uniformizer of $\tau$ in any
$(s',v')\in S'\times\bar V$. By the lemma, either $q'(s',v')$ is not contained
in the support of $\tau$, or there exists an open embedding from
a neighbourhood of $(s',v')$ into $S\times V$. Since both the slice map
$\gamma': S'\times\bar V' \hookrightarrow S'\times V'$ and the quotient
map $S\times V\rightarrow S\times \bar V$ are differentiable, it suffices
to discuss the unrigidified situation. With the notations and conventions
of Sections~\ref {sect_rigidification} and \ref{bundleE} (e.g.\ we identify
$\calc'\rightarrow S'$ with a locally analytic subset of $\calc\rightarrow S$ etc.),
the relevant term takes the form
\[
  {(\kappa'_s)^*}^{-1}
  \Gamma^\cz_{\ph'\circ\kappa'_s, \exp_{{T^0_s}'(\kappa'_s)^\star\dbar v}
  \ph'\circ\kappa'_s}
  \Gamma^\cz_{\exp_{T^0_s(\kappa_s^{-1})^\star\dbar\Theta
  {T^0_s}'(\kappa'_s)^\star\dbar v}
  \ph\circ\kappa_s,\ph\circ \kappa_s}\kappa_s^*\alpha_i\,,
\]
for $(s,v)$ in a neighbourhood of $(s',v')$ in $S'\times V'$.
Here we have used the fact that the trivializations $\Pi$ of ${q_1^p}_*\kappa^*
\ph^*T_M$ reduce to $T^0_s(\kappa_s^{-1})^\star\dbar$ away from the
$\eps$-thin part, on which
$\tau (x)$ is supported for any $x\in X'$ by construction. 
By differentiability of $\alpha=\sum \alpha_i d\bar z_i$ and of $\ph$ and $\ph'$
this expression depends manifestly differentiably on $(s,v)$.
\vspace{2ex}

\noindent
{\em Orientation and index.}\\
For pseudo-holomorphic curves $(C,{\bf x},\ph)$ with stable domains
the Riemann-Roch Theorem applied to $(\ph^*T_M, \dbar_{\ph^*T_M,J})$ yields
\[
  \frac{1}{2}\ind{\hat\sigma}_J(C,{\bf x},\ph)\ =\ 
  c_1(M,J)\cdot \ph_*[C]+(1-g)\cdot\dim_\cz M\,.
\]
The index (Definition~\ref{index_defi}) of our Kuranishi section at $(C,{\bf x},\ph)$
is this index plus the dimension $2(3g-3+k)$ of the base space of the universal
deformation of $(C,{\bf x})$, where $g$ is the genus of $C$ and $k$ the number of
marked points. If $(C,{\bf x})$ is unstable we have to restrict to a subspace
$\bar V\subset V$ of complex codimension equal to $\dim \Aut^0(C,{\bf x})$. But the
dimension of the base space $S$ of a semi-universal deformation of $(C,{\bf x})$
increases by the same number. So in any case, for any $(C,{\bf x},\ph)
\in \calc^\hol(M,J)\cap K$ the index of the Kuranishi section $s_{\dbar,J}$
depends only on $R=\ph_*[C]\in H_2(M;\gz)$, $g=g(c)$ and $k= \sharp{\bf x}$:
\[
  \ind_{(C,{\bf x},\ph)}s_{\dbar,J}\ =\ 2\Big(c_1(M,J)\cdot R+(1-g)\cdot\dim_\cz M
  +3g-3+k\Big)\,.
\]
Finally, complex linearity of the $\dbar$ operator and of $D_{\bar V}\gamma$
together with the complex structure on $S$ naturally orients the determinant line.
\vspace{2ex}

\begin{rem}\rm
With our choice of charts, $\hat s_{\dbar,J}$ is not differentiable in the horizontal
directions, the reason being singularities coming from the retraction $\kappa$
near the nodes. As pointed out to me by V.~Shevchishin one can remedy this
for charts centered in pseudo-holomorphic curves $(C,{\bf x},\ph)$ by
writing $\ph$ near the nodes as graph of $n$ complex valued functions $f_1,\ldots,
f_n$. For deformations supported near the nodes deform these functions as
before but using the alternative retraction from Remark~\ref{retraction},1. For
global deformations one can glue with deformations coming from exponential maps
and differentiable retractions on the nonsingular parts.

Charts of this form also show that in Theorem~\ref{orbifold} the assumption
of existence of an almost complex structure on $M$ can be removed.
\qed
\end{rem}


\section{Gromov-Witten invariants}
\label{chapter_GW}
In this chapter we consider closed symplectic manifolds $(M,\omega)$ with
an almost complex structure $J$ tamed by $\omega$. For $R\in H_2(M,\gz)$ and
natural numbers $g, k$ we restrict attention to the open suborbifold $\calc_{R,g,k}(M;p)
\subset \calc(M;p)$ of $k$-marked genus $g$ curves in $M$ representing the homology
class $R$. The subset $\calc_{R,g,k}^\hol(M,J)$ of $J$-holomorphic curves is then
compact (Theorem~\ref{cptness_thm}) and zero locus of the oriented
Kuranishi section $s_{\dbar,J}$ of $E$ of index
\[
  d(M,\omega,R,g)\ :=\ 2\Big(c_1(M,\omega)\cdot R+ (1-g)\cdot\dim_\cz M + 3g-3 +k\Big)\,.
\]
Note that $c_1(M,\omega)$ equals the first Chern class of the complex vector bundle
$(T_M,J)$ for any almost complex structure $J$ tamed by $\omega$.
The {\em virtual fundamental class} $[E, s_{\dbar,J}]$ of $\calc^\hol_{R,g,k}(M,J)$
obtained from Theorem~\ref{virfc} is thus an element in
$H_{d(M,\omega,R,g)}(\calc_{R,g,k}^\hol(M,J))$ that we denote $\GW^{M,J}_{R,g,k}$.
A priori there is an implicit dependence on
$p$ and the particular choice of Kuranishi structure, but
\begin{prop}
  $\GW^{M,J}_{R,g,k}$ is independent of $p$ and our choice of
  Kuranishi structure.
\end{prop}
\pf
We reexamine part of the proof of Theorem~\ref{virfc}. By construction of the finite rank
orbibundle $F$ we may assume that $F=F(p)$ can be constructed ``simultaneously''
for all $p$. This means that whenever $p<p'$ and $\iota_{pp'}:
\calc(M;p')\hookrightarrow\calc(M;p)$,
${\tilde\iota}_{pp'}:E(p')\hookrightarrow\iota_{pp'}^* E(p)$ are the natural inclusions,
then $F(p')=\iota_{pp'}^* F(p)$ and $\iota_{pp'}^*\tau(p)={\tilde\iota}_{pp'}\circ\tau(p')$.
Now $\tilde Z(p)\subset F(p)$ was defined as zero locus of $q^* s_{\dbar,J}+\tau(p)$, where
$q:F(p)\rightarrow\calc(M;p)$ is the projection. But from
\[
  \dbar_J\ph\ =\ -\tau(p)(C,{\bf x},\ph)\ \in\ L^\infty(C;\ph^*T_M\otimes{\bar\Omega}_C)
\]
and elliptic regularity we obtain $\ph\in\check L_1^p(C;M)$ for all $p$
(or even smoothness if the perturbation terms $\alpha$ were chosen smooth). So all
$\tilde Z(p)$ are identified via $\iota_{pp'}$ and ${\tilde\iota}_{pp'}$ and hence the
virtual fundamental classes all coincide.

As for independence of our choice of Kuranishi structure we just remark that local
uniformizers of $\tau$ are always differentiable, provided the sets $U_0$ are sufficiently
large. Thus for any two of such morphisms $\tau$, $\tau'$ we may find a system of
distinguished local trivializations making $\tau$ and $\tau'$ compatible Kuranishi
structures. 
\qed

There is a continuous map
\[
  \mbox{ev}:\calc_{R,g,k}^\hol(M,J)\longrightarrow M^k\,,\quad
  (C,{\bf x},\ph)\longmapsto(\ph(x_1),\ldots,\ph(x_k))\, .
\]
This is of course the restriction of a map with $\calc_{R,g,k}(M;p)$ as domain of
definition, but we will not need this. 
We are now ready to define GW-invariants, in the form proposed in \cite{kontmanin}.
Recall our convention $\calm_{g,k}=\{{\rm pt}\}$ if $2g+k<3$.
\begin{defi}\rm\label{GW_correspondence}
The maps
\begin{eqnarray*}
  GW_{R,g,k}^{M,J}:H^*(M)^{\otimes k}&\longrightarrow&H_*(\calm_{g,k}),\\
  \alpha_1\otimes\ldots\otimes\alpha_k&\longmapsto&p_*
  \Big(\GW^{M,J}_{R,g,k}\cap\mbox{ev}^*(\alpha_1\times\ldots \times\alpha_k)\Big)
\end{eqnarray*}
are called {\em GW-correspondences} (of $(M,J)$).
\qed
\end{defi}
\begin{rem}\rm
One can also incorporate other natural classes on $\calc^\hol_{R,g,k}(M,J)$, e.g.\ Chern
classes of the vector bundle formed by the tangent spaces at the $i$-th point,
cf.\ \cite{witten}, \cite{ruantian2}.

Moreover, it is not hard to see from our construction
that there is a $K$-theoretic index bundle $\ind s_{\dbar,J} =[T_{\tilde Z}]-[F]$, an element
in the Grothendieck group of orbibundles of finite rank over $\calc^\hol_{R,g,k}(M,J)$. Its
class is also independent of choices. If $s_{\dbar,J}$ is transverse, this is nothing but the
tangent bundle of the moduli space. So $\ind s_{\dbar, J}$ should be interpreted as
{\em virtual tangent bundle} of the moduli space. Characteristic classes of this bundle
can also be used in the construction of symplectic invariants.
\qed
\end{rem}
This still depends on $J$. However, the whole point of the theory is the large
invariance under changes of $J$, namely
\begin{theorem}\label{gw_are_invts}
  The GW-correspondences are invariants of the symplectic deformation type of $(M,\omega)$.
\end{theorem}
Recall that $(M,\omega)$, $(M',\omega')$ by definition are symplectic deformation equivalent
if there is a diffeomorphism $\Phi:M\rightarrow M'$ and a 1-parameter family of symplectic
forms $\omega_t$ on $M$, $t\in[0,1]$, with $\omega_0=\omega$, $\omega_1=\Phi^*\omega'$
($\omega_t$ may leave its cohomology class).
\vspace{1ex}

\noindent
\pf
We have to show independence of changes of $J$ inside the symplectic deformation class.
So let $\{J_t\}_{t\in(0,1)}$ be a continuous 1-parameter family of almost
complex structures on $M$ (of class $C^2$, say), each tamed
with respect to some symplectic form in the deformation class of $\omega$.
Then $\calc^\hol_{R,g,k}(M,J_t)$ is compact for any $t$. Let us
abbreviate $\calc:=\calc_{R,g,k}(M;p)$ and $\calc^\hol_t:=\calc^\hol_{R,g,k}(M,J_t)$. Now
consider the bundle $E\times(0,1)$ over $\calc\times(0,1)$ together with the section
$s_{\dbar,J_t}$ on the slice $\calc\times\{t\}$. For any $t_0\in(0,1)$ we may take
local uniformizing systems centered at $(C,{\bf x},\ph)\times\{t_0\}$ on
$\calc\times(0,1)$ using a fixed $J_{t_0}$, and local trivializations of $E\times(0,1)$
using the family $\{J_t\}$. One sees by reexamination of the formulas in the last chapter
that the $s_{\dbar,J_t}$ fit together to an oriented
Kuranishi section $s_{\dbar,\{J_t\}}$ of $E\times(0,1)$ (for the construction of
$\tau$ we might have to slightly shrink the intervall $(0,1)$). Theorem~\ref{virfc} now
shows that the various virtual fundamental classes $\GW_t\in
H_*(\calc_t^\hol)$ come from a single class $\GW\in
H_{*+1}(\calc^\hol)$, $\calc^\hol=\bigcup_t\calc_t^\hol\times\{t\}\subset\calc$, i.e.\
\[
  \GW_t= \GW\cap\prj^*\delta_t\,,
\]
if $\prj:\calc^\hol\rightarrow(0,1)$ is the projection and $\delta_t\in H_{\{t\}}^1((0,1))$ is
Poincar\' e dual to $[\{t\}]$. Letting $\tilde{\mbox{ev}}:\calc^\hol
\rightarrow M^k$ be the evaluation map and $\alpha_1,\ldots,\alpha_k\in H^*(M)$ we get
\[
  p_*\Big(\GW_t\cap\mbox{ev}^*(\alpha_1\times\ldots\times\alpha_k)\Big)
  =p_*\Big(\GW\cap(\prj^*\delta_t\cup\mbox{ev}^*
  (\alpha_1\times\ldots\times\alpha_k))\Big)\, .
\]
This is manifestly independent of $t\in(0,1)$.
\qed
\vspace{1ex}

The axioms of Behrend, Kontsevich, Manin \cite{kontmanin}, \cite{behrendmanin}
can be easily verified within our approach. We will not carry this out in this already
too long paper though. A sketch of the set-up and arguments needed is contained
in \cite{si2}.
\vspace{2ex}

\small
{\em Acknowledgements.} I thank my colleagues at Bochum
University, especially S.\ Schr\"oer for an open ear when needed most,
and H.\ Flenner, W.-D.\ Ruan, R.\ Ye, V.\ Shevchishin for occasional
discussions. Thanks also to R.\ Brussee for providing some detail to
\cite{brussee}. I am indepted to H.\ Hofer, M.\ Lehn, Y.\ Ruan and the referee
for comments and corrections concerning a previous version of this paper.
Special thanks to G.\ Tian for introducing me to the subject at an early
stage of the development of the field.


\vspace{2ex}

\noindent
{\small Fakult\"at f\"ur Mathematik\\
Ruhr-Universit\"at Bochum, D-44780 Bochum\\
e-mail: Bernd.Siebert@ruhr-uni-bochum.de} 

\addcontentsline{toc}{section}{References}
\begin{thebibliography}{XXXXXX}
\small\sloppy
\setlength{\itemsep}{0ex}\setlength{\baselineskip}{0.5ex}
\bibitem[Ab]{abikoff} W.\ Abikoff: {\sl The real analytic theory of
        Teichm\"uller space},
        Lecture Notes in Math.\ {\bf 820}, Springer 1976/1989
\bibitem[Be]{behrend} K.\ Behrend: {\em GW-invariants in algebraic
        geometry},
        Inv.\ Math.\ {\bf 127} (1997) 601--617
\bibitem[BeFa]{behrendfantechi} K.\ Behrend, B.\ Fantechi: {\em The
        intrinsic normal cone},
        Inv.\ Math.\ {\bf 128} (1997) 45--88
\bibitem[BeMa]{behrendmanin} K.\ Behrend, Y.\ Manin: {\em Stacks of stable
        maps and Gromov-Witten invariants},
        Duke.\ Math.\ Journ.\ {\bf 85} (1996) 1--60
\bibitem[Bi]{bingener} J.\ Bingener: {\em Offenheit der Versalit\"at in der analytischen
       Geometrie}, Math.\ Z.\ 173 (1980) 241--281
\bibitem[Br]{bredon} G.\ E.\ Bredon: {\sl Sheaf theory}, McGraw-Hill 1965
\bibitem[Bs]{brussee} R.\ Brussee: {\em The canonical class and the
        $C^\infty$-properties of K\"ahler surfaces},
        New York Journ.\ Math.\ {\bf 2} (1996) 103--146
        (available from {\tt http://nyjm.albany.edu/})
\bibitem[CfFe]{coifman} R.\ R.\ Coifman, C.\ Fefferman: {\em Weighted norm
        inequalities for maximal functions and singular integrals},
        Studia Math.\ {\bf 51} (1974) 241--250
\bibitem[Cn]{conway} J.\ Conway: {\sl Functions of one complex
        variable II}, Springer 1995
\bibitem[DeMu]{delignemumford} P.\ Deligne, D.\ Mumford: {\em The
        irreducibility of the space of curves of given genus},
        Publ.\ Math.\ IHES {\bf 36} (1996) 75--110
\bibitem[Do]{douady} A.\ Douady: {\em Le probl\'eme des modules locaux pour
        les espaces $\cz$-analytiques compacts},
        Ann.\ Sci.\ \'Ecole Norm.\ Sup.\ (4) {\bf 7} (1974) 569--602
\bibitem[Fe]{federer} H.\ Federer: {\sl Geometric measure theory},
       Springer 1969
\bibitem[Fl]{flenner} H.\ Flenner: {\em \"Uber Deformationen holomorpher
       Abbildungen"},
       Habilitationsschrift, Univ. Osnabr\"uck 1978
\bibitem[FoKn]{forsterknorr} O.\ Forster, K.\ Knorr: {\sl Konstruktion verseller
        Familien kompakter komplexer R\"aume},
        Lecture Notes in Math.\ {\bf 705},
        Springer 1979
\bibitem[FuOn]{fukaya} K.\ Fukaya, K.\ Ono: {\em Arnold conjecture
        and Gromov-Witten invariant}, 
        Warwick preprint 29/1996
\bibitem[Gr]{grauert} H.\ Grauert: {\em Der Satz von Kuranishi f\"ur kompakte
        komplexe R\"aume},
        Invent.\ Math.\ {\bf 25} (1974) 107--142
\bibitem[GrRe]{grauertremmert} H.\ Grauert, R.\ Remmert: {\sl Coherent
        analytic sheaves},
        Springer 1984
\bibitem[Gv]{gromov} M.\ Gromov: {\em Pseudo holomorphic curves
        in symplectic manifolds},
        Inv.\ Math.\ {\bf 82} (1985) 307--347
\bibitem[HoLiSk]{holizan} H.\ Hofer, V.\ Lizan, J.-C.\ Sikorav: {\em On
        genericity of holomorphic curves in 4-dimensional almost-complex
        manifolds},
        preprint Toulouse III 1994
\bibitem[IvSh]{ivshev} S.\ Ivashkovich, V.\ Shevchishin: {\em 
        Pseudo-holomorphic curves and envelopes of meromorphy of
        two-spheres in ${\bf CP}^2$},
        preprint Bochum 1995
\bibitem[Iv]{iversen} B.\ Iversen: {\sl Cohomology of Sheaves},
        Springer 1986
\bibitem[Kn]{knudson} F.\ Knudson: {\em The projectivity of the moduli
        space of stable curves, II: The stacks $M_{g,n}$},
        Math.\ Scand.\ {\bf 52} (1983) 161--199
\bibitem[KoMa]{kontmanin} M.\ Kontsevich, Y.\ Manin: {\em
        Gromov-Witten classes, quantum cohomology, and enumerative
        geometry},
        Comm.\ Math.\ Phys.\ {\bf 164} (1994) 525--562
\bibitem[KpKp]{kaupkaup} B.\ Kaup, L.\ Kaup: {\sl Holomorphic functions
        of several variables},
        de Gruyter 1983
\bibitem[Kw]{kawasaki} T.\ Kawasaki: {\em The signature theorem for
        V-manifolds},
        Topology {\bf 17} (1978) 75--83
\bibitem[Ku]{kuranishi} M.\ Kuranishi: {\em New proof for the existence
        of locally complete families of complex structures},
        in: {\sl Proceedings of the conference on complex analysis,
        Minneapolis 1964}, A.\ Aeppli, E.\ Calabi, H.\ R\"ohrl (eds.) 142--154,
        Springer 1965
\bibitem[Le]{lehto} O.\ Lehto: {\sl Univalent functions and Teichm\"uller
        spaces},
        Springer 1987
\bibitem[LiTi1]{litian} J.\ Li, G.\ Tian: {\em Virtual moduli cycles and
        GW-invariants of algebraic varieties},
        Journal Amer.\ Math.\ Soc.\ {\bf 11} (1998) 119--174
\bibitem[LiTi2]{litian2} J.\ Li, G.\ Tian: {\em Virtual moduli cycles and
        Gromov-Witten invariants of general symplectic manifolds},
        preprint {\tt alg-geom 9608032}.
\bibitem[LuTi]{liu} G.\ Liu, G.\ Tian: {\em Floer homology and Arnold
        conjecture}, preprint 1996
\bibitem[LcMc]{lockmac} R.\ Lockhardt, R.\ Mc Owen: {\em Elliptic
        differential operators on noncompact manifolds},
        Ann.\ Sc.\ Norm.\ Sup.\ Pisa {\bf 12} (1985) 409--447
\bibitem[Lj]{looij} E.\ Looijenga: {\em Cellular decompositions of
        compactified moduli spaces of pointed curves}, in:
        {\sl The moduli space of curves, Texel Island 1994}, 369--400,
        Progr.\ Math.\ {\bf 129}, Birkh\"auser
\bibitem[MaPl]{mazya} V.\ G.\ Maz'ya, B.\ A.\ Plamenevski\u{\i}:
        {\em Estimates in $L_p$ and in H\"older classes and the
        Miranda-Agmon maximum principle for solutions of elliptic
        boundary value problems in domains with singular points on the
        boundary},
        Math.\ Nachr.\ {\bf 81} (1978) 25--82 (in Russian),
        engl.\ transl.\ in: AMS Transl., Ser.\ 2 {\bf 123} (1984) 1--56
\bibitem[McSa]{duffsal} D.\ McDuff, D.\ Salamon: {\sl $J$-holomorphic
        curves and quantum cohomology}, Amer.\ Math.\ Soc.\ 1994
\bibitem[Mu]{mumford} D.\ Mumford: {\em Stability of projective varieties},
        L'Enseignement Mathématique {\bf 24}, Geneva 1977
\bibitem[Pa]{palamodov} V.\ P.\ Palamodov: {\em Deformations of complex
        spaces} (Russian),
        Uspekhi Mat.\ Nauk (3) {\bf 31} (1976) 129--194;
        English translation: Russian Math.\ Surveys (3) {\bf 31} (1976) 129--197
\bibitem[Pn]{pansu} P.\ Pansu: {\em Compactness}, in: {\sl Holomorphic
        curves in symplectic geometry},
        M.\ Audin, J.\ Lafontaine (eds.) 233--249, Birkh\"auser 1994
\bibitem[PrWo]{parkerwolf} T.\ Parker, J.\ Wolfson: {\em
        Pseudoholomorphic maps and bubble trees},
        Journ.\ Geom.\ Anal.\ {\bf 3} (1993) 63--98
\bibitem[Ra]{ran} Z.\ Ran: {\em Deformations of maps},
        in: {\sl Algebraic curves and projective geometry},
        E.\ Ballico, C.\ Ciliberto (eds.) 246--253, Lecture Notes in Math.\ 1389,
        Springer 1989
\bibitem[Ru1]{ruan} Y.\ Ruan: {\em Topological sigma model and
        Donaldson type invariants in Gromov theory},
        Duke Math.\ Journ. {\bf 83} (1996) 461--500
\bibitem[Ru2]{ruan2} Y.\ Ruan: {\em Virtual neighborhoods and
        pseudo-holomorphic curves},
        preprint {\tt alg-geom 9611021}
\bibitem[RuTi1]{ruantian1} Y.\ Ruan, G.\ Tian: {\em A mathematical
        theory of quantum cohomology},
        Journ.\ Diff.\ Geom.\ {\bf 42} (1995) 259--367
\bibitem[RuTi2]{ruantian2} Y.\ Ruan, G.\ Tian: {\em Higher genus
        symplectic invariants and sigma model coupled with gravity},
        Inv.\ Math.\ {\bf 130} (1997) 455--516
\bibitem[SaUh]{sacks} J.\ Sacks, K.\ Uhlenbeck: {\em The existence of minimal
        immersions of $2$-spheres}, 
        Ann.\ of Math.\ {\bf 113} (1981) 1--24
\bibitem[Sa]{satake} I.\ Satake: {\em The Gauss-Bonnet theorem
        for V-manifolds},
        Journ.\ Math.\ Soc.\ Japan {\bf 9} (1957) 464--492
\bibitem[Sc]{schlessinger} M.\ Schlessinger: {\em Functors of Artin rings},
        Trans.\ Amer.\ Math.\ Soc.\ {\bf 130} (1968) 208--222
\bibitem[SeSi]{sesinger} R.\ Seeley, I.\ Singer: {\em Extending $\dbar$ to
        singular Riemann surfaces},
        Journ.\ Geom.\ Phys.\ {\bf 5} (1988) 121--136
\bibitem[Si1]{si2} B.\ Siebert: {\em Symplectic Gromov-Witten invariants},
        to appear in: {\sl Proceedings of the AGE conference, Warwick 1996},
        K.\ Hulek, M.\ Reid (eds.), Cambridge Univ.\ Press
\bibitem[Si2]{si} B.\ Siebert: {\em Symplectic and algebraic Gromov-Witten
        invariants coincide}, in preparation
\bibitem[Sk]{sikorav} J.-C.\ Sikorav: {\em Some properties of holomorphic
         curves in almost complex manifolds},
         in: {\sl Holomorphic curves in symplectic geometry},
         M.\ Audin, J.\ Lafontaine (eds.) 165--189, Birkh\"auser 1994
\bibitem[Sl]{sklya} E.\ G.\ Sklyarenko: {\em Homology and cohomology
        theories of general spaces},
        in: {\sl General topology II}, A.\ V.\ Arhangel'skii (ed.),
        Encyclopedia of mathematical sciences, Springer 1996
\bibitem[Va]{vafa} C.\ Vafa: {\em Topological Mirrors and Quantum Rings},
        in: {\sl Essays on Mirror Manifolds}, S.-T.\ Yau (ed.) 96--119, International
        Press 1992
\bibitem[Ve]{vekua} I.\ N.\ Vekua: {\sl Generalized analytic functions},
        Pergamon 1962
\bibitem[Wh]{whitehead} G.\ W.\ Whitehead: {\em Homotopy theory},
        Springer 1978
\bibitem[Wi1]{witten0} E.\ Witten: {\em Topological sigma models},
        Comm.\ Math.\ Phys.\ {\bf 118} (1988) 411--449
\bibitem[Wi2]{witten} E.\ Witten: {\em Two-dimensional gravity and
        intersection theory on moduli space},
        in: Surveys in Differential Geometry {\bf 1} (1991) 243--310
\bibitem[Wo]{wolpert} S.\ Wolpert: {\em The hyperbolic metric and the
        geometry of the universal curve},
        Journ.\ Diff.\ Geom.\ {\bf 31} (1990) 417--472
\bibitem[Ye]{ye} R.\ Ye: {\em Gromov's compactness theorem for
        pseudoholomorphic curves},
        Trans.\ Amer.\ Math.\ Soc.\ {\bf 342} (1994) 671--694
\bibitem[Ze]{zeidler} E.\ Zeidler: {\sl Nonlinear functional analysis and
        its applications I},
        Springer 1986
\bibitem[Zi]{ziemer} W.\ Ziemer: {\sl Weakly differentiable functions},
        Springer 1989
\end{thebibliography}
\end{document}